\documentclass[11pt]{article}

\usepackage{graphicx}
\usepackage{caption}
\usepackage{subcaption}

\usepackage[margin = 1in]{geometry}
\usepackage{amsmath,amssymb,amsfonts,amsthm}
\usepackage{color}
\usepackage{hyperref} 
\def \p{\partial}

\def \ds{\displaystyle}
\def \div{\nabla \cdot}
\def \grad{\nabla}

\def \bM{\mathbb{M}}

\def \cTh{\mathcal{T}_h}
\def \cThtil{\mathcal{\tilde{T}}_h}
\def \cT{\mathcal{T}_{c}}

\def \cTch{\mathcal{T}_{c,h}}

\newcommand{\vect}[1]{\boldsymbol{#1}}
\newcommand{\trace}{\text{trace}}

\newcommand{\bey}{\begin{eqnarray}}
\newcommand{\eey}{\end{eqnarray}}

\newcommand{\beq}{\begin{equation}}
\newcommand{\eeq}{\end{equation}}
\theoremstyle{plain}

\theoremstyle{definition}

\theoremstyle{remark}
\newtheorem{exam}{\hspace{6mm}Example}[section]

\newenvironment{keywords}%
   {\begin{trivlist}\item[]{\bfseries\sffamily Key words:}~}
   {\end{trivlist}}
\newenvironment{AMS}%
   {\begin{trivlist}\item[]{\bfseries\sffamily AMS subject classifications:}~}
   {\end{trivlist}}

\allowdisplaybreaks

\title{%
Adaptive Finite Element Solution of the Porous Medium Equation in Pressure Formulation
}

\author{%
   Cuong~Ngo%
   \thanks{%
      Bethany College, 335 East Swensson Street, Lindsborg, KS~67456, U.S.A.
      (\href{mailto:ngoc@bethanylb.edu}{\nolinkurl{ngoc@bethanylb.edu}}).%
   }%
   \and
   Weizhang~Huang%
   \thanks{%
      Department of~Mathematics, University of~Kansas, Lawrence, KS~66045, U.S.A.
      (\href{mailto:whuang@ku.edu}{\nolinkurl{whuang@ku.edu}}).%
   }%
 }

\date{}

\begin{document}
\maketitle

\begin{abstract}

A lack of regularity in the solution of the porous medium equation
poses a serious challenge in its theoretical and numerical studies.
A common strategy in theoretical studies
is to utilize the pressure formulation of the equation where a new variable called
the mathematical pressure is introduced. It is known that the new variable has much better regularity than the original
one and Darcy's law for the movement of the free boundary can be expressed naturally
in this new variable. The pressure formulation has not been used in numerical studies.
The goal of this work is to study its use in the adaptive finite element
solution of the porous medium equation.
The MMPDE moving mesh strategy is employed for adaptive mesh movement while
linear finite elements are used for spatial discretization.
The free boundary is traced explicitly by integrating Darcy's law with the Euler scheme.
Numerical results are presented for three two-dimensional examples.
The method is shown to be second-order in space and first-order in time in the pressure variable.
Moreover, the convergence order of the error in the location of the free boundary is almost second-order
in the maximum norm. However, numerical results also show that the convergence order of the method
in the original variable stays between first-order and second-order in the $L^1$ norm or between
0.5th-order and first-order in the $L^2$ norm. Nevertheless, the current method can offer some advantages over
numerical methods based on the original formulation for situations with large exponents
or when a more accurate location of the free boundary is desired.
\end{abstract}

\begin{keywords}
porous medium equation, adaptive moving mesh method, MMPDE method, finite element method,
free boundary, Darcy's law, pressure formulation
\end{keywords}

\begin{AMS}
65M60, 65M50, 35Q35
\end{AMS}

\section{Introduction}
\label{SEC:intro}

We consider the numerical solution of the two-dimensional (2D) porous medium equation (PME)
\begin{equation}
u_t = \div(|u|^m \grad u),
\label{eqn:PME-U}
\end{equation}
where $m$ is a physical parameter whose value and physical meaning depend on a specific application.
PME is a nonlinear and nontrivial extension of the classic heat equation and arises in many
applications including the modeling of gas flow through porous media, boundary layer theory for fluid dynamics,
image processing, and particle physics. An important feature of PME is its degeneracy: the diffusion coefficient
vanishes wherever the solution becomes zero. This induces the finite propagation property in the sense that
the solution of PME is compactly supported at all time if it is so initially. The boundary of the compact support
becomes a free boundary, propagating at a finite speed governed by Darcy's law.
Moreover, the boundary is becoming smoother while the solution itself is losing its regularity
at the free boundary as time evolves. Mathematical studies for PME began in 1950s; for example, see
Ole{\u\i}nik et al. \cite{Oleinik1958}, Kala{\v{s}}nikov \cite{Kalashnikov1967}, and
Aronson \cite{Aronson1969}. Significant progress has been made since the late 1980s; e.g. see
Herrero et al. \cite{Herrero1987}, Shmarev \cite{Shmarev2005}, and
the monograph of V{\'a}zquez \cite{Vazquez2007}.

PME has been a challenge for numerical simulation due to its lack of regularity and free boundary.
The finite element method (FEM) has been used on fixed uniform or quasi-uniform meshes and
error estimates have been established in, e.g.,  Rose \cite{Rose1983},
Nochetto and Verdi \cite{Nochetto-1988}, Rulla and Walkington \cite{Rulla1996},
Ebmeyer \cite{Ebmeyer1998}, and Wei and Lefton \cite{Wei1999}.
Although applicable to two and higher dimensions, most of these results are not optimal.
In particular, they state that the convergence order for P1 finite elements is at most first-order and decreases
as the parameter $m$ increases.
On the other hand, due to the dynamic nature of free boundary and the low solution regularity
thereon, it is natural to employ dynamically adaptive meshes for the numerical solution of PME.
For this reason, there has been a trend in using adaptive moving mesh methods for PME
since 1990s.  For example, Budd et al. \cite{BCHR99} study the numerical simulation of self-similarity solutions of
one-dimensional PME using the Moving Mesh PDE (MMPDE) method \cite{HRR94a,HR11}
with a monitor function designed to preserve scaling invariance. 
Baines et al.  \cite{BHJ05,BHJ05a,BHJJ06} develop a moving mesh linear FEM for PME that
locally conserves the mass in the initial solution. For the situation with $m = 1$, their method
achieves an optimal convergence order for both the 1D and 2D cases.
However, it does not produce an optimal convergence order for $m > 1$.
An exception can be made in 1D with a special initial mesh,
which leads to second-order convergence for the case $m = 3$.
Unfortunately,  this initial mesh has not yet been attempted in 2D due to its high computational cost.
More recently, Duque et al. \cite{Duque2013,Duque2014,Duque2015} apply an MMPDE method to PME
having variable exponent (where $m$ depends on location and time) and with or without absorption terms
and show first-order convergence of the method for 2D problems. 
Ngo and Huang \cite{NH2016} study the MMPDE method with linear finite elements
for a general form of PME (with constant or variable exponents and with or without absorption terms)
in a large domain that contains the solution support for the whole time period under consideration.
The method is shown to be able to handle problems with complex, emerging, or splitting free boundaries
and be second-order in convergence in space when the mesh adaptation is based on the Hessian
of the approximate solution.
While FEM has been a common choice for PME, other methods have also been studied,
including the finite difference method by Socolovsky \cite{Socolovsky1984},
the local discontinuous Galerkin method by Zhang and Wu \cite{Zhang2009}, and
the spectral Galerkin method by Barrera \cite{Barrera2011}. 
Especially, the local discontinuous Galerkin method of \cite{Zhang2009} is very successful,
which gives high order of convergence within the solution support (away from the free boundary)
and completely eliminates unwanted oscillations at the free boundary,
although it has yet to be applied in 2D.

The form of (\ref{eqn:PME-U}), which will hereafter
be referred to as the original formulation or simply PME-U,
has been used in all of the existing numerical works so far. A drawback of this formulation
is that its solution $u$ has a steep or infinite slope at the free boundary.
On the other hand, it is a common practice to use the mathematical pressure $v := {u^m}/{m}$
in theoretical studies of PME
(e.g., see \cite{Angenent1988,Caffarelli1980,Caffarelli1987,Daskalopoulos1998a})
since $v$ has much better regularity than $u$.
For the Barenblatt-Pattle solution of PME, for example, $v$ is Lipschitz continuous in $\vect{x}$ and $t$ and
$\nabla v$ is bounded in the support of the solution (see \cite{Vazquez2007}).
Using the new variable $v$, we can rewrite (\ref{eqn:PME-U}) into
\begin{equation}
v_t = \grad \cdot(m v \grad v) - (m-1) |\grad v|^2 ,
\label{eqn:PME-V}
\end{equation}
which will be referred to as PME-V hereafter.
Compared to (\ref{eqn:PME-U}), (\ref{eqn:PME-V}) is no longer in the divergence form for $m \neq 1$.
However, its solution has higher regularity than that of PME-U, as mentioned earlier.
Moreover, as will be seen later, the equation governing the movement of the free boundary
can be expressed naturally in $v$.

    \begin{figure}
        \centering
        \begin{subfigure}[b]{0.4\linewidth}
                \includegraphics[width=\textwidth]{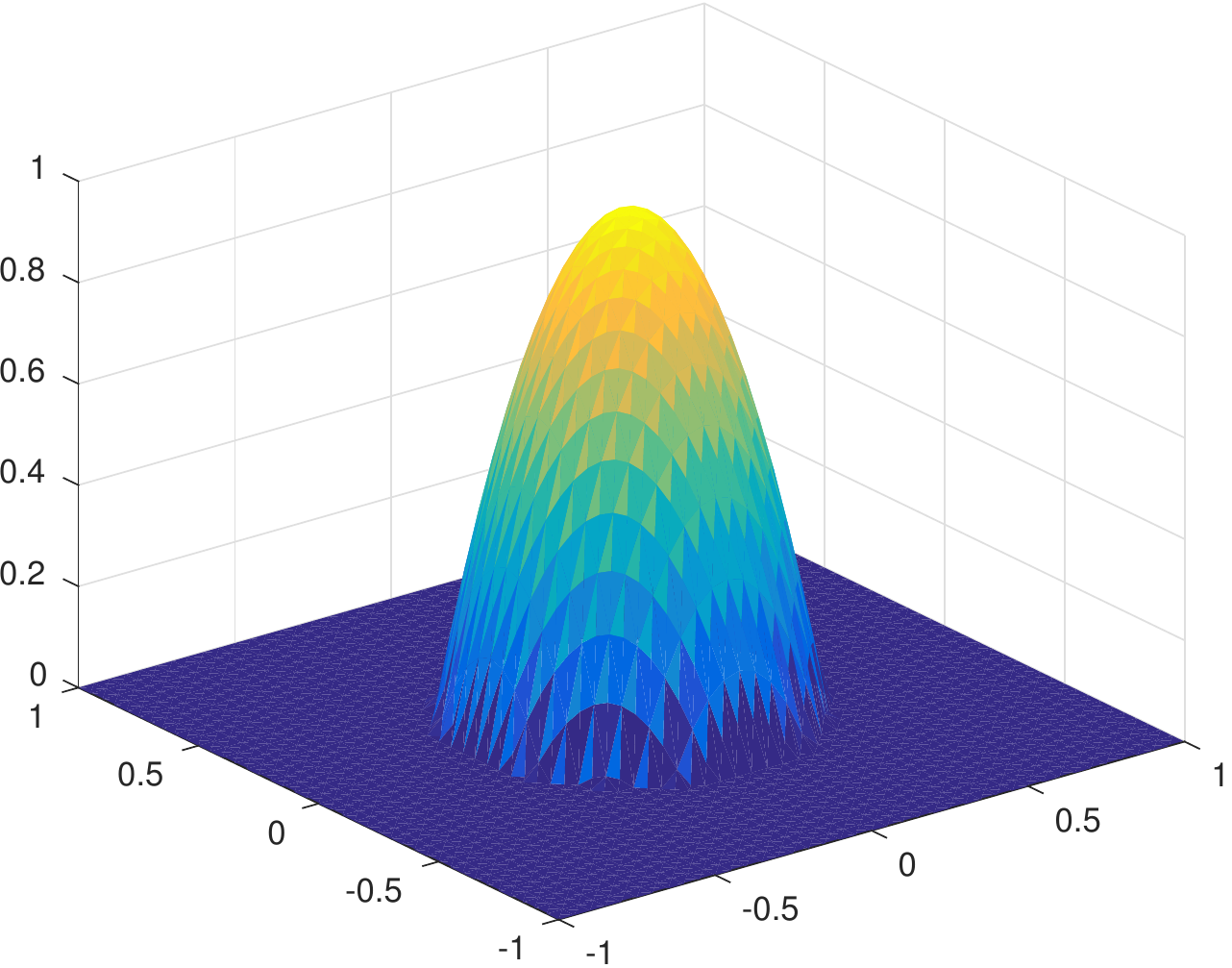}
                \caption{Embedding approach.}
                \label{subfig:pme2D-exact-m1-a}
        \end{subfigure}\hspace{5mm}%
        \begin{subfigure}[b]{0.4\linewidth}
                \includegraphics[width=\textwidth]{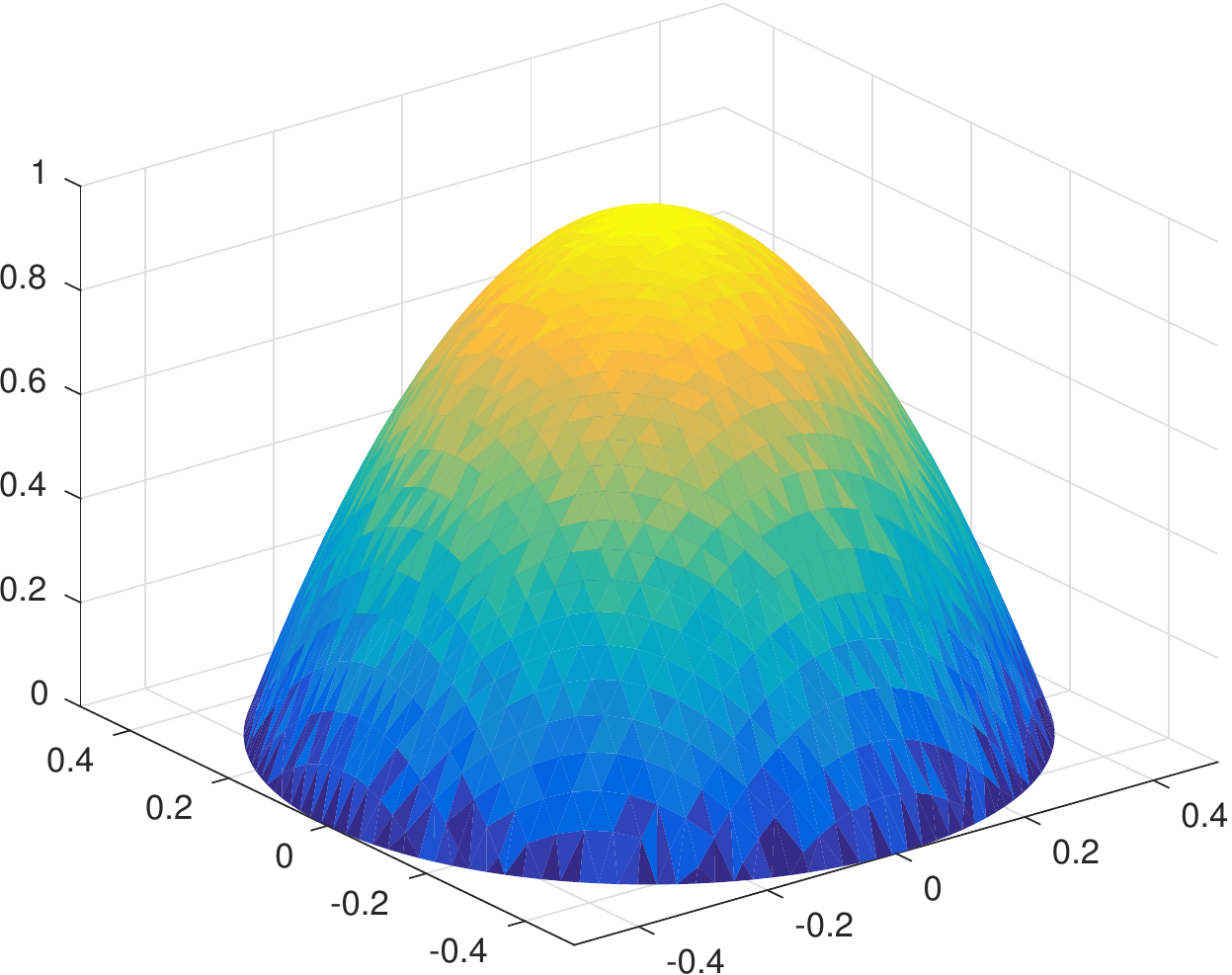}
                \caption{Nonembedding approach.}
                \label{subfig:pme2D-exact-m1-b}
        \end{subfigure}
        \caption{Two approaches for solving free boundary problems.}
    \end{figure}
    
When dealing with free or moving boundary problems, there are two main numerical approaches in general.
The first, which can be called the ``embedding'' or ``immersed boundary'' approach, performs computations on a fixed
but generally large domain which contains the support of the solution at all time of the simulation
(cf. Fig. \ref{subfig:pme2D-exact-m1-a}). Although more memory and more CPU time are required
for covering the extra domain outside the solution support and the solution usually has lower regularity
across the free boundary, this approach does not need to explicitly trace the free boundary
and is very robust for handling situations with complex free boundaries or topological changes (such as
merging or splitting) in free boundaries. This approach has been used in, for example, 
\cite{Duque2013,NH2016,Socolovsky1984, Zhang2009} for the numerical solution of PME.
The other approach, which can be called the ``nonembedding'' approach, solves PME
only within the solution support (see an illustration in Fig. \ref{subfig:pme2D-exact-m1-b}).
The main challenge of this approach is to trace the free boundary (denoted by $\Gamma(t)$)
via Darcy's law (e.g. see Shmarev \cite{Shmarev2005})
\begin{equation}
    \Gamma'(t) = - \lim_{\vect{x} \to \Gamma(t)^-} \grad \left (\frac{u^m}{m}\right ) \cdot {\vect{n}},
    \label{eqn:u-darcy}
\end{equation}
where $\Gamma'(t)$ is the velocity of the free boundary and ${\vect{n}}$ is its unit outward normal.
Note that the boundary is part of the solution to be sought. Moreover, once the boundary moves,
the mesh points should be redistributed or adjusted accordingly.
This approach has been used in \cite{BHJ05, BHJ05a,BHJJ06,BCHR99,Duque2014, Duque2015}
with moving mesh methods.

The objective of this paper is to study the numerical solution of PME in the PME-V formulation
using the nonembedding approach.
The PME-V formulation has not been used in the past in the numerical simulation of PME.
As mentioned earlier, its solution has better regularity. Moreover, noticing that the term in the parentheses
in (\ref{eqn:u-darcy}) is actually $v$, we can express the equation naturally in $v$ as
\begin{equation}
    \Gamma'(t) = - \lim_{\vect{x} \to \Gamma(t)^-} \grad v \cdot {\vect{n}} .
    \label{eqn:v-darcy}
\end{equation}
With this, the initial-boundary value problem of PME in the PME-V formulation reads as
\begin{equation}
  \begin{cases}
      v_t = \grad \cdot(m v \grad v) - (m-1) |\grad v|^2 \,, \quad & \text{in} \quad \Omega(t),\quad t \in (t_0,T] \\
      v(\vect{x},t_0) = v_0(\vect{x})  \,, \quad & \text{in} \quad \Omega(t_0) \\
      v(\vect{x},t) = 0 \,, \quad & \text{on} \quad   \p \Omega(t)\,, \quad t \in (t_0,T]  \\
      \Gamma'(t) = \ds -\lim_{\vect{x} \to \Gamma(t)^-}  \grad v(\vect{x},t) \cdot {\vect{n}}, & 
      \text{on} \quad   \p \Omega(t)\,, \quad t \in (t_0,T]
  \end{cases}
\label{eqn:PME-V-Darcy}
\end{equation}
where $\Omega(t)$ denotes the support of the solution at time $t$ and $\p \Omega(t) = \Gamma(t)$. 
We will explore the advantages and disadvantages of the PME-V formulation in
the numerical solution of PME and study how well the MMPDE moving mesh method \cite{HRR94a,HR11}
fares with this formulation. Recall that the latter has been used in \cite{NH2016} for the PME-U formulation
with the embedding approach. It has been shown that the MMPDE method is able to concentrate
the mesh points around the free boundary for the PME-U formulation and can lead to second-order
convergence for some cases with linear finite elements. We will investigate if this is still true for
the PME-V formulation.

The outline of the paper is as follows. Section~\ref{SEC:mmfem} describes the moving mesh finite element
method, involving the mesh generation via the MMPDE approach and
the linear finite element discretization of PME. Numerical results obtained for three 2D examples
are given in Section \ref{sec:pme-v-numerics}. Conclusions are drawn in
Section~\ref{SEC:conclusion}. The reader is referred to \cite{NH2016} for a summary of mathematical
properties of PME that are relevant to the numerical study or \cite{Vazquez2007} for a more complete
list of mathematical properties of PME.

\section{The moving mesh finite element method}
\label{SEC:mmfem}

In this section we describe the moving mesh FE solution of PME-V \eqref{eqn:PME-V-Darcy}
based on the MMPDE method. The discretization is based on the nonembedding approach
which requires to trace the free boundary explicitly at each time level using Darcy's law.

\subsection{The MMPDE method} 
\label{subsec:MMPDE}

We describe a dynamical mesh generation method based on the MMPDE method \cite{HRR94a}.
The computation is performed on the moving domain $\Omega(t)$ which is partitioned at each time $t$
by a moving mesh $\cTh(t)$, always having $N_v$ vertices (denoted by $\vect{x}_i(t) : i = 1, \ldots, N_v$),
$N$ elements, and a fixed connectivity. In addition, assuming that the first $N_{vi}$ vertices are interior vertices,
we can then denote $\Gamma_h(t) = \{ \vect{x}_i : i = N_{vi} + 1, \ldots, N_v \}$ as a discretization of the boundary $\Gamma(t)$. We further assume that the numerical solution is obtained at the time instants 
\begin{equation}
    t_0 < t_1 < \ldots < t_{n_f} \equiv T ,
    \label{time-grid}
\end{equation}
and denote related quantities at time level $t = t_n$ by $\Omega^n$, $\cTh^n$, $\Gamma_h^n$, and
$\{ \vect{x}_i^n : i = 1, \ldots, N_v\}$. Let the approximate solution at $t= t_n$ be $v_h^n$.
Our goal is to generate a new mesh $\cTh^{n+1}$ for the computation of the solution $v_h^{n+1}$
at the new time level $t = t_{n+1}$. To this end, we first move the free boundary according
to Darcy's law \eqref{eqn:v-darcy}. For simplicity, we use the Euler scheme to move boundary vertices
in the normal direction, i.e.,
\begin{equation}
\vect{x}_i^{n+1} = \vect{x}_i^n - \left. \left ((\nabla_h v_h^n ) \cdot {\vect{n}} \right ) {\vect{n}}
\right |_{\vect{x}_i^n} (t_{n+1}-t_n),
\quad i = N_{vi} + 1, \ldots, N_v
    \label{eqn:fw-euler}
\end{equation}
where $\nabla_h v_h^n$ is an average of the gradient of $v_h$ in the neighboring elements of $\vect{x}_i^n$.
After the boundary nodes are updated, we obtain a new intermediate physical mesh $\cThtil^{n+1}$.

We apply the MMPDE method to generate the new physical mesh $\cTh^{n+1}$
based on $\cThtil^{n+1}$ and $v_h^n$. The goal of the method is to make $\cTh^{n+1}$ to be as uniform
as possible when viewed in the Riemannian metric specified by a matrix-valued function $\bM$.
We call such a mesh an $\bM$-uniform mesh. The metric tensor $\bM$ can be chosen based on error estimates,
physical or geometric considerations, or combined. For our situation, we want more mesh points to concentrate
around the free boundary. Recall that the solution $v$ is smooth around the free boundary and
using the gradient or curvature of $v$ to define the metric tensor will not serve the purpose.
For this reason, we define 
\begin{equation}
    \mathbb{M} = \frac{1}{\sqrt{(v_h^n)^2 + 10^{-5}}} \; \mathbb{I} ,
    \label{eqn:M}
\end{equation}
with which more mesh points will be placed in regions where $v_h^n$ is small (i.e., regions around
the free boundary). The threshold $10^{-5}$ is chosen by trial and error. Its choice is
not very sensitive but cannot be too large (which will concentrate
less mesh points near the free boundary) or too small (which will concentrate too many mesh
points near the free boundary).

Having defined the metric tensor, we now discuss the MMPDE formulation.
For the moment, we denote the computational mesh by $\cT$ and the targeting physical mesh by $\cTh$,
which are assumed to have the same numbers of elements and vertices and the same
connectivity as $\cTh^{n}$.
It has been shown (e.g., see \cite{HR11}) that an $\bM$-uniform mesh satisfies the equidistribution and
alignment conditions
\begin{equation}
    |K| \sqrt{\det(\bM_K)} = \frac{\sigma_h |K_c|}{|\Omega_c|} , \quad \forall K \in \cTh
\label{eqn:equidist}
\end{equation}
\begin{equation}
\frac{1}{d} \trace 
      \left( 
	\left( F'_K \right)^{-1}
	  \bM_K^{-1} \left (F'_K \right )^{-T}
      \right)
      =
      \det \left( 
	\left( F'_K \right)^{-1}
	  \bM_K^{-1} \left (F'_K \right )^{-T}
          \right)^{\frac{1}{d}} , \quad \forall K \in \cTh
\label{eqn:alignment}
\end{equation}
where $d$ is the dimension of $\Omega$ ($d=2$ for the current situation),
$|K|$ and $|K_c|$ are the volumes of $K$ and its counterpart $K_c \in \cT$, respectively,
$\det(\cdot)$ and $\trace(\cdot)$ denote the determinant and trace of a matrix, respectively,
$\bM_K$ is the average of $\bM$ over $K$, $F'_K$ is the Jacobian matrix
of the affine mapping $F_K: K_c \to K$, and 
\[
    \sigma_h = \sum_{K \in \cTh} |K| \sqrt{\det(\bM_K)}, \quad 
|\Omega_c | = \sum_{K_c \in \cT} |K_c| .
\]
The equidistribution condition (\ref{eqn:equidist}) requires that all of the elements have the same
volume in the metric specified by $\bM$. It determines the size of the elements through $\bM$.
On the other hand, the alignment condition (\ref{eqn:alignment}) requires that $K$, when viewed
in the metric $\bM$, be similar to $K_c$ which is viewed in the Euclidean metric. It determines
the shape and orientation of $K$ through $\bM$ and the shape and orientation of $K_c$.

A mesh that closely satisfies the equidistribution and alignment condition can be obtained
by minimizing
\begin{align}
    I_h & = \theta \sum_{K \in \cTh} |K| \sqrt{\det(\bM_K)}
			\left( 
                        \trace({(F_K')}^{-1} {\bM}_{K}^{-1} {(F_K')}^{-T})
			\right)^{\frac{d p}{2}} \nonumber \\
			&  \quad \quad + (1-2\theta) d^{\frac{dp}{2}}
                        \sum_{K \in \cTh} |K| \sqrt{\det (\bM_K)}
			\left( 
			\frac{|K_c|}{|K| \sqrt{\det(\bM_K)}} 
                        \right)^{p},
\label{eqn:Ih}
 \end{align}
where $\theta \in (0, \frac{1}{2}]$ and $p > 1$ are non-dimensional parameters.
The function \eqref{eqn:Ih} is a discrete form of the continuous functional developed in \cite{Hua01b}
that combines the equidistribution and alignment conditions. Numerical experiment shows that
the choice $\theta = 1/3$ and $p = 2$ works for most problems, which we will use in our computation.
Due to its high nonlinearity, solving the minimization problem for \eqref{eqn:Ih} is a nontrivial task.
We use here the MMPDE method (a time-transient approach) to solve the problem.
Notice that $I_h$ is a function of the coordinates of the computational vertices ($\vect{\xi}_j,\; j = 1, ..., N_v$)
and the physical vertices ($\vect{x}_j,\; j = 1, ..., N_v$).
Moreover, recall that our goal is to generate a new physical mesh $\cTh^{n+1}$.
It is straightforward to solve for $\vect{x}_j,\; j = 1, ..., N_v$ directly
for a given reference computational mesh $\cT$.
The formulas for this so-called $\vect{x}$-formulation of the MMPDE method are given in
\cite[Equation (39)-(41)]{HK2014} but omitted here to save space. It is proven in \cite{HK2015}
that the mesh governed by this formulation stays nonsingular all the time if it is nonsingular initially.
This result holds at semi- and fully-discrete levels and for convex and concave domains in any dimension.
However, the metric tensor $\bM$ is defined on the physical domain and a function of $\vect{x}_j,\; j = 1, ..., N_v$,
which increases the nonlinearity of $I_h$. Even worse, $\bM$ is typically known only on the previous mesh
$\cTh^{n}$. As a result, $\bM$ needs to be updated constantly during the minimization process when
the physical mesh varies all the time.
To avoid these difficulties, we use here an indirect approach (the $\vect{\xi}$-formulation of the MMPDE method)
with which we take $\cTh = \cThtil^{n+1}$
(which is basically the physical mesh at the previous step and fixed during the minimization process),
minimize $I_h$ to obtain a new computation mesh, and finally obtain $\cTh^{n+1}$ through interpolation.
Notice that $I_h$ now is a function of $\vect{\xi}_j,\; j = 1, ..., N_v$ only.
Then the MMPDE for the computational vertices is defined as the gradient system of $I_h$ as 
\begin{equation}
\frac{d \vect{\xi}_j}{d t} = - \frac{P_j}{\tau} \left [ \frac{\partial I_h} {\partial \vect{\xi}_j} \right ]^T,
\quad j = 1, ..., N_v 
\label{eqn:mmpde-1}
\end{equation}
where the row vector ${\partial I_h}/{\partial \vect{\xi}_j}$ is the derivative of $I_h$ with respect to $\vect{\xi}_j$,
$\tau > 0$ is a parameter used to adjust how fast the mesh movement reacts to any change in the metric tensor,
and $P_j = \det(\bM(\vect{x}_j))^{\frac{p-1}{2}} $ is chosen such that \eqref{eqn:mmpde-1} is invariant
under the scaling transformation of $\bM$: $\bM \to c \, \bM$ for any positive constant $c$.
The derivative ${\partial I_h}/{\partial \vect{\xi}_j}$ can be found analytically using the notion of
scalar-by-matrix differentiation \cite{HK2014}. Using this, we can rewrite (\ref{eqn:mmpde-1}) into
\begin{equation}
\frac{d \vect{\xi}_j}{d t} = \frac{P_j}{\tau} \sum_{K \in \omega_j} |K| \vect{v}_{j_K}^K , \quad j = 1, \dotsc, N_v
\label{eqn:mmpde-2}
\end{equation}
where $\omega_j$ denotes the element patch associated with the vertex $\vect{x}_j$, $j_K$ is
the local index of the same vertex on the element $K$, and $\vect{v}_{j_K}^K$ is the velocity contributed
by $K$ to the vertex $\vect{x}_j$. Each element $K$ contributes a velocity $\vect{v}_i^K$
to its vertex $\vect{x}_i^K$ ($i = 0, \ldots, d$). These local velocities are given by 
\begin{equation}
\label{mmpde-3}
\begin{bmatrix} {(\vect{v}_1^K)}^T \\ \vdots \\ {(\vect{v}_d^K)}^T \end{bmatrix}
= - E_K^{-1} \frac{\partial G}{\partial \mathbb{J}} - \frac{\partial G}{\partial \det(\mathbb{J})}
\frac{\det(\hat{E}_K)}{\det(E_K)} \hat{E}_K^{-1}, 
\quad
\vect{v}_0^K = - \sum_{i=1}^d \vect{v}_d^K ,
\end{equation}
where 
$E_K = [\vect{x}_1^K-\vect{x}_0^K, ..., \vect{x}_d^K-\vect{x}_0^K]$ and
$\hat{E}_K = [{\vect{\xi}}_1^K-{\vect{\xi}}_0^K, ..., {\vect{\xi}}_d^K-\vect{\xi}_0^K]$ 
are the edge matrices of $K$ and $K_c$, respectively, $\mathbb{J} = (F_K')^{-1}$,
the function $G$ that is associated with \eqref{eqn:Ih} has the form 
\[
G(\mathbb{J}, \det(\mathbb{J}), \bM)
= \theta \sqrt{\det(\bM)} \left ( \trace(\mathbb{J} \bM^{-1} \mathbb{J}^T) \right )^{\frac{dp}{2}}
+ (1-2\theta) d^{\frac{d p}{2}} \sqrt{\det(\bM)} \left (\frac{\det(\mathbb{J})}{\sqrt{\det(\bM)}} \right )^p,
\]
and the derivatives $\frac{\partial G}{\partial \mathbb{J}}$ and $\frac{\partial G}{\partial \det(\mathbb{J})}$  are given by
\begin{align*}
& \frac{\partial G}{\partial \mathbb{J}} = d p \theta \sqrt{\det(\bM)} \left ( \trace(\mathbb{J} \bM^{-1} \mathbb{J}^T)
\right )^{\frac{dp}{2}-1}\bM^{-1} \mathbb{J}^T,
\\
& \frac{\partial G}{\partial \det(\mathbb{J})} = p (1-2\theta) d^{\frac{d p}{2}} \det(\bM)^{\frac{1-p}{2}}
\det (\mathbb{J})^{p-1} .
\end{align*} 
Note that these derivatives are evaluated at
$(\mathbb{J}, \det(\mathbb{J}), \bM) = ((F_K')^{-1}, \det(F_K')^{-1}, \bM_K)$.  

The above mesh velocities should be modified for boundary vertices. For instance, they should be set to be zero
for fixed boundary vertices. If boundary vertices are allowed to slid on the boundary, 
additional constraints should be applied so that they do not move inside or outside the domain.

The system \eqref{eqn:mmpde-2}, with proper modifications for the boundary vertices, can
be integrated from $t = t_n$ to $t = t_{n+1}$ to get the new computational mesh $\cT^{n+1}$,
with the initial mesh at $t = t_n$ taken to be the reference computational mesh $\cTch$
which can be chosen as a quasi-uniform mesh for $\Omega(t_0)$.
We use Matlab's ODE solver \texttt{ode15s} for this purpose. 
Since $\cT^{n+1}$ has the same structure (i.e., the same number of elements and vertices and
the same connectivity) as the physical mesh $\cThtil^{n+1}$, there is a piecewise affine map $\Phi_h^{n+1}$
such that $\cThtil^{n+1} = \Phi_h^{n+1}( \cT^{n+1})$. From this, we can compute the new physical
mesh as $\cTh^{n+1} = \Phi_h^{n+1}(\cTch)$ using linear interpolation.

\subsection{Linear finite element discretization on moving meshes} 
\label{subsec:LFEM}

In this subsection we describe the moving mesh FE method and the procedure for solving (\ref{eqn:PME-V-Darcy})
on moving meshes. We first notice that we now have the old mesh $\cTh^{n}$, the new mesh $\cTh^{n+1}$,
and the physical solution $v_h^n$ defined on $\cTh^{n}$.
For any $t \in [t_n, t_{n+1}]$, we define the mesh $\cTh(t)$ as the one having the same structure as
$\cTh^{n}$ and $\cTh^{n+1}$ and with the coordinates of the vertices as
\[
 \vect{x}_j(t) =
    \frac{t-t_n}{t_{n+1}-t_n} \vect{x}_j^{n+1} +
    \frac{t_{n+1}-t}{t_{n+1}-t_n} \vect{x}_j^{n} , \quad j = 1, ..., N_v.
\]
The corresponding nodal mesh velocities are given by
\[
\dot{\vect{x}}_j (t) = \frac{\vect{x}_j^{n+1}-\vect{x}_j^{n}}{t_{n+1}-t_n},\quad j = 1, ..., N_v.
\]
The linear basis function $\phi_j(\vect{x},t)$ associated with vertex $\vect{x}_j$ satisfies
\[ 
\phi_j|_K \in \mathbb{P}_1 \quad \forall K \in \cTh(t),
\quad \text{and} \quad \phi_j(\vect{x}_i,t) = \delta_{ij}, \quad i = 1, \ldots, N_v .
\]
We can then define the moving finite element space as
\[
    V_h(t) = \text{span}\{\phi_1(\cdot, t), ..., \phi_{N_{vi}}(\cdot, t)\}\,, \quad t \in (t_n,t_{n+1}) 
\]
and the linear finite approximation to (\ref{eqn:PME-V-Darcy}) as $v_h(\cdot, t) \in V_h(t), \, t \in (t_n, t_{n+1}]$
such that
\begin{equation}
  \int_{\Omega(t)} \frac{\p v_h}{\p t} \psi~d\vect{x} = -m \int_{\Omega(t)} v_h \grad v_h
        \cdot \grad \psi~d\vect{x}  + (1-m) \int_{\Omega(t)} |\grad v_h|^2 \psi~d\vect{x},
        \quad \forall \psi \in V_h(t) .
    \label{eqn:fem0}
\end{equation}
We need to pay special attention to the time derivative $\frac{\p v_h}{\p t}$ in the above equation.
Expressing $v_h$ into
\begin{equation}
    v_h(\cdot,t) = \sum_{j = 1}^{N_{vi}} v_j(t) \phi_j(\cdot, t)
    \label{eqn:linear-combination}
\end{equation}
and differentiating it with respect to $t$, we have
\[
    \frac{\p v_h}{\p t}  = \sum_{j=1}^{N_{v_i}} \frac{d v_j}{d t} \phi_j(\vect{x},t) + 
	\sum_{j=1}^{N_{v_i}} v_j(t) \frac{\p \phi_j}{\p t} .
\]
It can be shown (e.g., see Jimack and Wathen \cite[Lemma 2.3]{Jimack-1991}) that
\[
\frac{\p \phi_j}{\p t} = - \nabla \phi_j \cdot  \dot{\vect{X}} ,
\]
where $\dot{\vect{X}}$ is a piecewise linear mesh velocity function defined as
$\dot{\vect{X}} (\vect{x},t) = \sum_{j=1}^{N_v} \dot{\vect{x}}_j(t) \phi_j (\vect{x}, t)$.
Combining these we get
\begin{equation}
    \frac{\p v_h}{\p t}  = \sum_{j=1}^{N_{v_i}} \frac{d v_j}{d t} \phi_j(\vect{x},t) - \nabla v_h\cdot  \dot{\vect{X}} .
    \label{eqn:vh_t}
\end{equation}

 Substituting (\ref{eqn:linear-combination}) and (\ref{eqn:vh_t}) into \eqref{eqn:fem0}
 and taking $\psi = \phi_i,\, i = 1, \ldots, N_{vi}$ successively, we get 
\begin{align*}
    \sum_{j=1}^{N_{vi}} \left(\int_{\Omega(t)} 
     \phi_j \phi_i~d\vect{x}\right) \frac{d v_j}{dt} 
     = & \int_{\Omega(t)} \grad v_h \cdot \left( \dot{\vect{X}} \phi_i
    - m v_h \grad \phi_i \right) ~d\vect{x} \\
    & + (1-m) \int_{\Omega(t)} |\grad v_h|^2 \phi_i~d\vect{x}, \quad i = 1, ..., N_{vi}
\end{align*}
which can be written into the matrix form as
\begin{equation}
B(\vect{X}) \dot{\vect{V}} = F(\vect{V}, \vect{X}, \dot{\vect{X}}),
\label{eqn:fem2}
\end{equation}
where $B$ is the mass matrix and $\vect{X}$ and $\vect{V}$ are vectors representing the mesh and solution,
respectively. The ODE system \eqref{eqn:fem2} is integrated from $t_n$ to $t_{n+1}$
using the fifth-order Radau IIA method (an implicit Runge-Kutta method) with a standard time step selection
procedure where the relative and absolute tolerance are chosen as $10^{-6}$ and $10^{-8}$, respectively, and
the error estimation is based on a two-step error estimator of Gonzalez-Pinto et al.~\cite{Montijano2004}.

\vspace{10pt}

To conclude this section, we summarize the entire procedure for solving (\ref{eqn:PME-V-Darcy}) in the following.
    \begin{enumerate}
        \item[(i)] Partition and approximate $\Omega(t_0)$ with the initial mesh $\cTh^0$. We choose
        $\cTch = \cTh^0$ as our reference computational mesh for the MMPDE moving mesh method.
        Assume that at time level $t = t_n$, we have the solution $v_h^n$ on the mesh $\cTh^n$.
        \item[(ii)] Use the scheme \eqref{eqn:fw-euler} for Darcy's law to the boundary $\Gamma^n$
        	to obtain the new boundary $\Gamma^{n+1}$. This new boundary reflects the new domain
	$\Omega_h^{n+1}$ and is incorporated into the physical mesh $\cThtil^{n+1}$.
        \item[(iii)] Apply the MMPDE method with the mesh $\cThtil^{n+1}$ and its corresponding solution $v_h^n$
        to get the new physical mesh $\cTh^{n+1}$ at time level $t = t_{n+1}$.
        \item[(iv)] Apply the moving mesh FE method with $v_h^n$, $\cTh^n$, and $\cTh^{n+1}$ to
        obtain the new solution $v_h^{n+1}$ at time level $t = t_{n+1}$.
        \item[(v)] Repeat steps (ii) to (iv) until the final time $t = T$ is reached.
    \end{enumerate}

One can see that the boundary movement, the mesh movement, and the update of the physical solution
are split and performed sequentially. We thus expect that the method is first-order in time. On the other hand,
the physical PDE is discretized with linear finite elements and we expect the method to be second-order
in space when the solution is sufficiently smooth.
   
\section{Numerical results}
\label{sec:pme-v-numerics}

In this section we present numerical results obtained with the moving mesh method
described in the previous section for IBVP \eqref{eqn:PME-V-Darcy}. 
Unless otherwise stated, we choose the parameter $\tau$ of the MMPDE equation \eqref{eqn:mmpde-1}
to depend on the number of mesh elements $N$, i.e., 
\begin{equation}
    \tau = \min\left\{ 10^{-3}, \frac{10^{-1}}{N} \right\} .
    \label{eqn:pme-V-tau}
\end{equation}
Moreover, we restrict the time step to be no greater than $10^{-4}$, namely,
\begin{equation}
    \Delta t_{max} := \max_{n = 1,\,\ldots\,,\, n_f} (t_{n} - t_{n-1}) \leq 10^{-4}.
    \label{eqn:def-deltaT-max}
\end{equation}
These choices are made by trial and error. They produce satisfactory results for the problems we
have tested.

\begin{exam}[Barenblatt-Pattle solution] \label{ex:pme-V-BP}
In this example we consider the well-known Barenblatt-Pattle (BP) solution for PME-V given by
\begin{equation}
    v(r,t) = 
        \begin{cases}
            \frac{1}{m \lambda^{d m}(t)} 
            \left( 
            1 - \left( \frac{r}{r_0 \lambda(t)} \right)^2 
            \right) , \quad & \text{for} \quad |r| \leq r_0 \lambda(t)  \\
                                0 , \quad & \text{for} \quad |r| > r_0 \lambda(t)  
        \end{cases}
    \label{eqn:pme-V-BP}
\end{equation}
where
    \begin{equation}
        r = | \vect{x} | , \quad 
        \lambda(t) = \left(\frac{t}{t_0} \right)^{\frac{1}{2+d m}} ,
        \quad 
        t_0 = \frac{r_0^2 m}{2 (2 + d m)} \,,
        \label{eqn:pme-V-BP-par}
    \end{equation}
and $r_0>0$ is a parameter representing the radius of the solution support at the initial time $t_0$.
This solution is radially symmetric, self-similar, and compactly supported for any finite time.
It has been used as a benchmark example to test numerical algorithms for PME.
Notice that the BP solution (in terms of $v$) is smooth on its support (including the regions near the boundary)
for any parameter $m$, as can be seen in Fig. \ref{fig:BP-pressure}.

    \begin{figure}
        \centering
        \begin{subfigure}[b]{0.3\linewidth}\includegraphics[width=\textwidth]{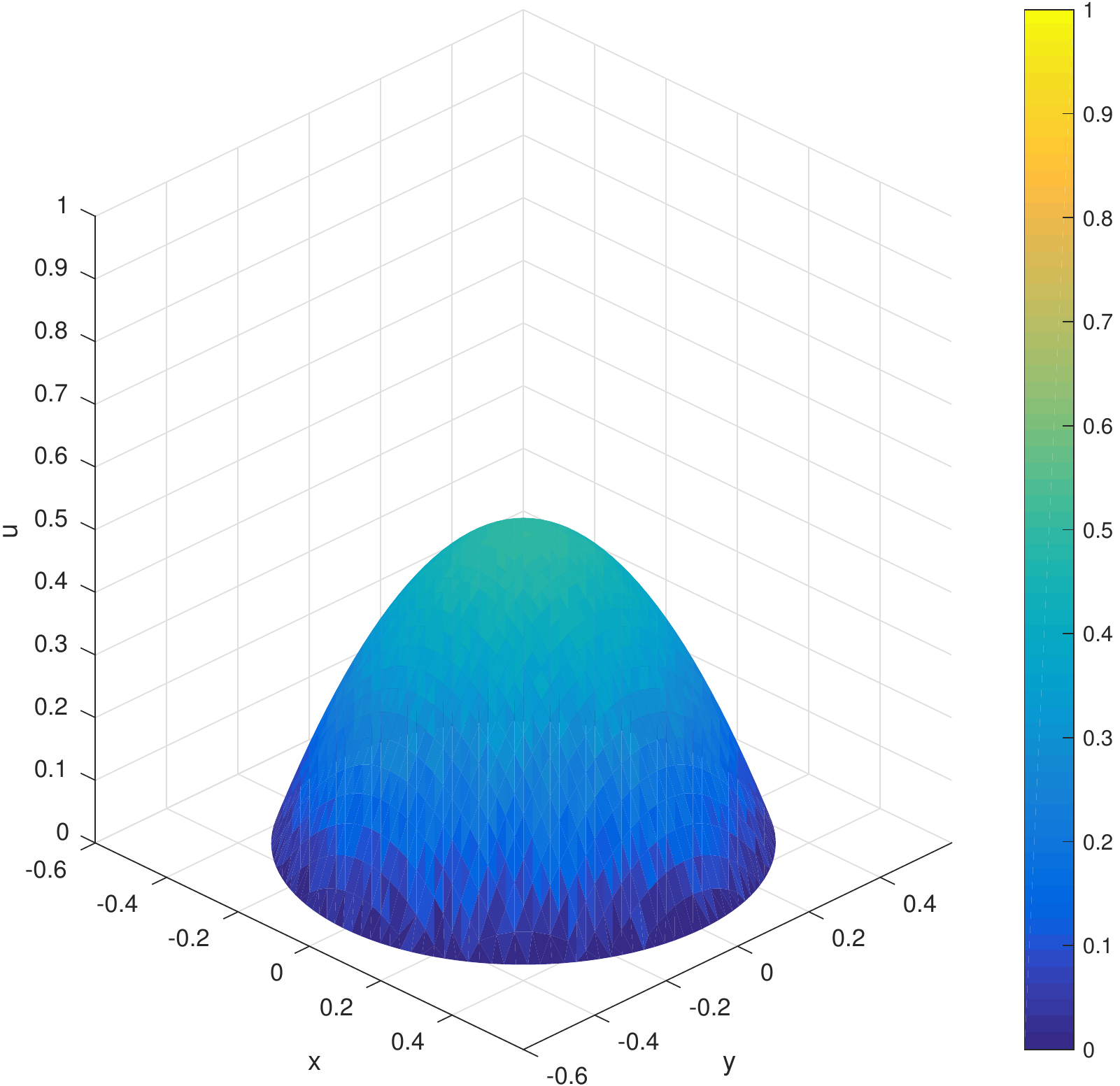}\caption{$m = 2$}\end{subfigure}%
        \begin{subfigure}[b]{0.3\linewidth}\includegraphics[width=\textwidth]{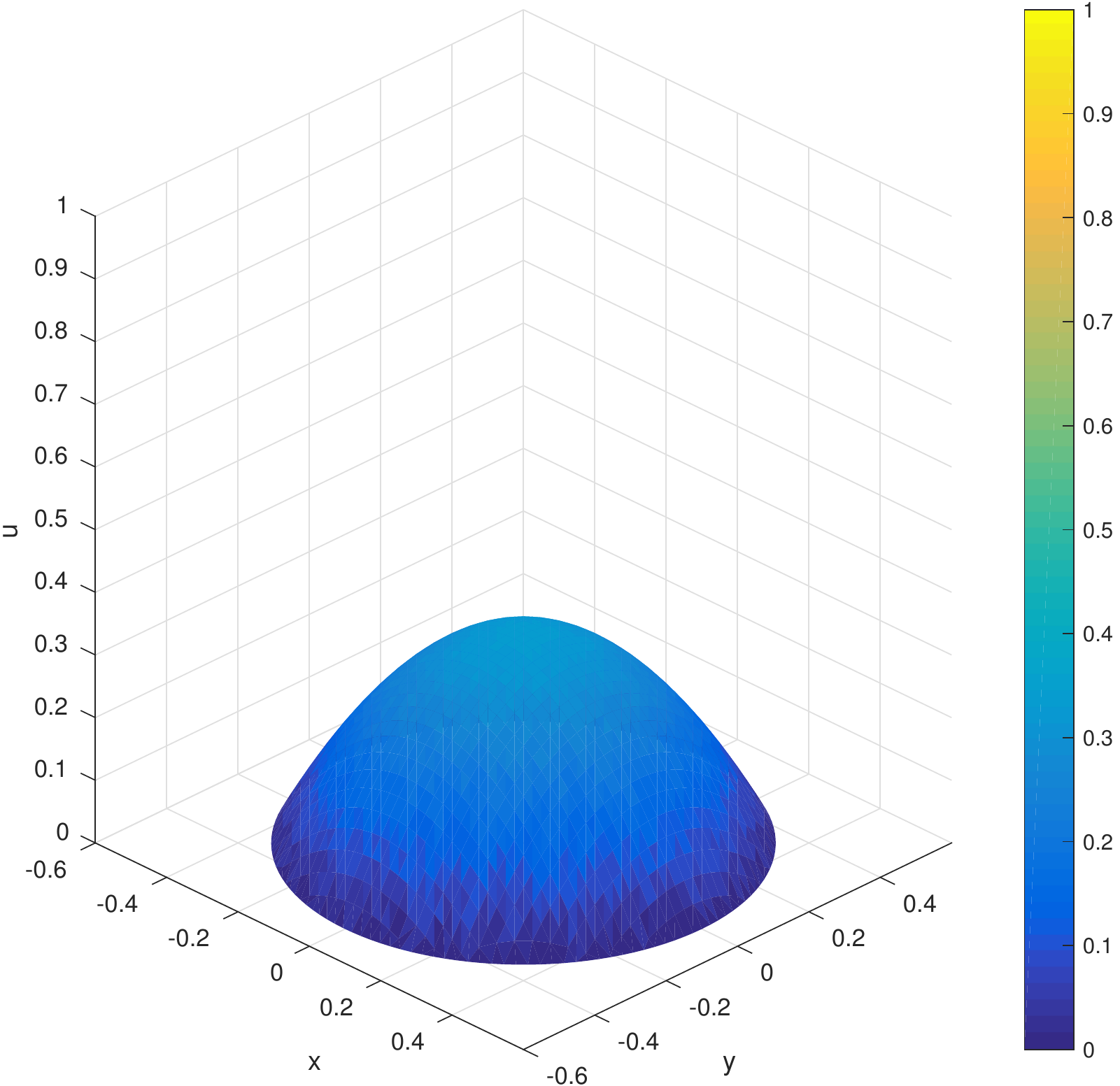}\caption{$m = 3$}\end{subfigure}%
        \begin{subfigure}[b]{0.3\linewidth}\includegraphics[width=\textwidth]{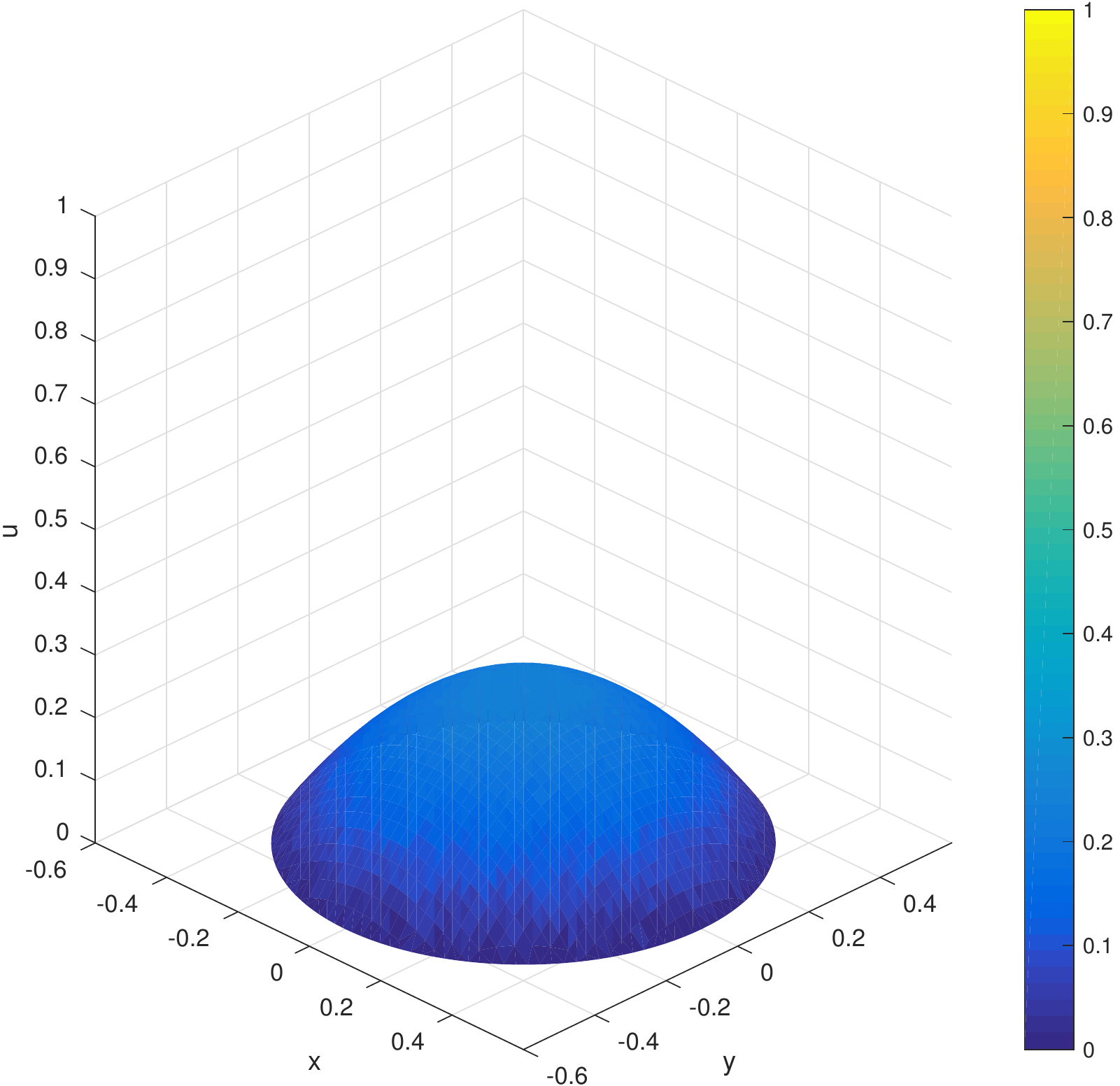}\caption{$m = 4$}\end{subfigure}%
        \caption{The Barenblatt-Pattle special solution for PME-V \eqref{eqn:PME-V-Darcy} is plotted for
        $m=2$, $3$, and $4$.}
        \label{fig:BP-pressure}
    \end{figure}
    
The moving mesh method is applied with $T = (t_0 + 0.1)/2$ and $r_0 = 0.5$
(with $t_0$ defined in \eqref{eqn:pme-V-BP-par}). 
A typical mesh and its associated solution are plotted in Fig. \ref{fig:pmeV-m3}.
One can see that the mesh points are concentrated correctly around the free boundary.
This is a demonstration of the mesh adaptation capability of the numerical method.
A convergence history (as the number of mesh elements $N$ increases)
for both the solution and the boundary is given in Fig. \ref{fig:pme-V-conv} for $m = 2,\, 3,\, 4$, and $5$.
The $L^2$-norm of the error converges with a second-order rate (in terms of the average element diameter
$h = 1/\sqrt{N}$) for all considered values of $m$. (A similar convergence behavior is also observed in the $L^1$-norm.)
The convergence order of the error in the boundary location, measured in the $L^\infty$-norm,
is almost second-order.
One may notice that these errors are smaller for higher $m$.
A closer examination shows that there are no oscillations in the computed solutions.
To see how accurate the computed solution in the original variable is, we obtain $u$ using the computed $v$ and
compare it with the exact solution. The convergence history in $u$ is shown in
Fig. \ref{fig:pme-UV-conv} for both the $L^1$ and $L^2$-norms. 
It can be seen that the convergence order in this variable in $L^2$-norm is about first-order for $m=2$ and
deteriorates as $m$ increases, which is indeed consistent with theoretical estimates in literature (e.g. \cite{Ebmeyer1998,Nochetto-1988,Rose1983,Rulla1996,Wei1999}).
Interestingly, the error in the $L^1$-norm has a better convergence rate, almost second-order for $m=2$
and 1.5th-order for $m = 5$.

We now study the effects of large $m$ on the computation. Recall that the solution $u$ of PME-U formulation
has a steeper or infinite slope at the free boundary as $m$ increases. Due to this, numerical methods
based on PME-U generally have difficulties to deal with large $m$ situations. For example,
the MMPDE moving mesh FE method in \cite{NH2016} can lead to second-order convergence
in the $L^2$ norm of error in $u$ for small $m$ ($m < 5$) but has difficulty to maintain the same
convergence order for large $m$ ($m \ge 5$). This is because, in the latter situation, the mesh elements
have to be extremely stretched near the free boundary which can cause the time integration of PME to stop
due to too small time steps.
On the other hand, $m$ has nearly no effects on the computation of PME-V and the convergence order
in $v$ stays almost the same (second-order) for all $m$; see Fig.~\ref{fig:pme-V-conv}(a) for $m = 2,\, 3,\, 4$, and $5$
and Fig.~\ref{fig:pme-V-conv-Bigm}(a) for much higher $m=8,\; 9, \; 10$, and $15$.
The convergence order of the error in the boundary location stays almost second-order too for large $m$
(cf. Fig.~\ref{fig:pme-V-conv-Bigm}(b)).
However, the convergence order in the $u$ computed through $v$ becomes close to 0.5th-order in $L^2$
norm and first-order in $L^1$ norm for very large $m$; see Fig.~\ref{fig:pme-UV-conv-Bigm}(a,b).

It is interesting to see how the time step size affects the accuracy of the computed solution.
We use a relatively fine mesh with $N = 40459$ and a sequence of $\Delta t_{max}$. The error is shown
in Fig.~\ref{fig:pme-V-dtFixed}. One can see that the error remains constant for small $\Delta t_{max}$
because the total error is dominated by the spatial discretization error when the time step is small.
On the other hand, for larger $\Delta t_{max}$, the error decreases at a rate between first-order and second-order.
The overall convergence is more like first-order, which is consistent with our expectation
from the construction of the numerical method.
\qed
\end{exam}

\begin{figure}[ht]
    \begin{center}
    \begin{subfigure}[b]{0.35\linewidth}
       \includegraphics[scale=0.34]{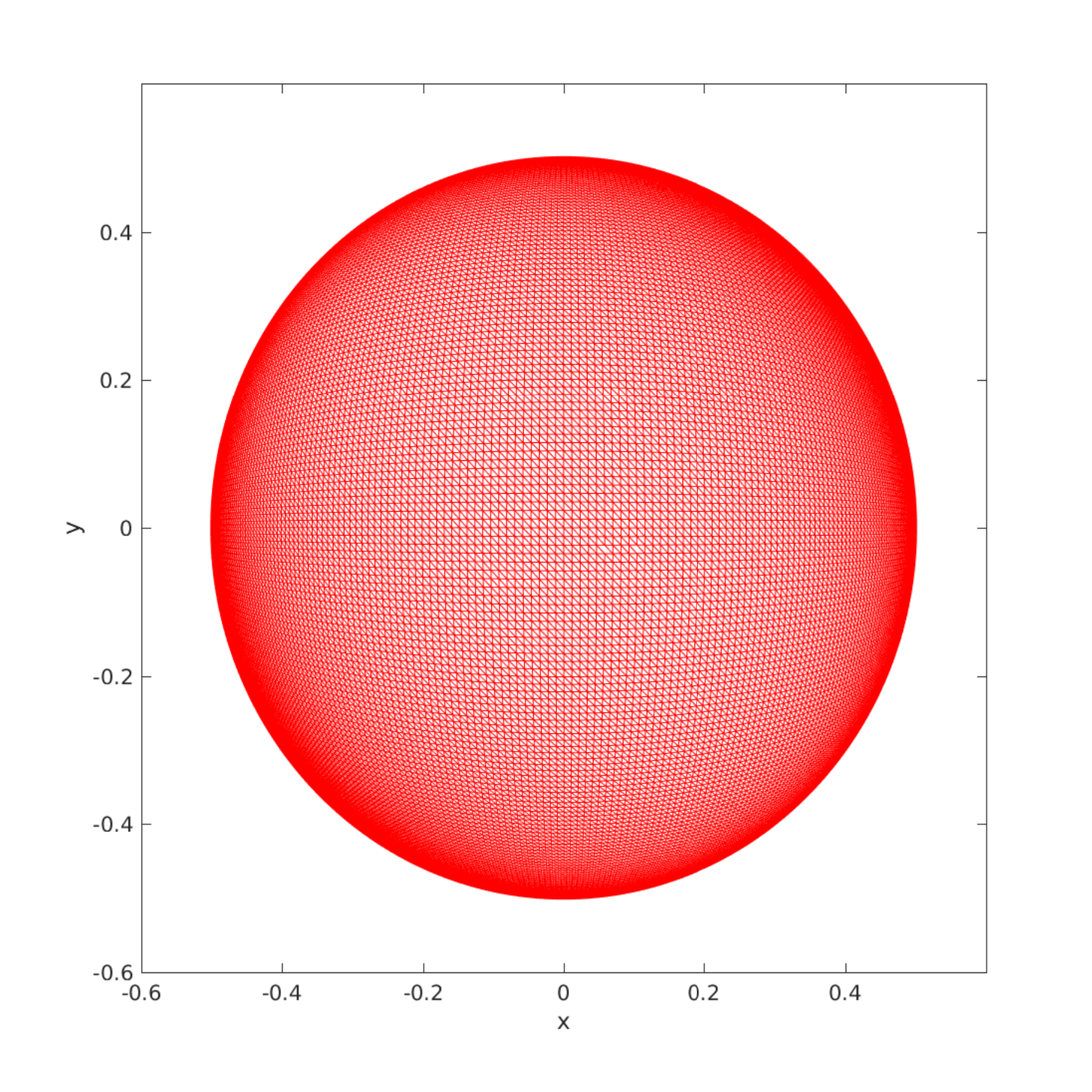}
       \caption{Mesh}
       \label{subfig:pmeV-mesh}
   \end{subfigure}\hspace{10mm}%
    \begin{subfigure}[b]{0.35\linewidth}
       \includegraphics[scale=0.30]{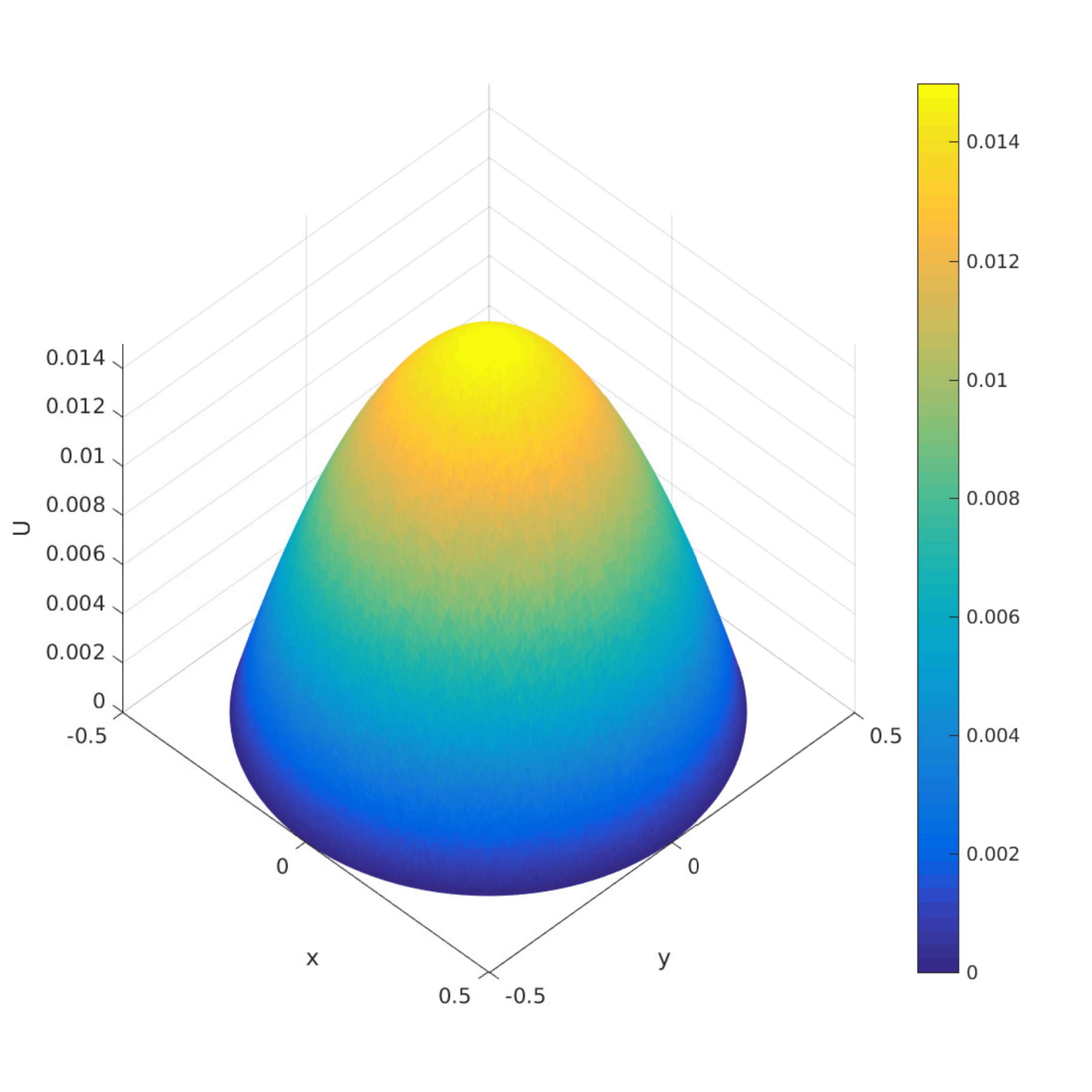}
       \caption{Computed solution}
       \label{subfig:pmeV-soln}
    \end{subfigure}%
    \caption{Example~\ref{ex:pme-V-BP}. The final mesh and computed solution are plotted for $m = 3$ ($N = 40459$).}
    \label{fig:pmeV-m3}
    \end{center}
\end{figure} 

\begin{figure}
    \centering
    \begin{subfigure}[b]{0.35\linewidth}
       \includegraphics[width=\textwidth]{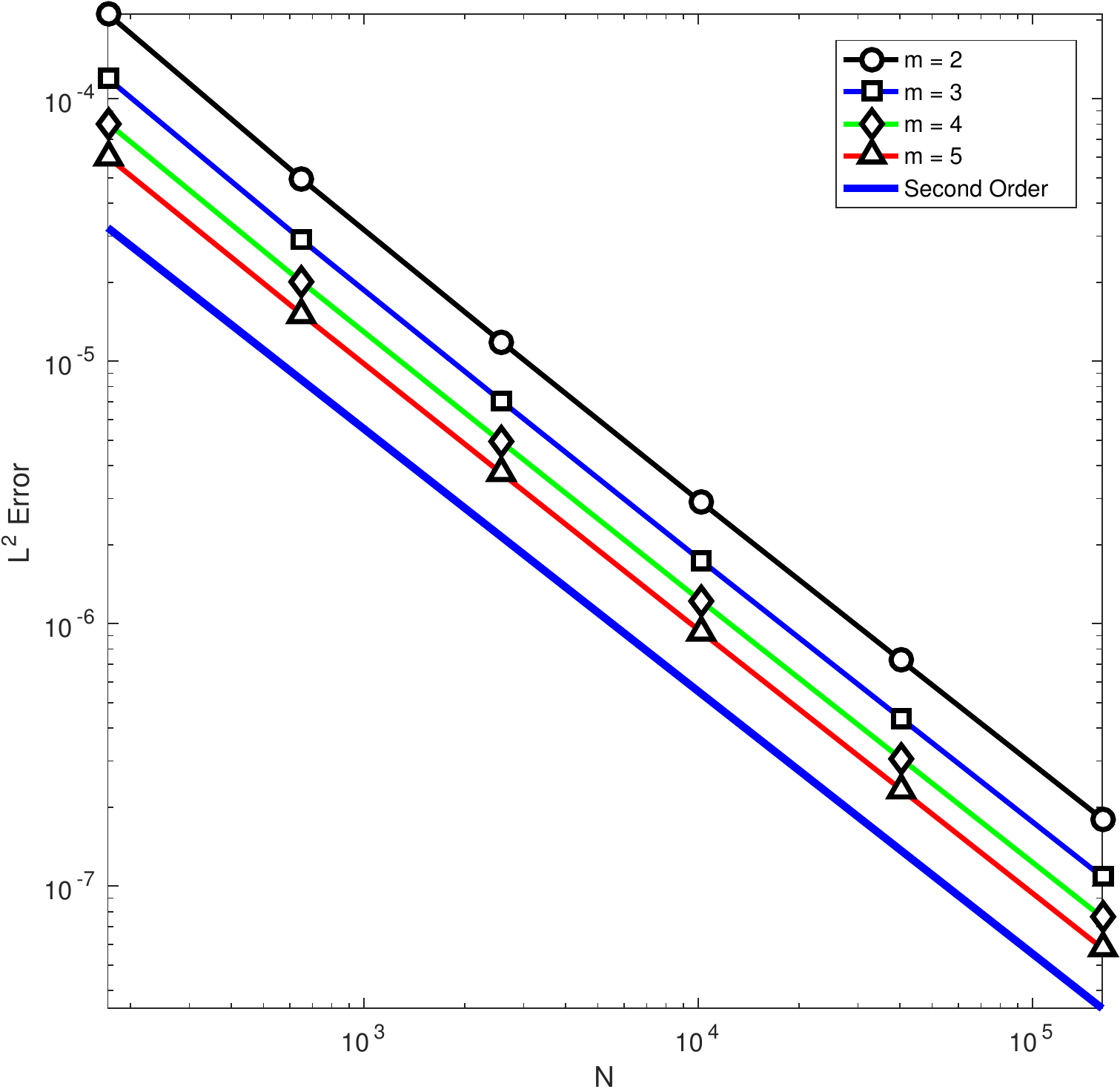}
       \caption{$L^2$ error of solution}
       \label{subfig:pmeV-m}
    \end{subfigure}\hspace{5mm}%
    \begin{subfigure}[b]{0.35\linewidth}
       \includegraphics[width=\textwidth]{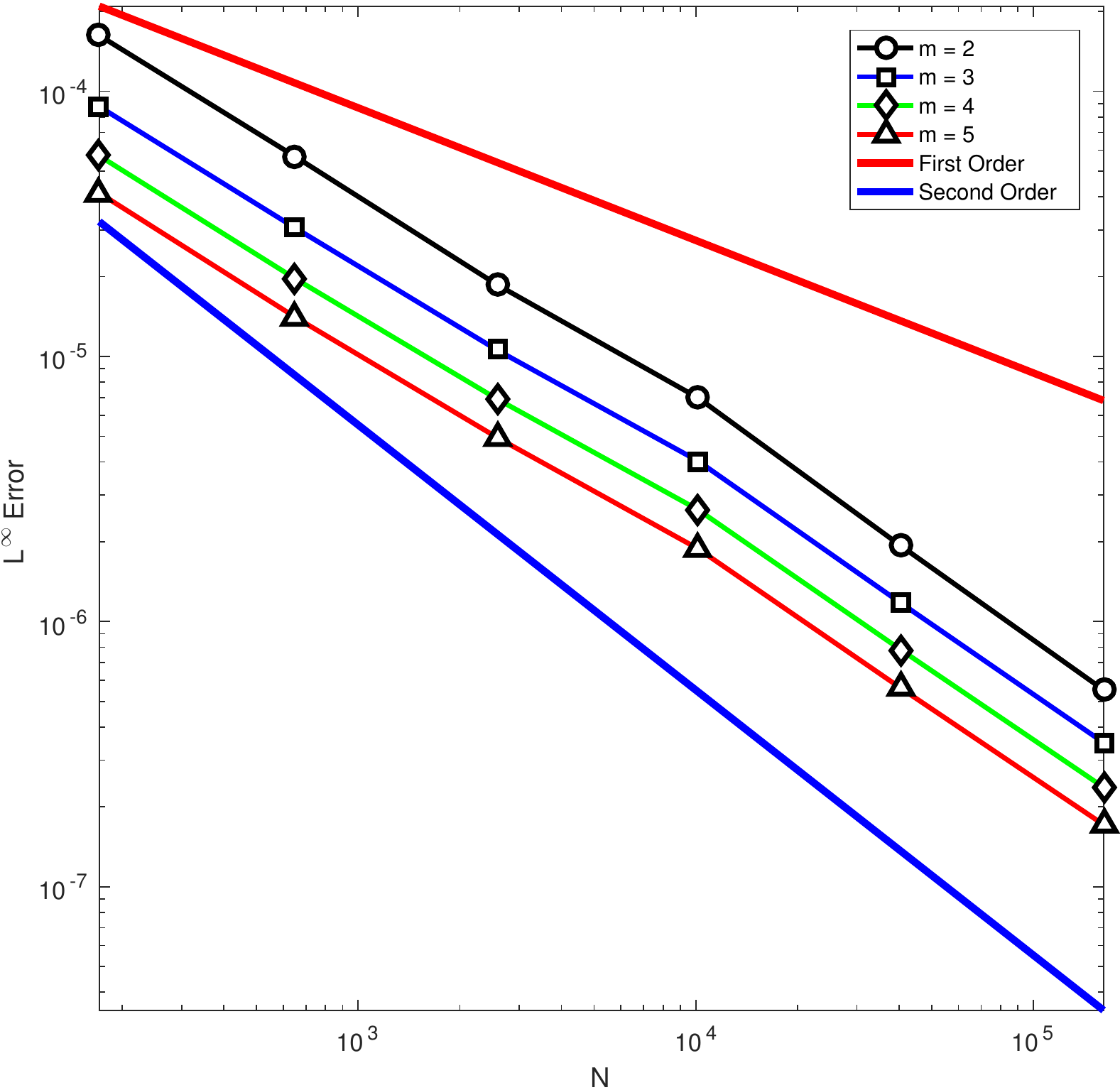}
       \caption{$L^{\infty}$ error of the boundary}
       \label{subfig:pmeV-m-bdry}
    \end{subfigure}%
    \caption{Example \ref{ex:pme-V-BP}. Convergence histories of the moving mesh method ($v$-variable).}
    \label{fig:pme-V-conv}
\end{figure}

\begin{figure}
    \centering
    \begin{subfigure}[b]{0.35\linewidth}
       \includegraphics[width=\textwidth]{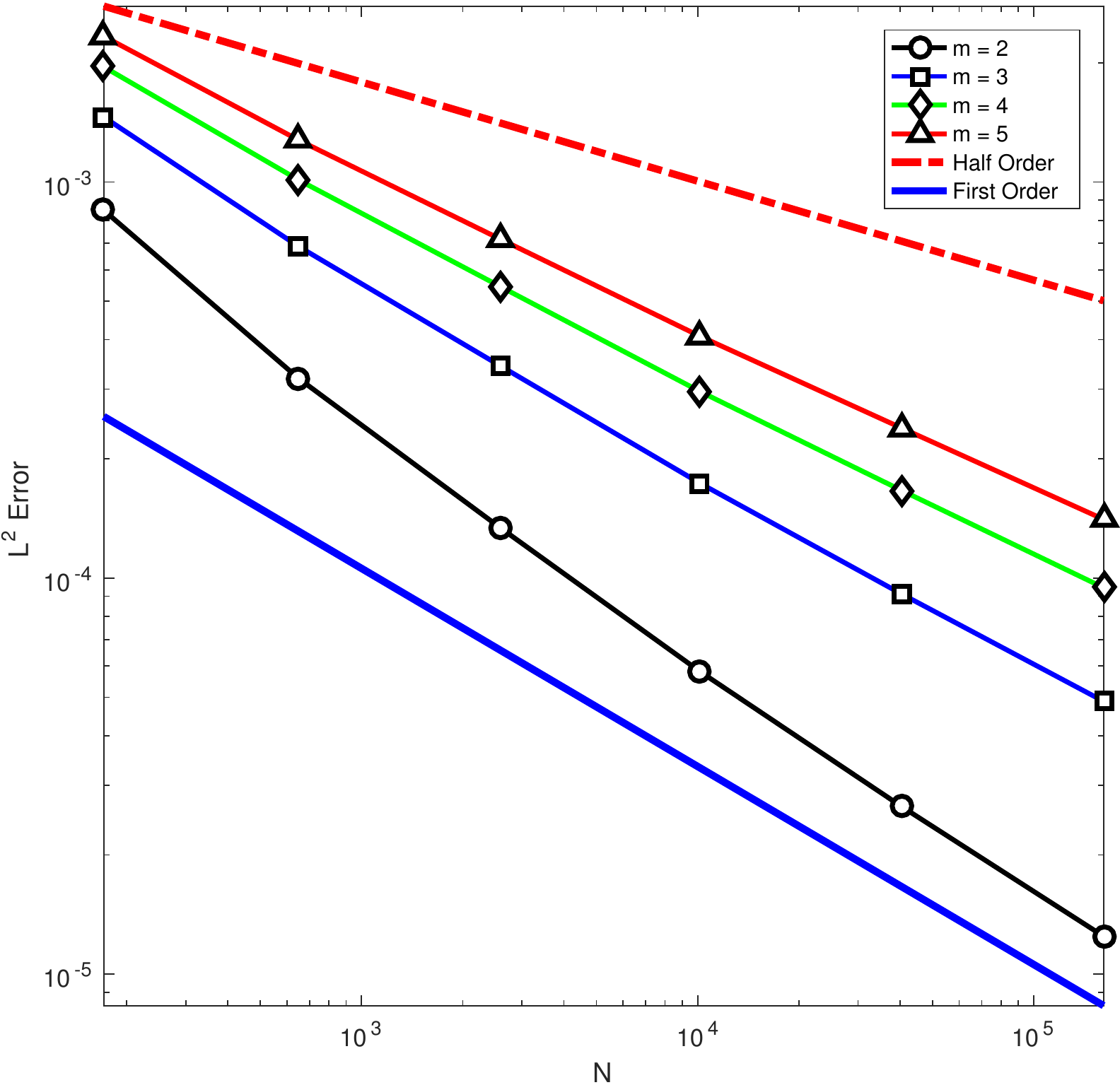}
       \caption{$L^2$ Error}
       \label{subfig:pmeUV-m-L2}
    \end{subfigure}\hspace{5mm}%
    \begin{subfigure}[b]{0.35\linewidth}
       \includegraphics[width=\textwidth]{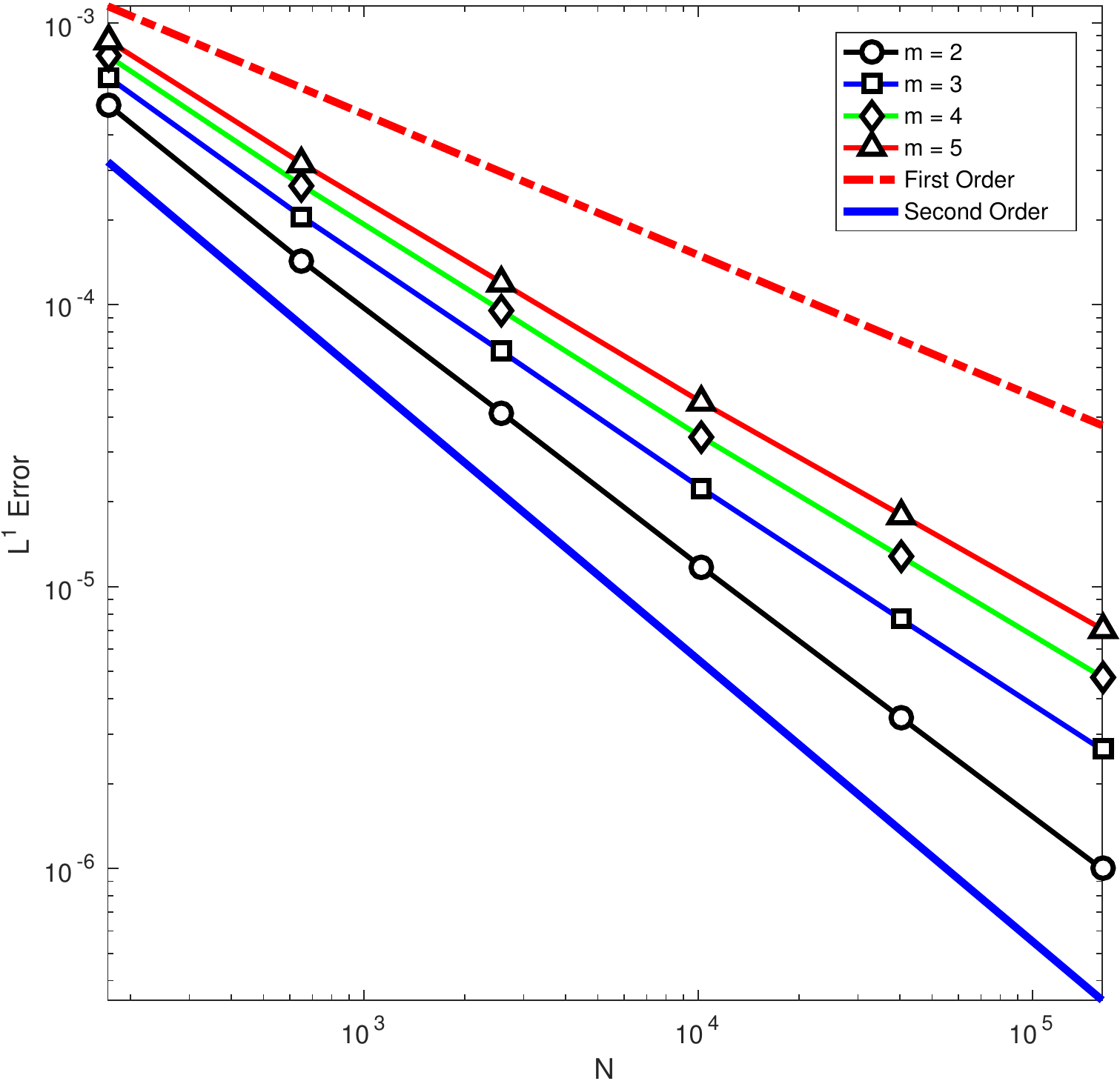}
       \caption{$L^1$ Error}
       \label{subfig:pmeUV-m-L1}
    \end{subfigure}%
    \caption{Example \ref{ex:pme-V-BP}. Convergence history of the method in the original $u$-variable.}
    \label{fig:pme-UV-conv}
\end{figure}

\begin{figure}
    \centering
    \begin{subfigure}[b]{0.35\linewidth}
       \includegraphics[width=\textwidth]{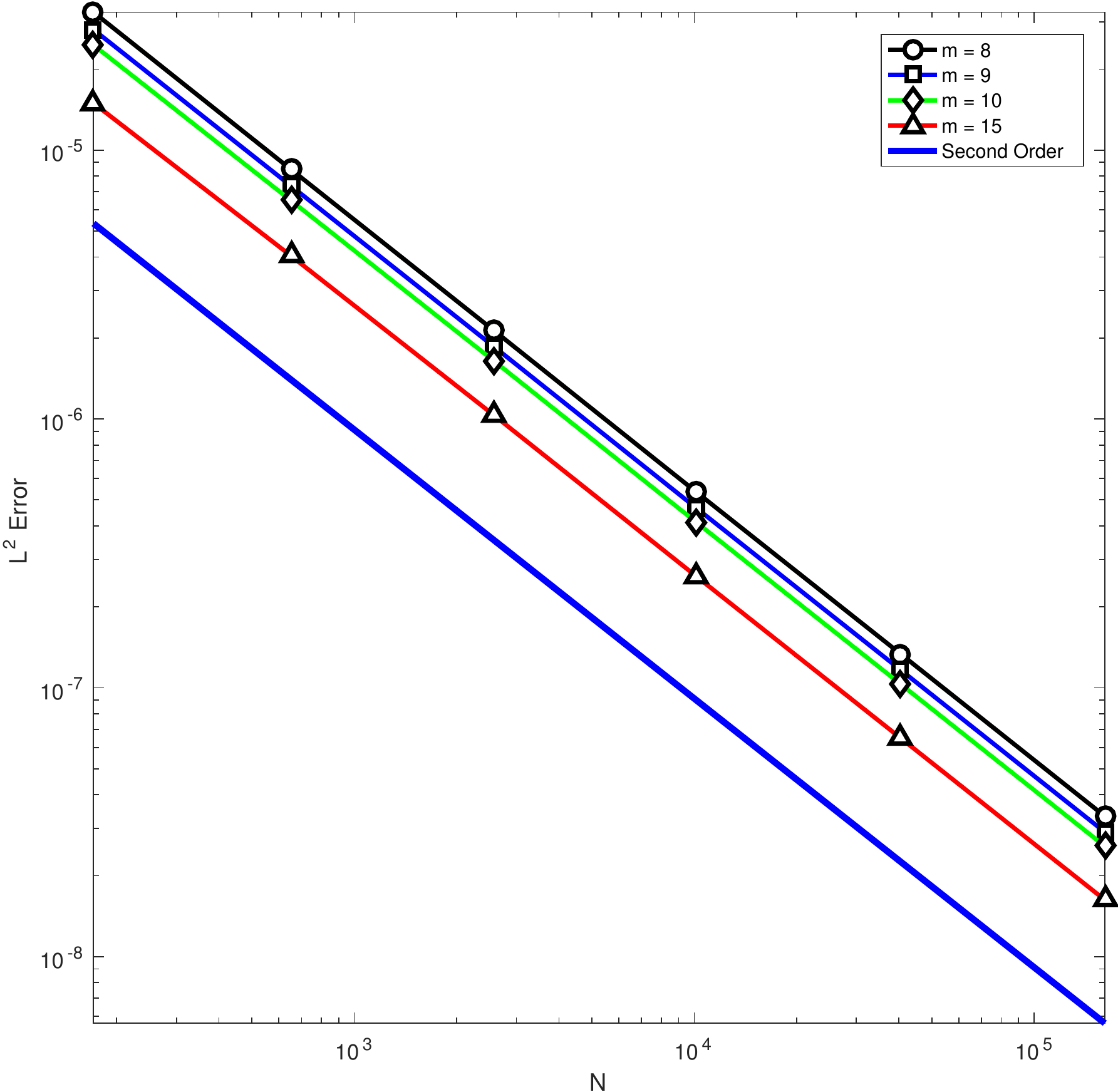}
       \caption{$L^2$ error of solution}
       \label{subfig:pmeV-m}
    \end{subfigure}\hspace{5mm}%
    \begin{subfigure}[b]{0.35\linewidth}
       \includegraphics[width=\textwidth]{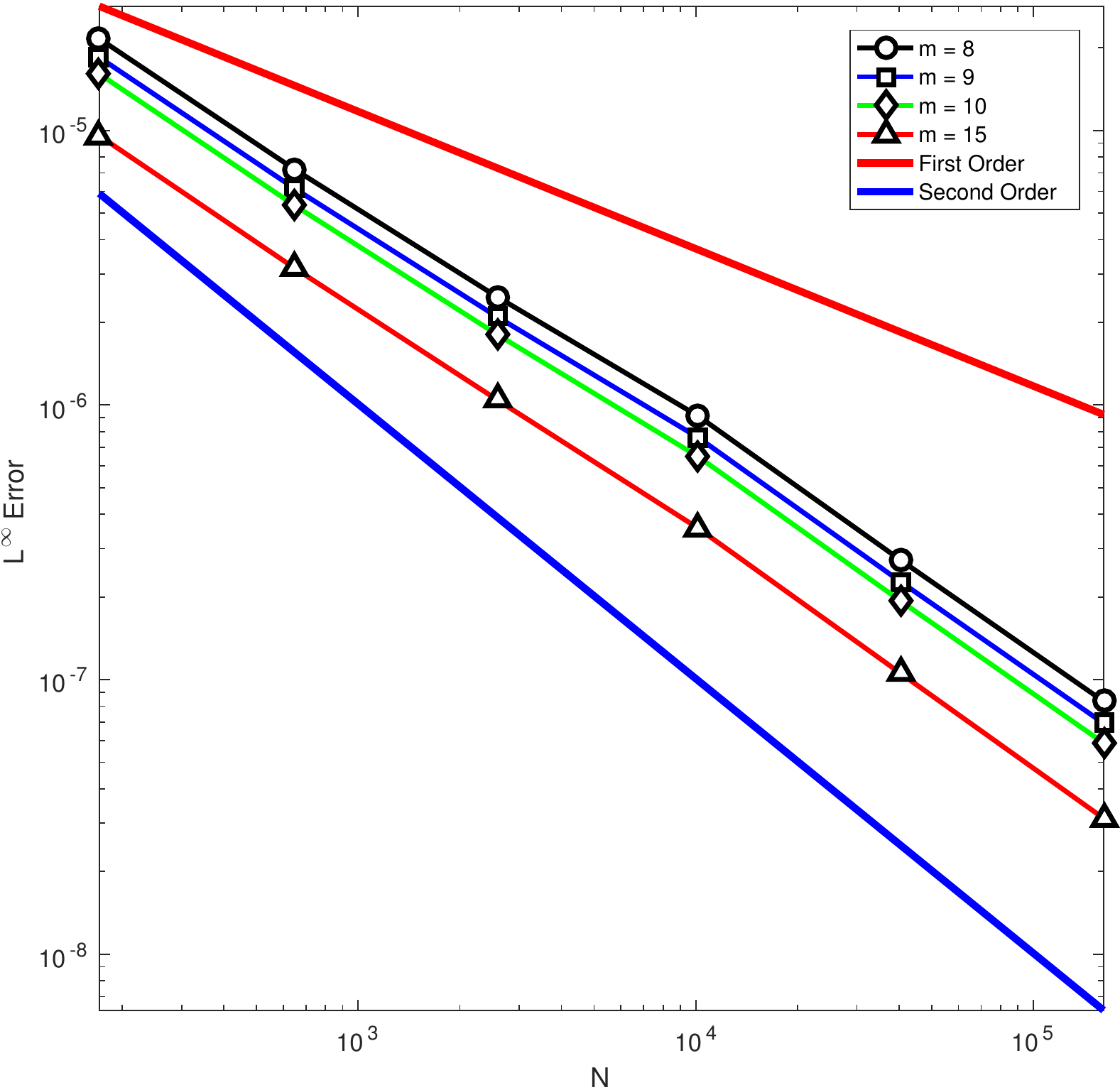}
       \caption{$L^{\infty}$ error of the boundary}
       \label{subfig:pmeV-m-bdry}
    \end{subfigure}%
    \caption{Example \ref{ex:pme-V-BP}. Convergence histories of the moving mesh method ($v$-variable)
    for $m=8,\; 9, \; 10$, and $15$.}
    \label{fig:pme-V-conv-Bigm}
\end{figure}

\begin{figure}
    \centering
    \begin{subfigure}[b]{0.35\linewidth}
       \includegraphics[width=\textwidth]{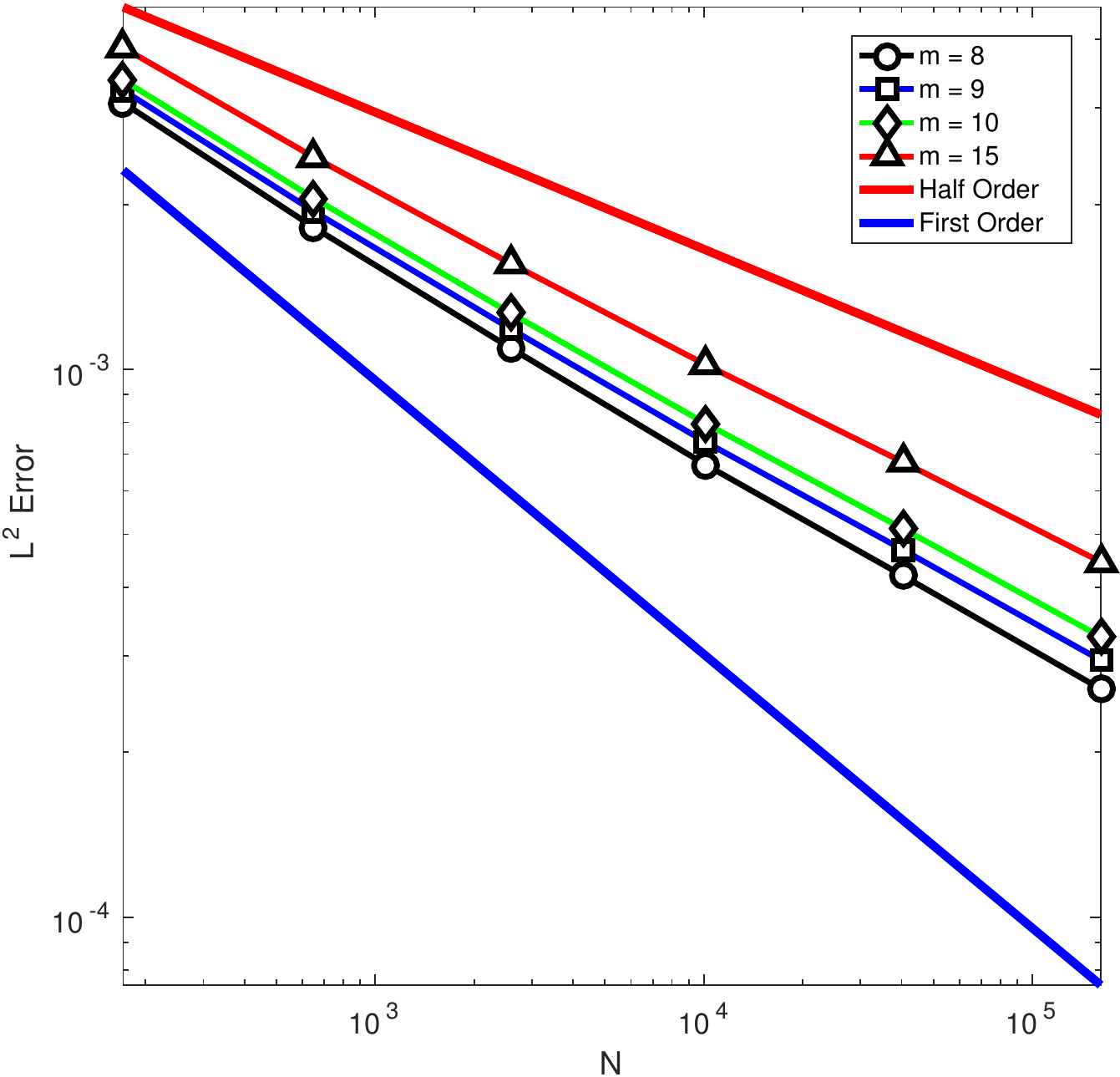}
       \caption{$L^2$ Error}
       \label{subfig:pmeUV-m-L2}
    \end{subfigure}\hspace{5mm}%
    \begin{subfigure}[b]{0.35\linewidth}
       \includegraphics[width=\textwidth]{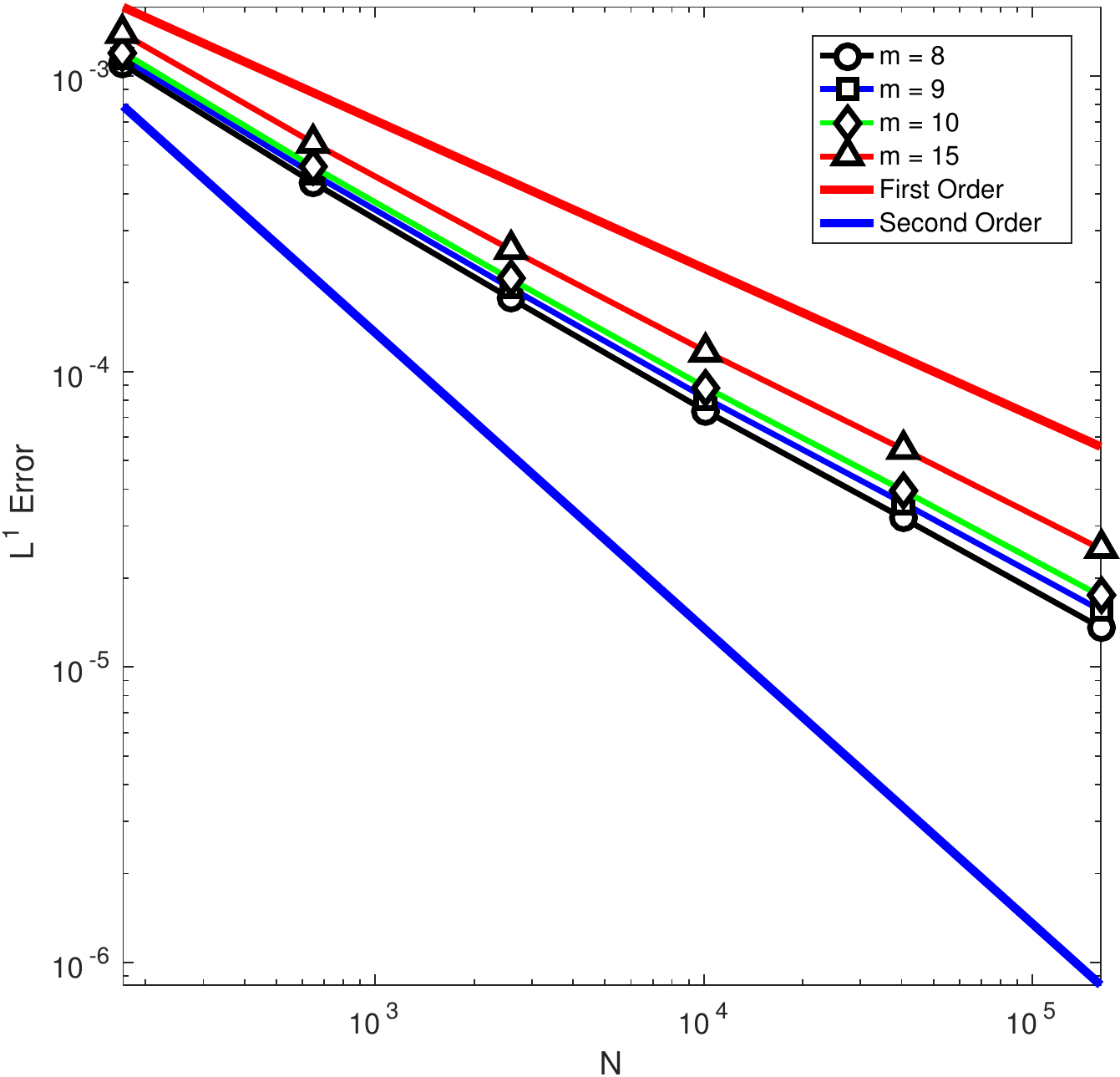}
       \caption{$L^1$ Error}
       \label{subfig:pmeUV-m-L1}
    \end{subfigure}%
    \caption{Example \ref{ex:pme-V-BP}. Convergence history of the method in the original $u$-variable
    for $m=8,\; 9, \; 10$, and $15$.}
    \label{fig:pme-UV-conv-Bigm}
\end{figure}

\begin{figure}
    \centering
    \begin{subfigure}[b]{0.35\linewidth}
       \includegraphics[width=\textwidth]{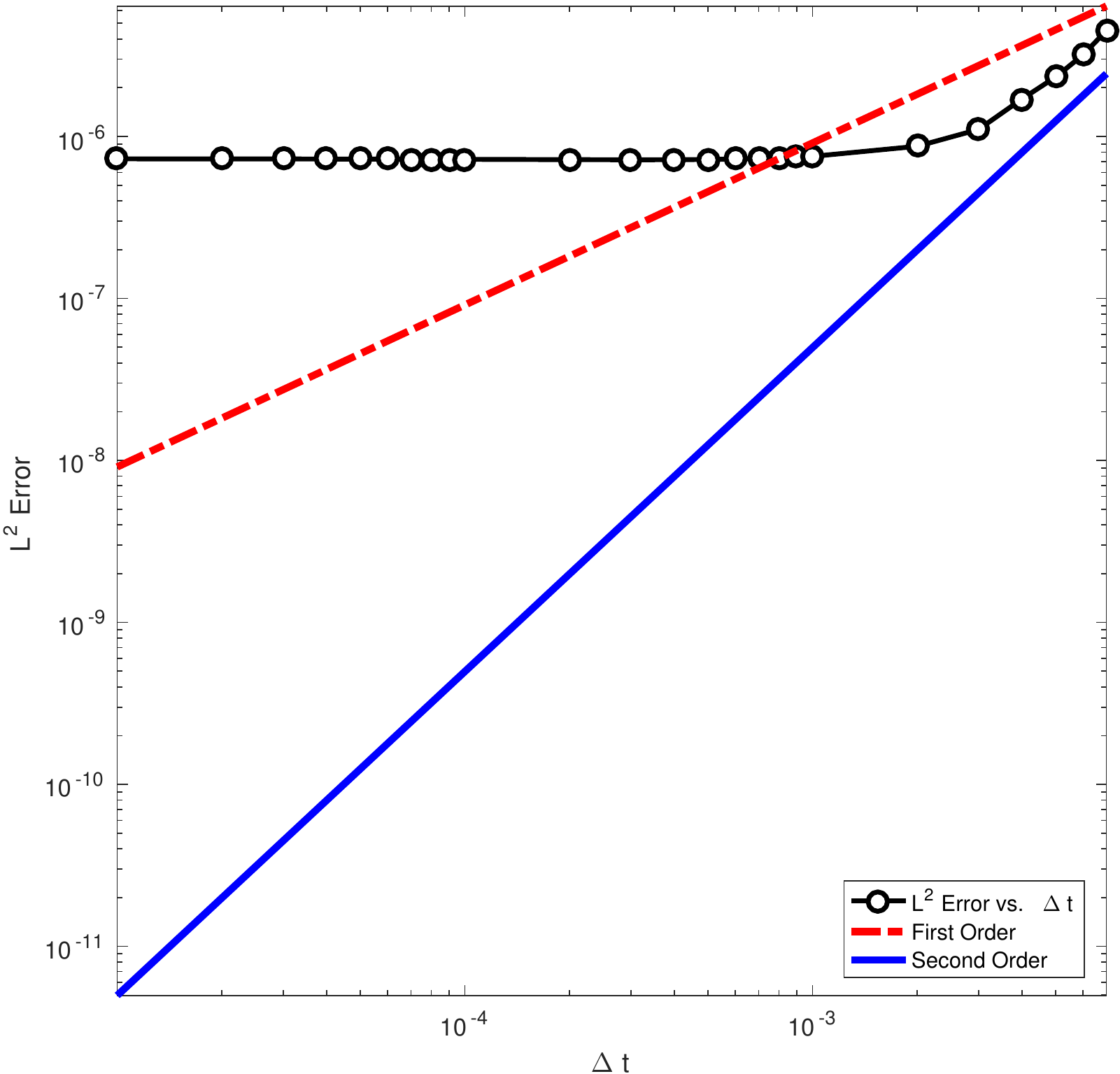}
       \caption{$L^2$ Error of solution}
       \label{subfig:dtFixed-1}
    \end{subfigure}\hspace{5mm}%
    \begin{subfigure}[b]{0.35\linewidth}
       \includegraphics[width=\textwidth]{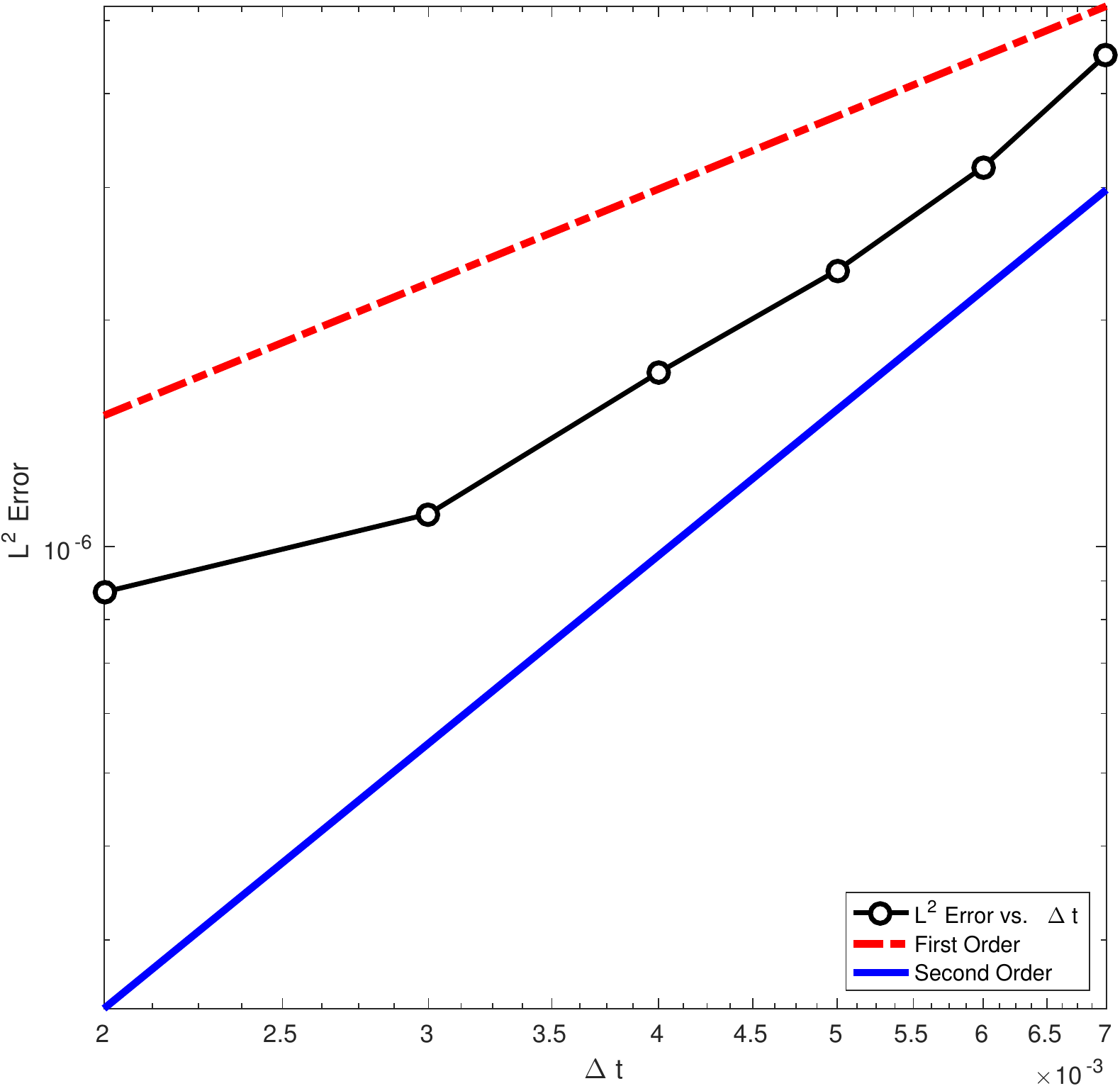}
       \caption{A closer look}
       \label{subfig:dtFixed-2}
    \end{subfigure}%
    \caption{Example \ref{ex:pme-V-BP} with $m=2$.
    Error history of the moving mesh method ($v$-variable) when compared with fixed time steps.}
    \label{fig:pme-V-dtFixed}
\end{figure}

\begin{figure}[ht]
        \vspace{-2cm}
        \centering
        \begin{subfigure}[b]{0.35\linewidth}\includegraphics[scale=0.28]{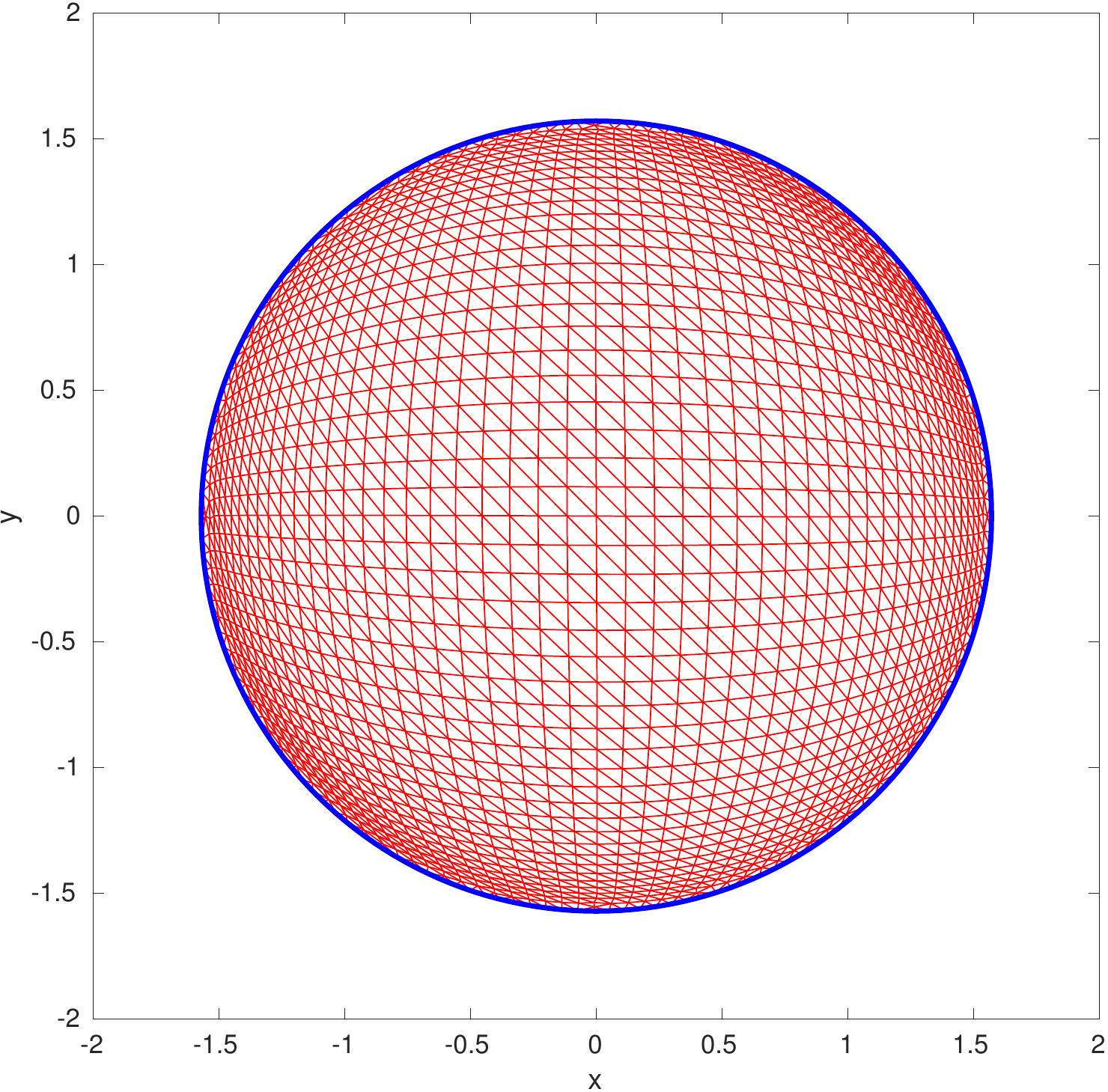}\caption{$t = 0.11$}\end{subfigure}\hspace{5mm}
        \begin{subfigure}[b]{0.35\linewidth}\includegraphics[scale=0.28]{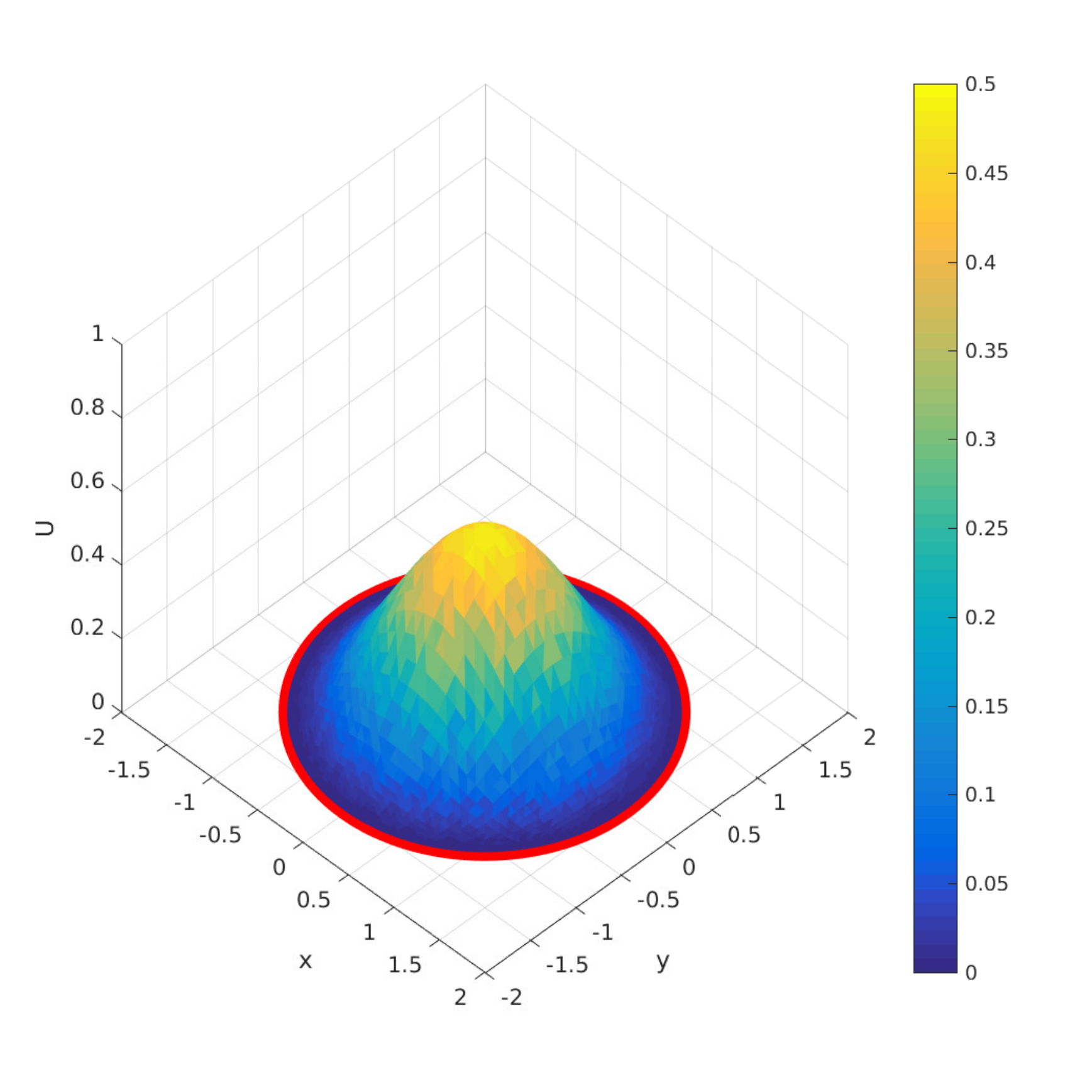}\caption{$t = 0.11$}\end{subfigure}\\%
        \begin{subfigure}[b]{0.35\linewidth}\includegraphics[scale=0.28]{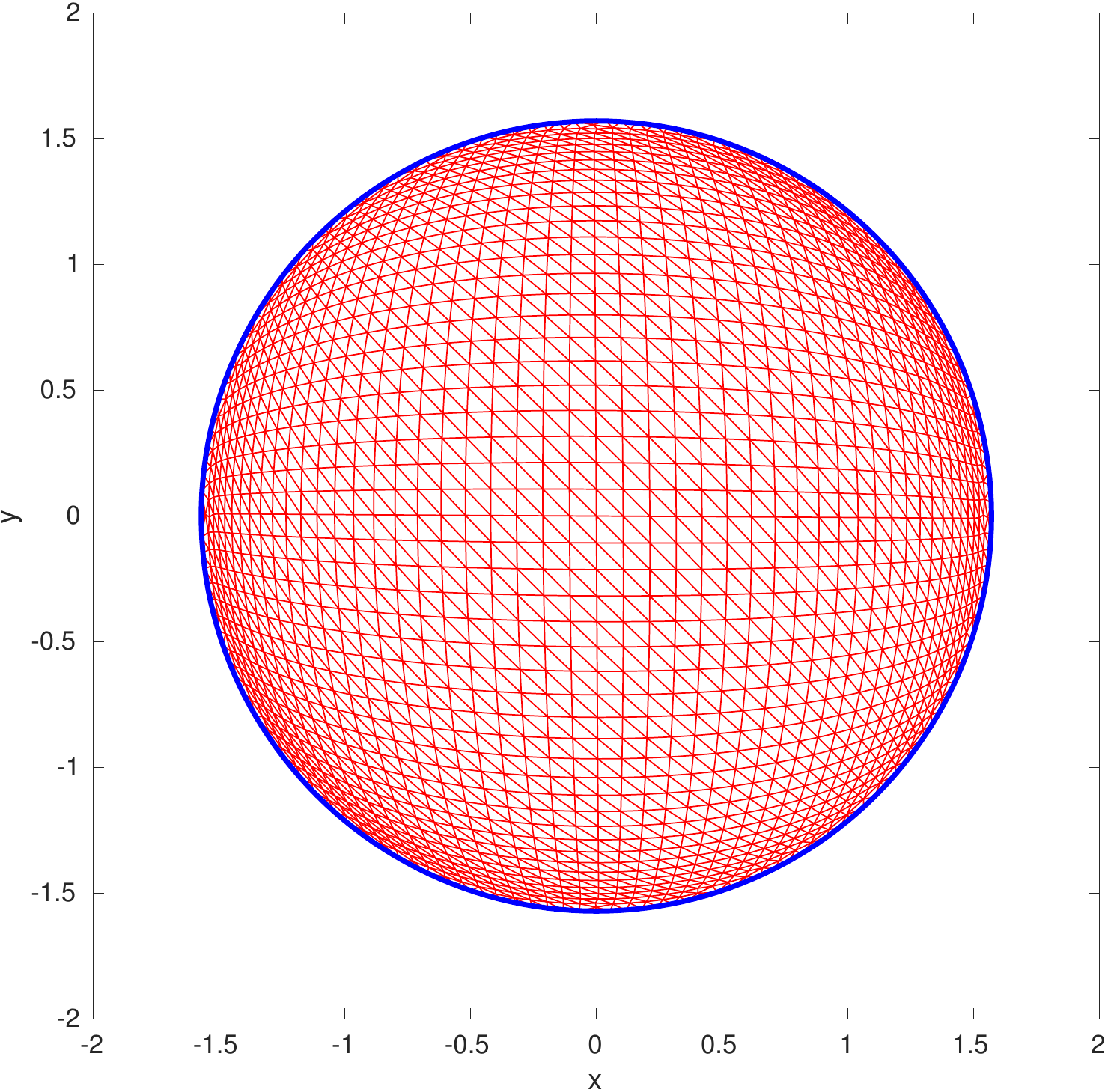}\caption{$t = 0.2$}\end{subfigure}\hspace{5mm}
        \begin{subfigure}[b]{0.35\linewidth}\includegraphics[scale=0.28]{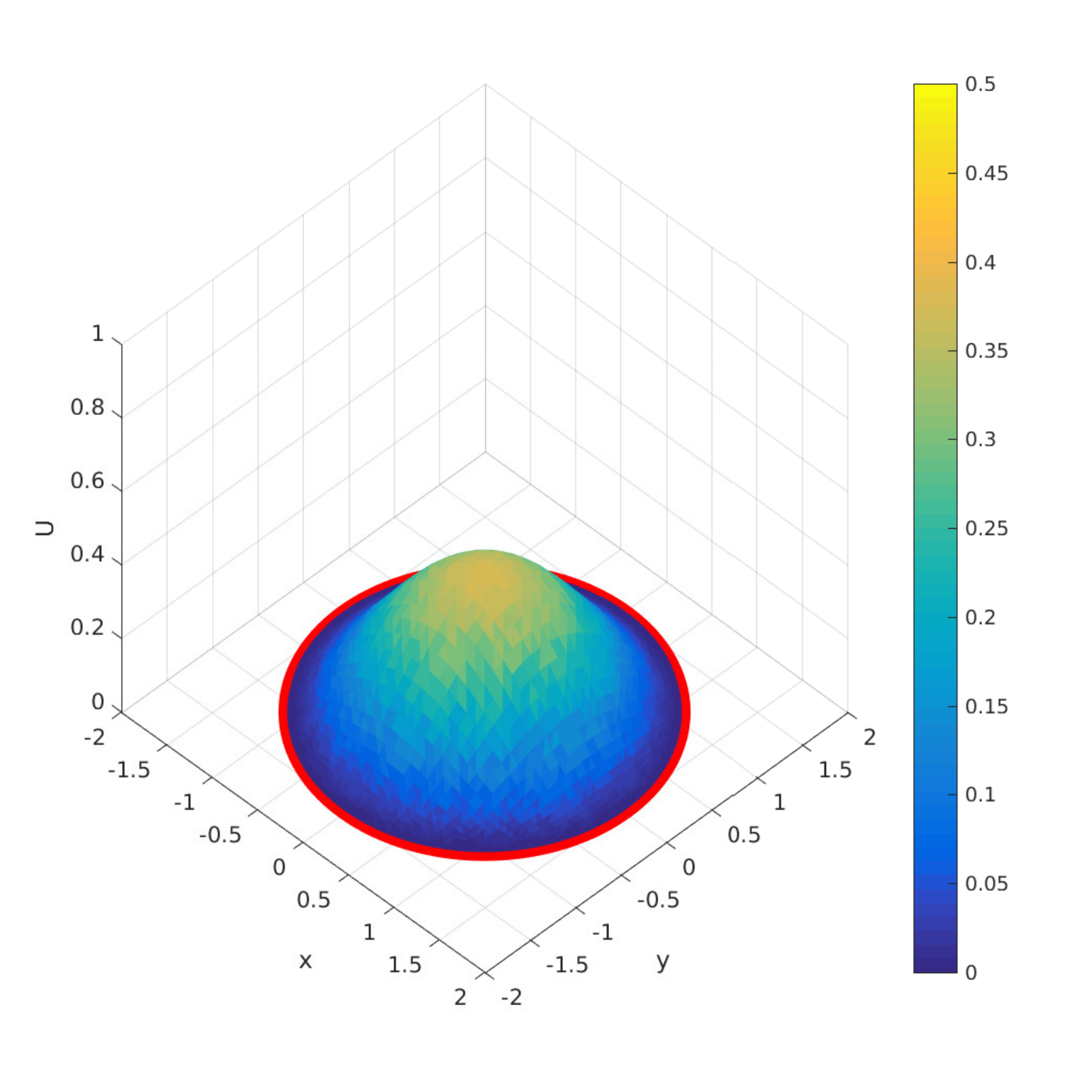}\caption{$t = 0.2$}\end{subfigure}\\%
        \begin{subfigure}[b]{0.35\linewidth}\includegraphics[scale=0.28]{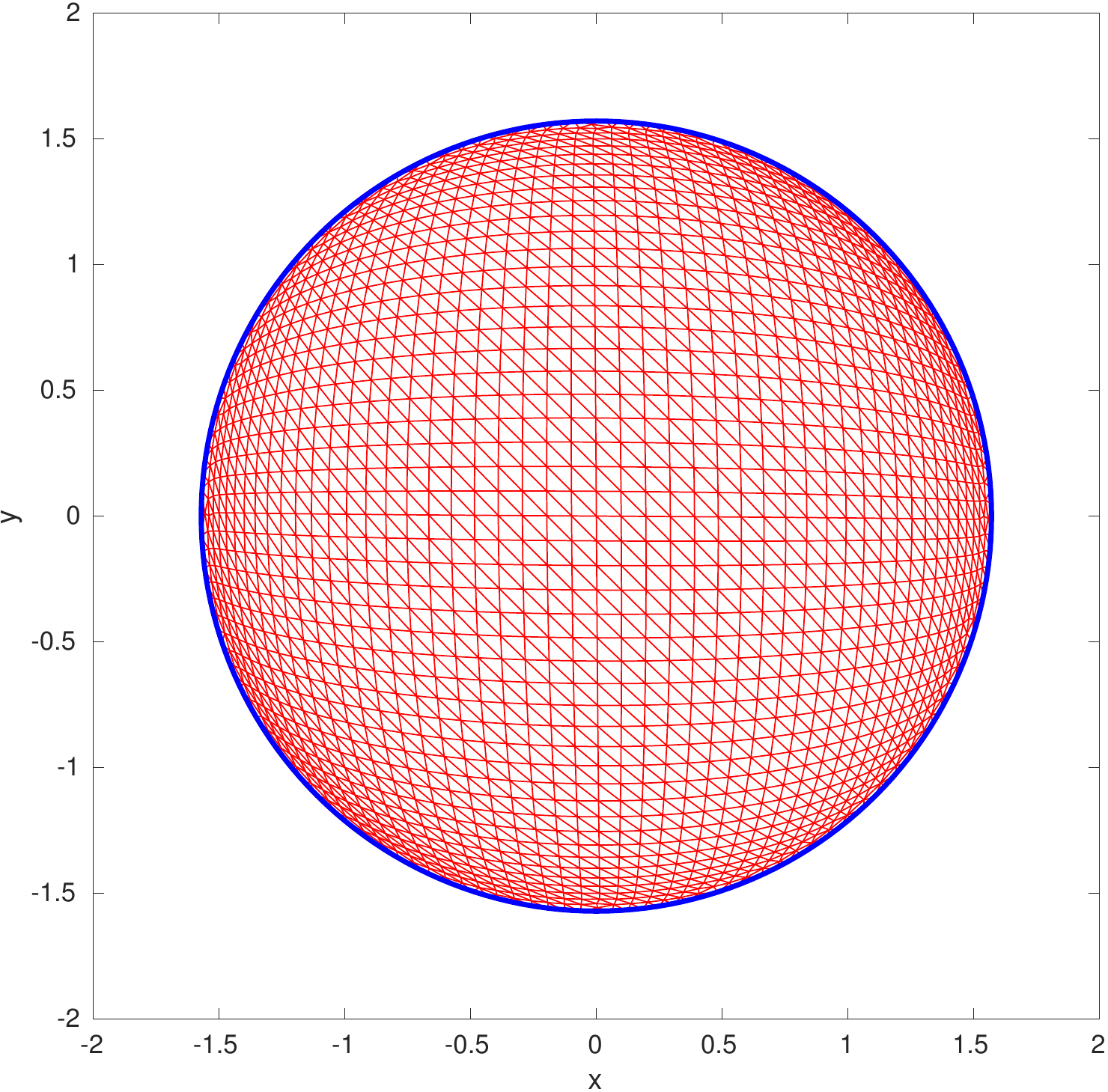}\caption{$t = 0.3$}\end{subfigure}\hspace{5mm}
        \begin{subfigure}[b]{0.35\linewidth}\includegraphics[scale=0.28]{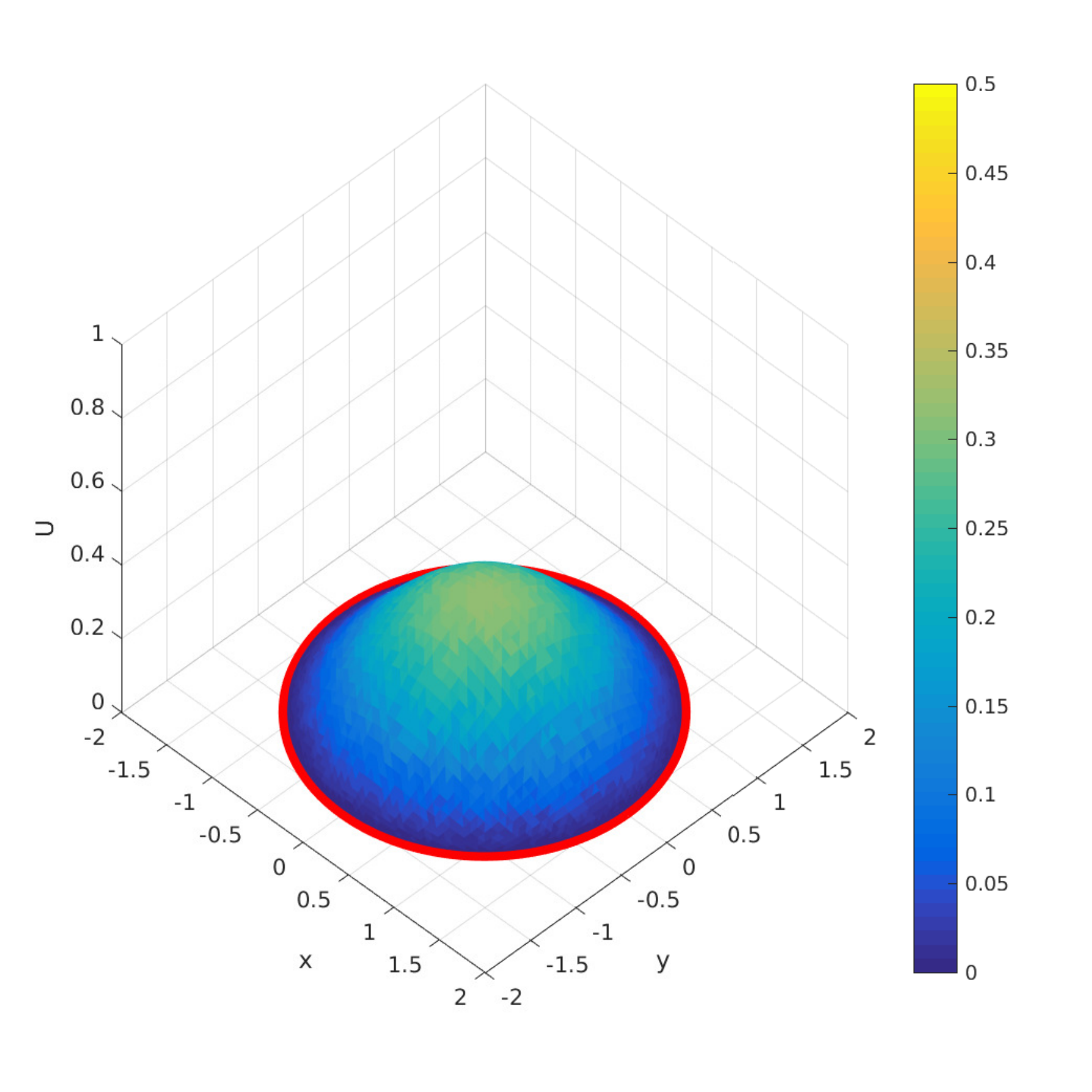}\caption{$t = 0.3$}\end{subfigure}\\%
        \begin{subfigure}[b]{0.35\linewidth}\includegraphics[scale=0.28]{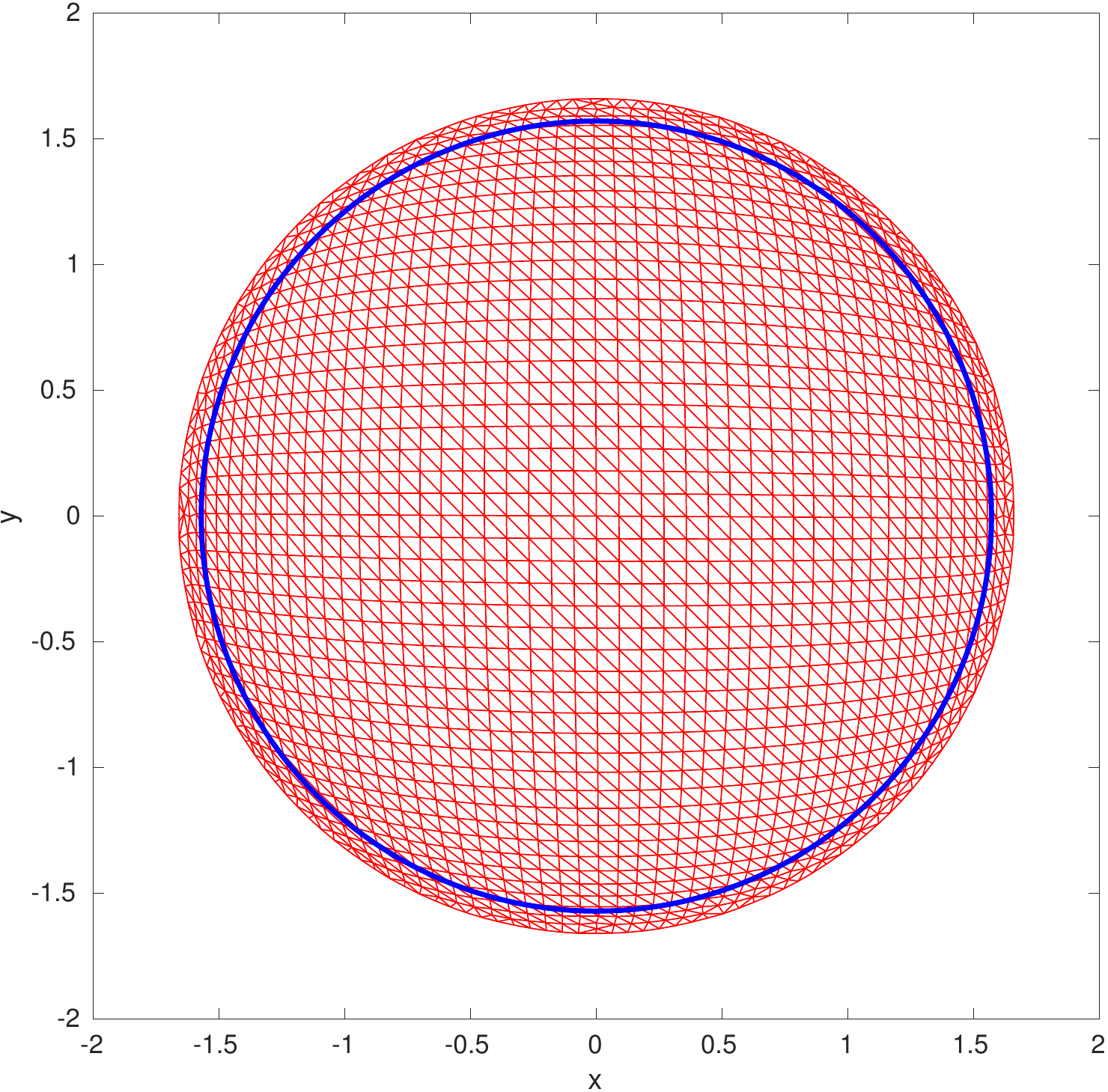}\caption{$t = 0.8$}\end{subfigure}\hspace{5mm}
        \begin{subfigure}[b]{0.35\linewidth}\includegraphics[scale=0.28]{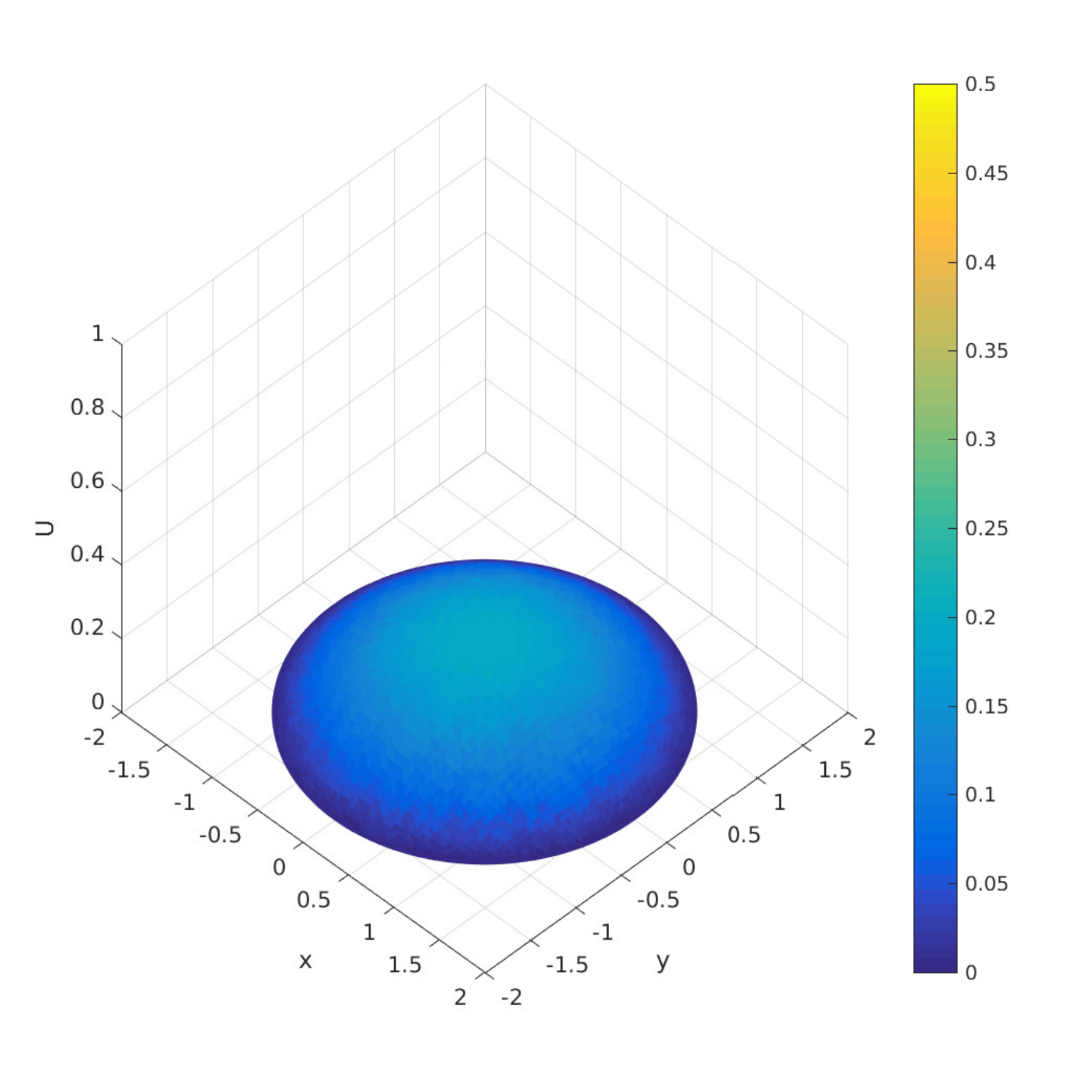}\caption{$t = 0.8$}\end{subfigure}%
        \caption{Example~\ref{ex:pme-V-waiting}. The mesh and computed solution are
        shown at various time instants ($N = 4011$). The contrasting circles represent the initial boundary.}
        \label{fig:pme-V-waiting}
    \end{figure}

\begin{figure}[ht]
        \centering
        \begin{subfigure}[b]{0.23\linewidth}\includegraphics[width=\textwidth]{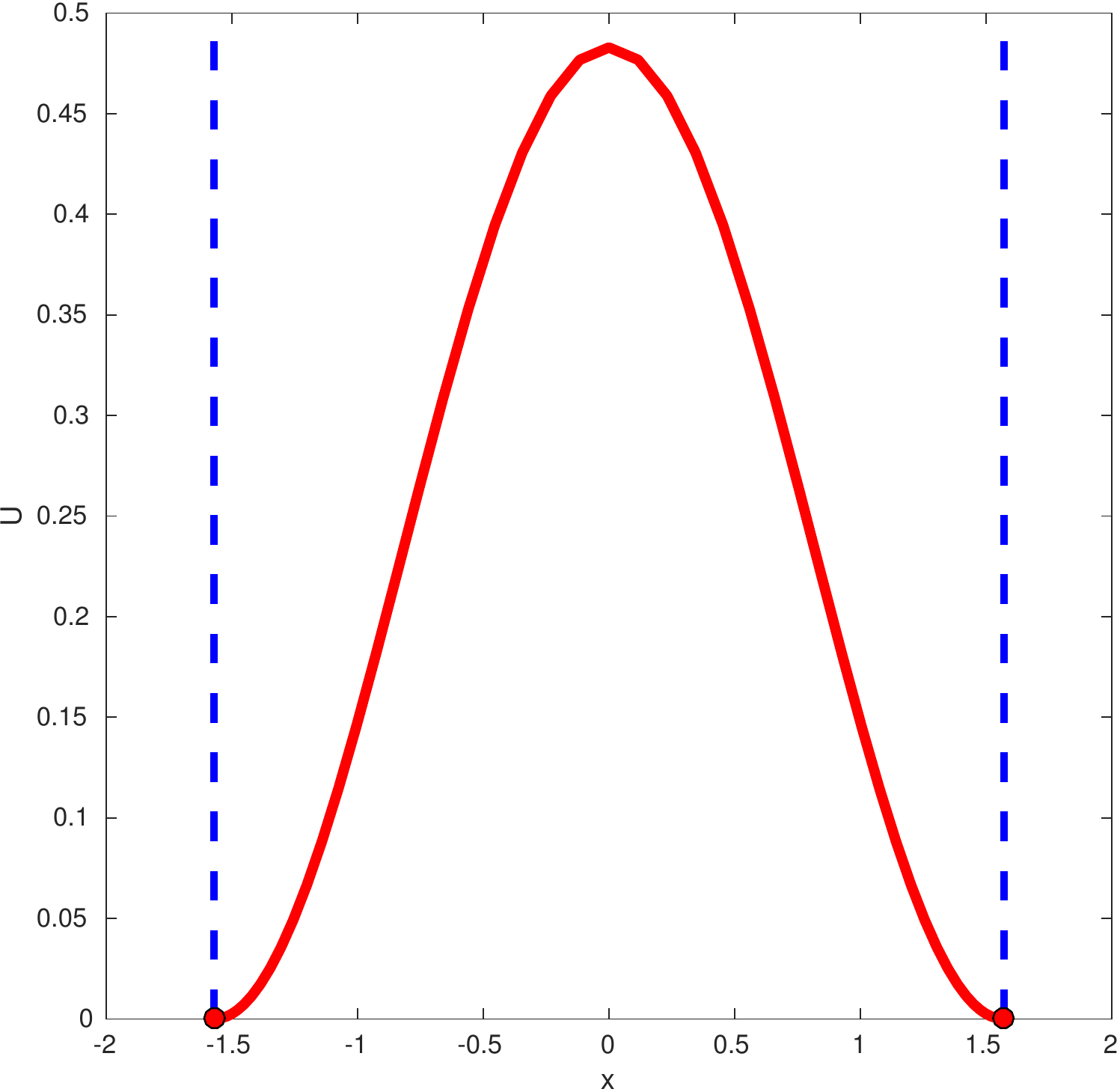}\caption{$t = 0.11$}\end{subfigure}\hspace{2mm}%
        \begin{subfigure}[b]{0.23\linewidth}\includegraphics[width=\textwidth]{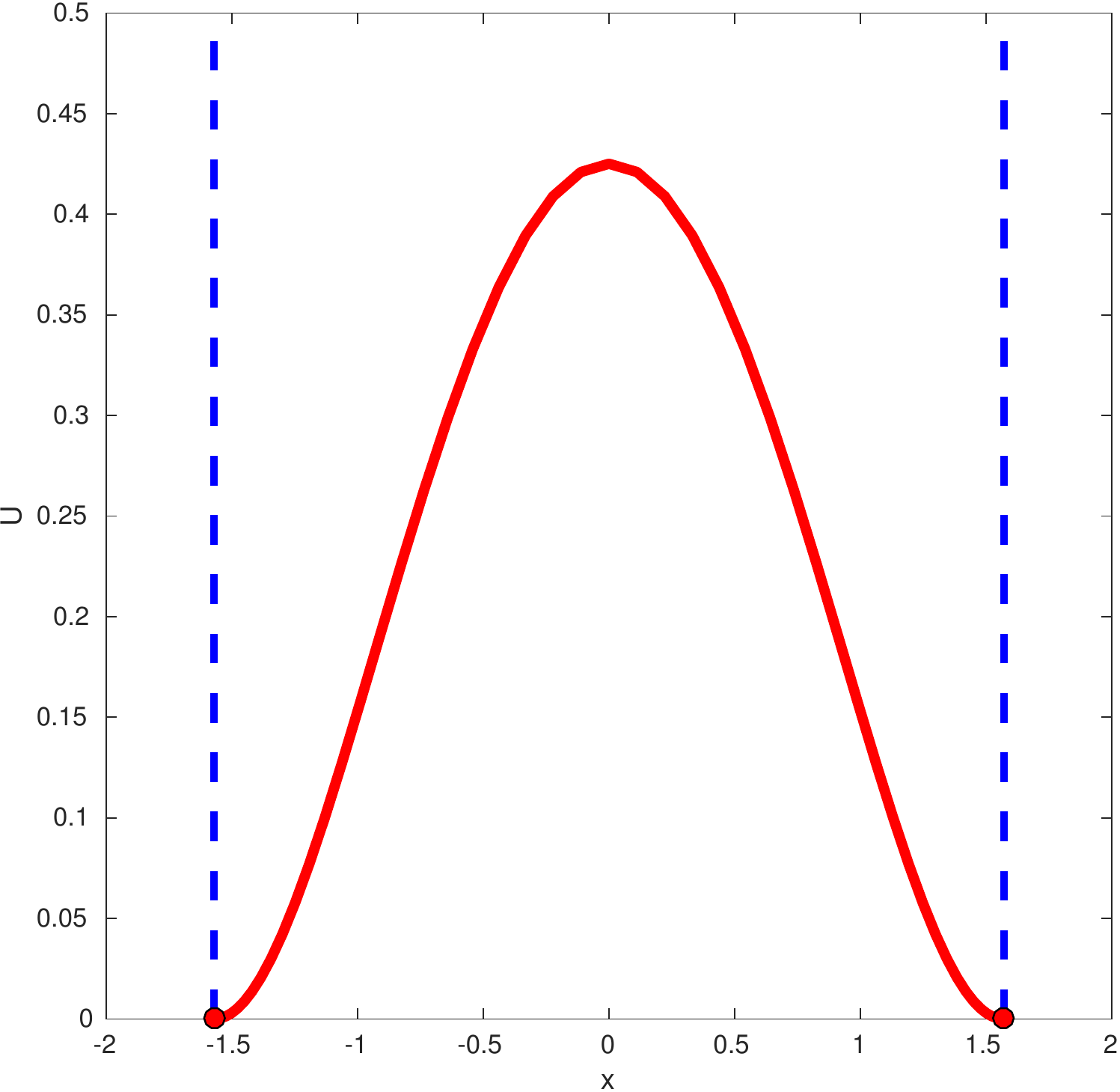}\caption{$t = 0.15$}\end{subfigure}\hspace{2mm}%
        \begin{subfigure}[b]{0.23\linewidth}\includegraphics[width=\textwidth]{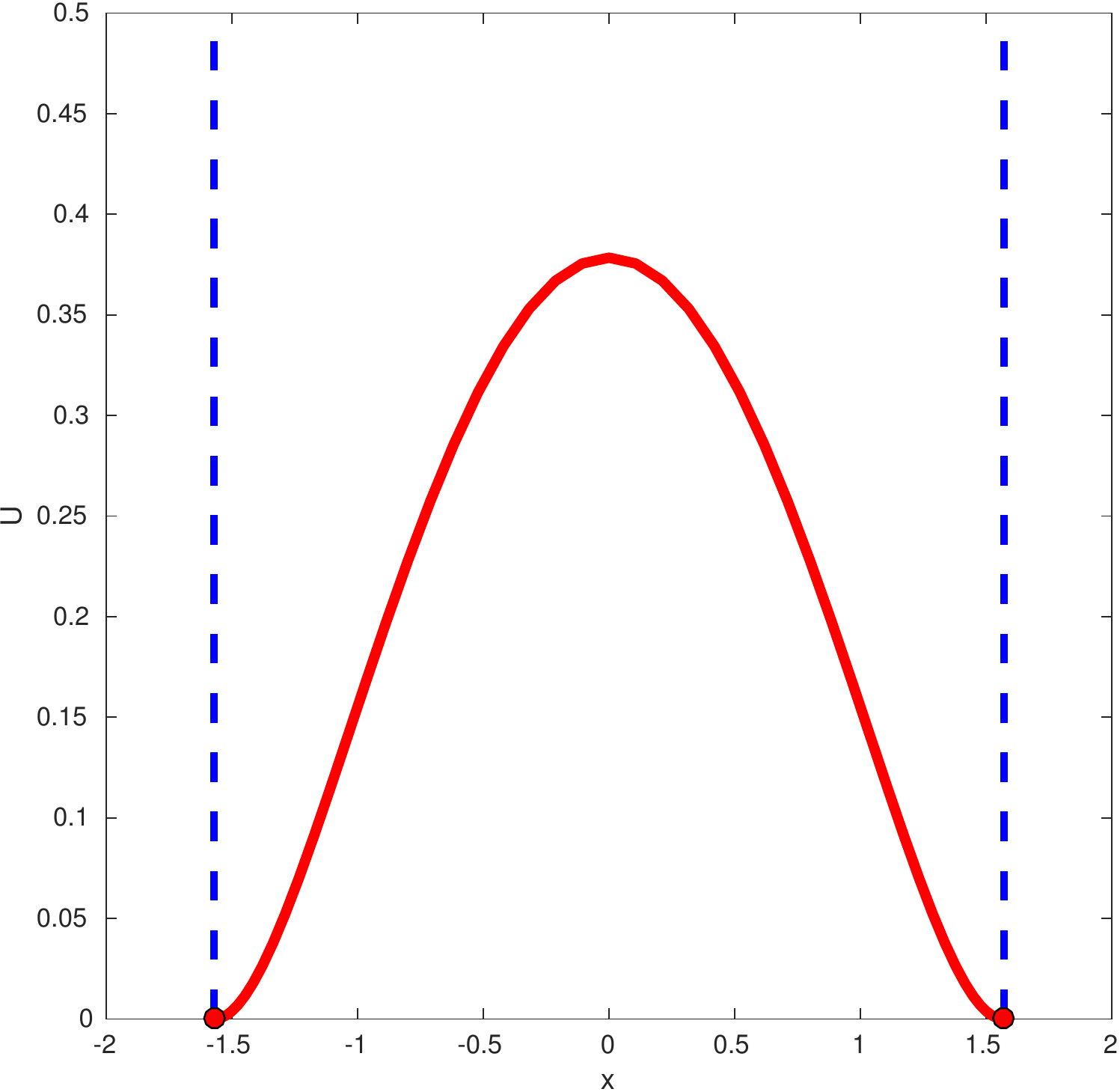}\caption{$t = 0.20$}\end{subfigure}\hspace{2mm}%
        \begin{subfigure}[b]{0.23\linewidth}\includegraphics[width=\textwidth]{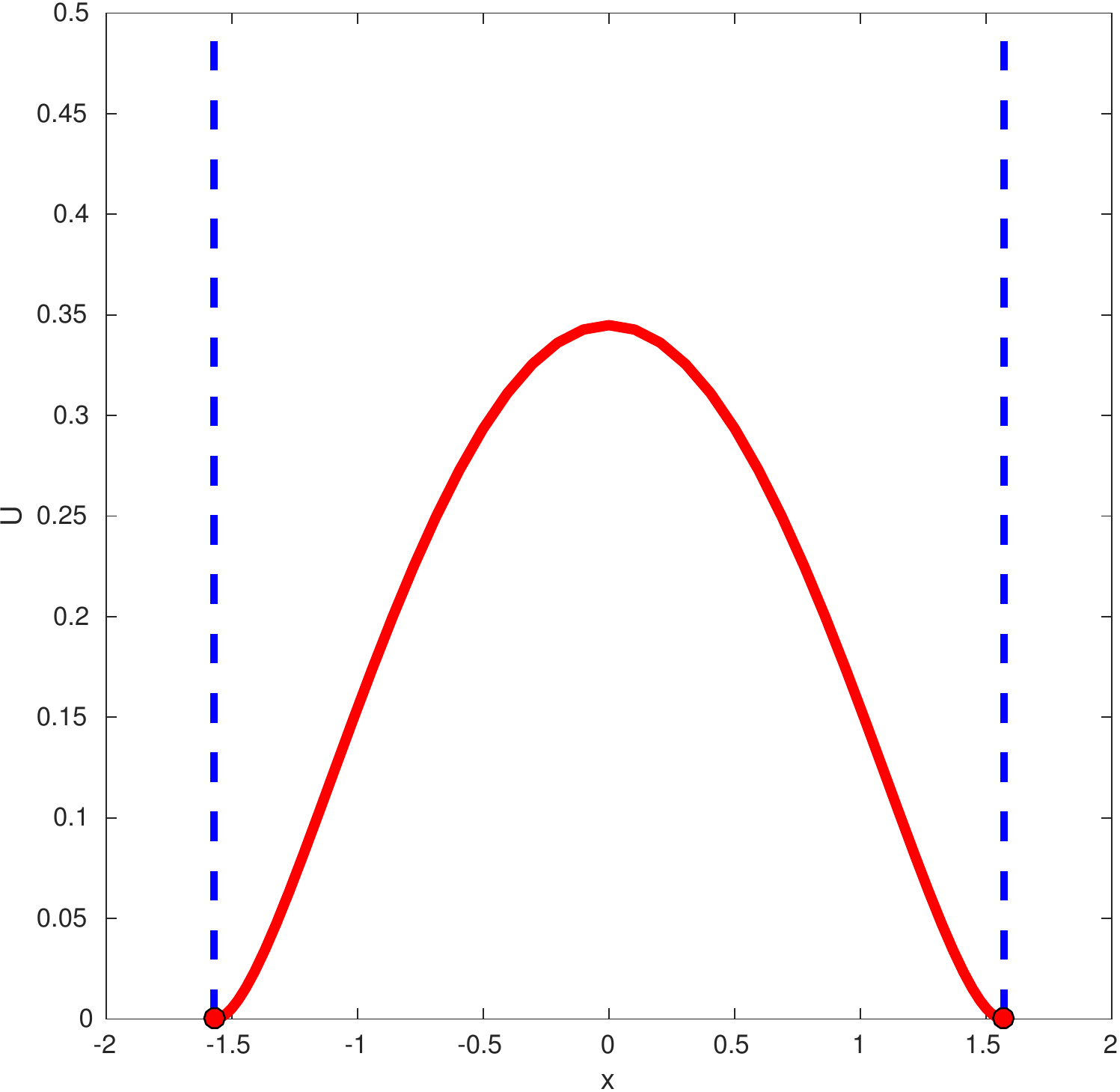}\caption{$t = 0.25$}\end{subfigure}\\%
        \begin{subfigure}[b]{0.23\linewidth}\includegraphics[width=\textwidth]{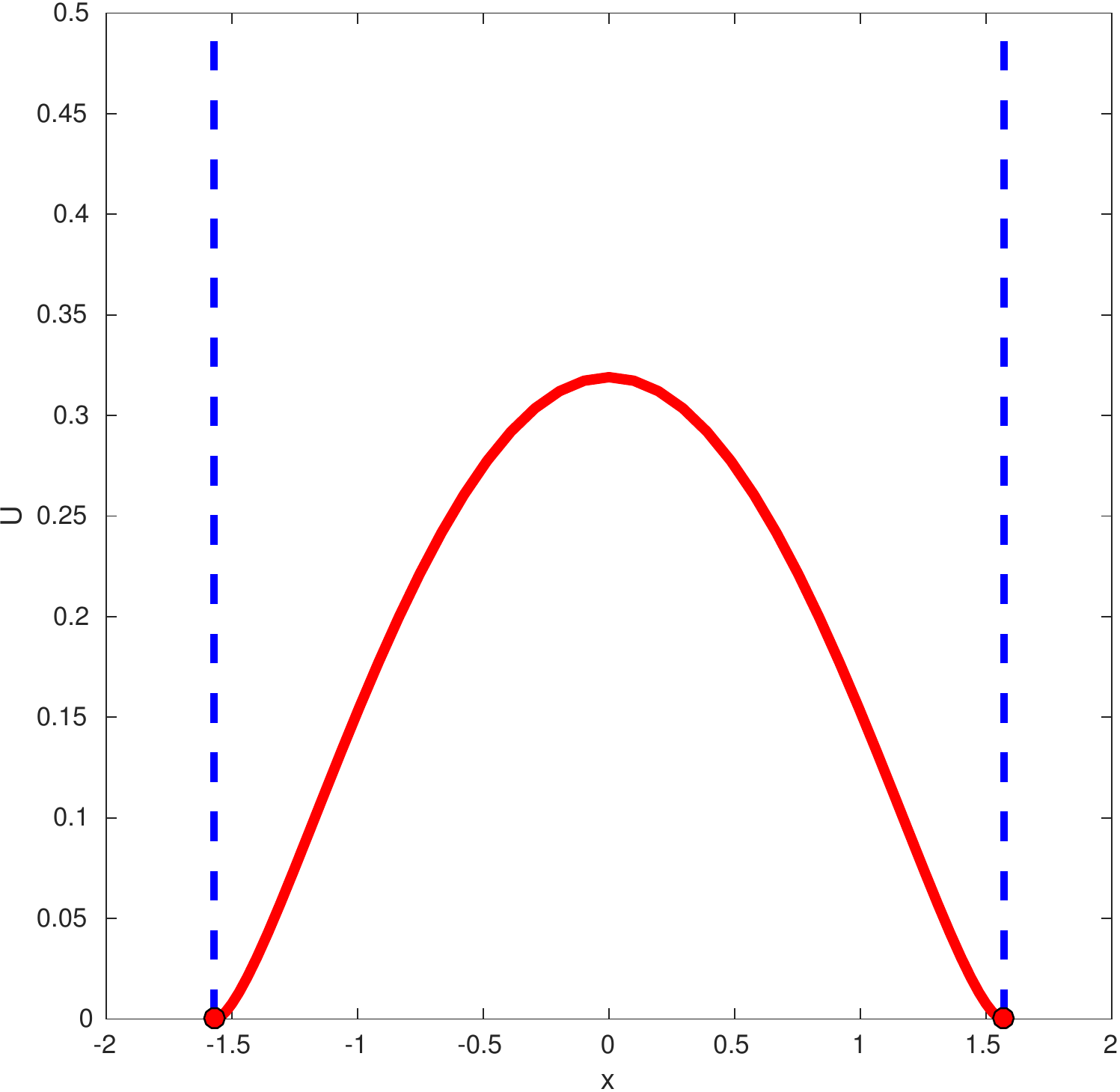}\caption{$t = 0.30$}\end{subfigure}\hspace{2mm}%
        \begin{subfigure}[b]{0.23\linewidth}\includegraphics[width=\textwidth]{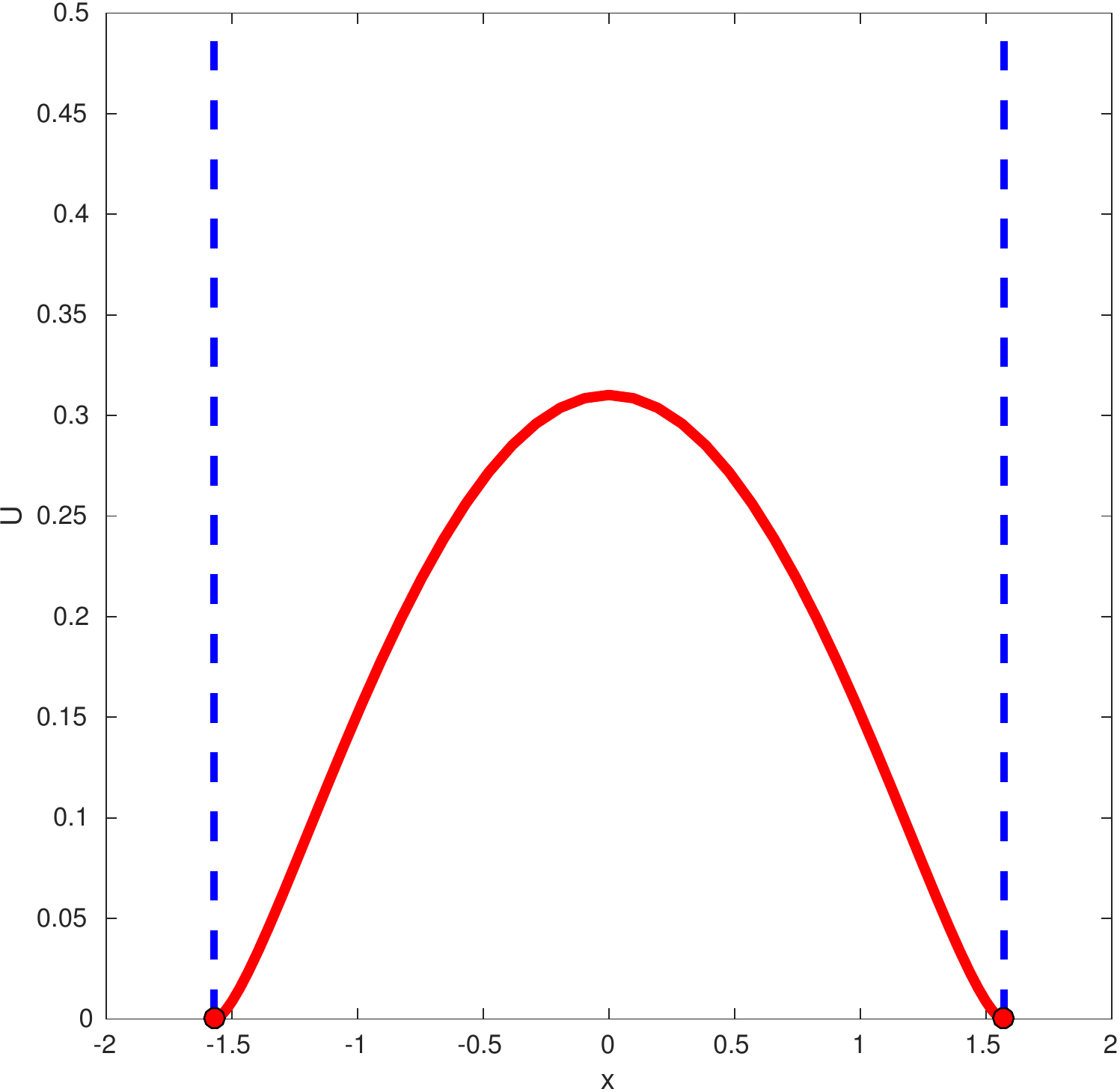}\caption{$t = 0.32$}\end{subfigure}\hspace{2mm}%
        \begin{subfigure}[b]{0.23\linewidth}\includegraphics[width=\textwidth]{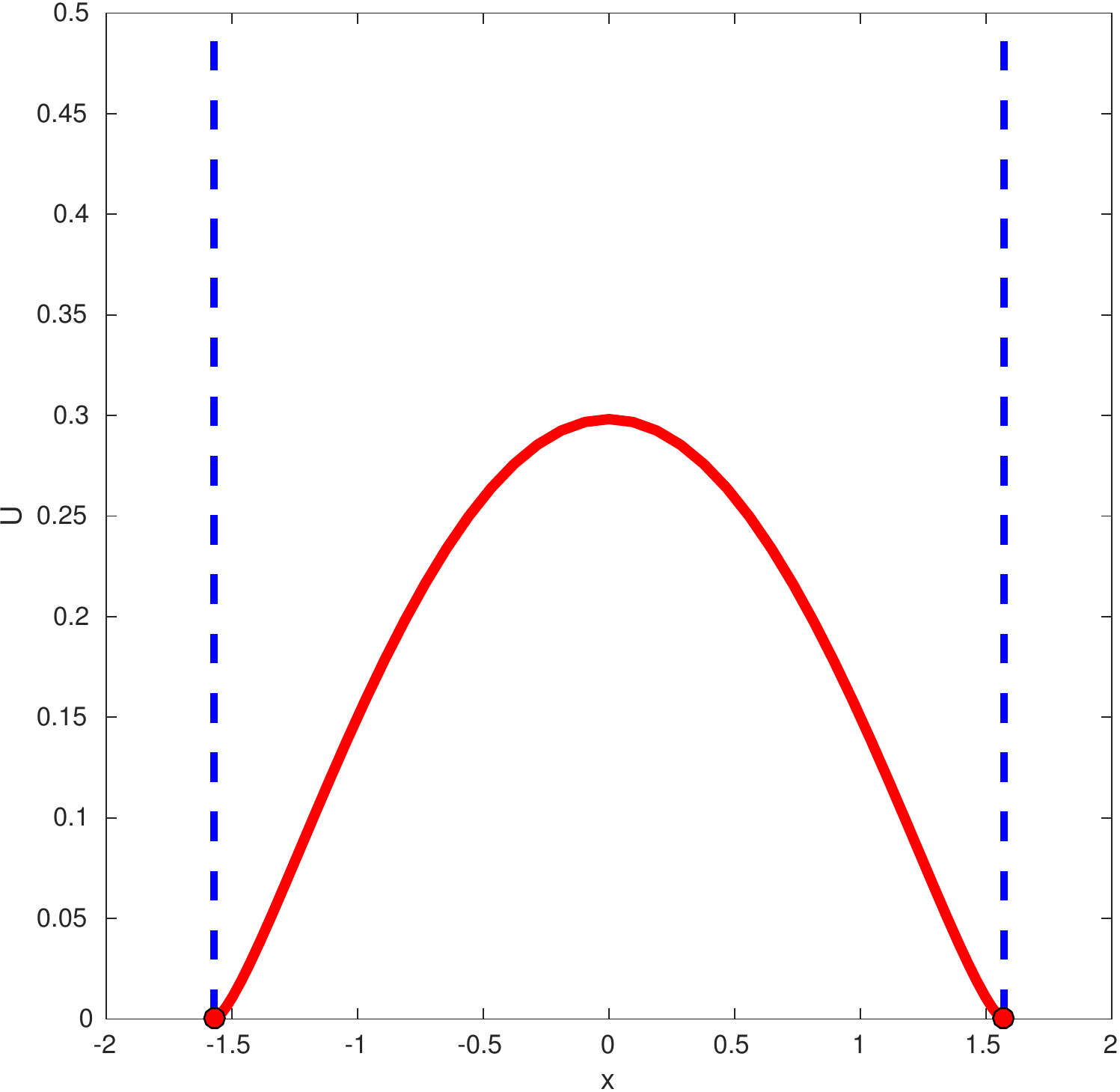}\caption{$t = 0.35$}\end{subfigure}\hspace{2mm}%
        \begin{subfigure}[b]{0.23\linewidth}\includegraphics[width=\textwidth]{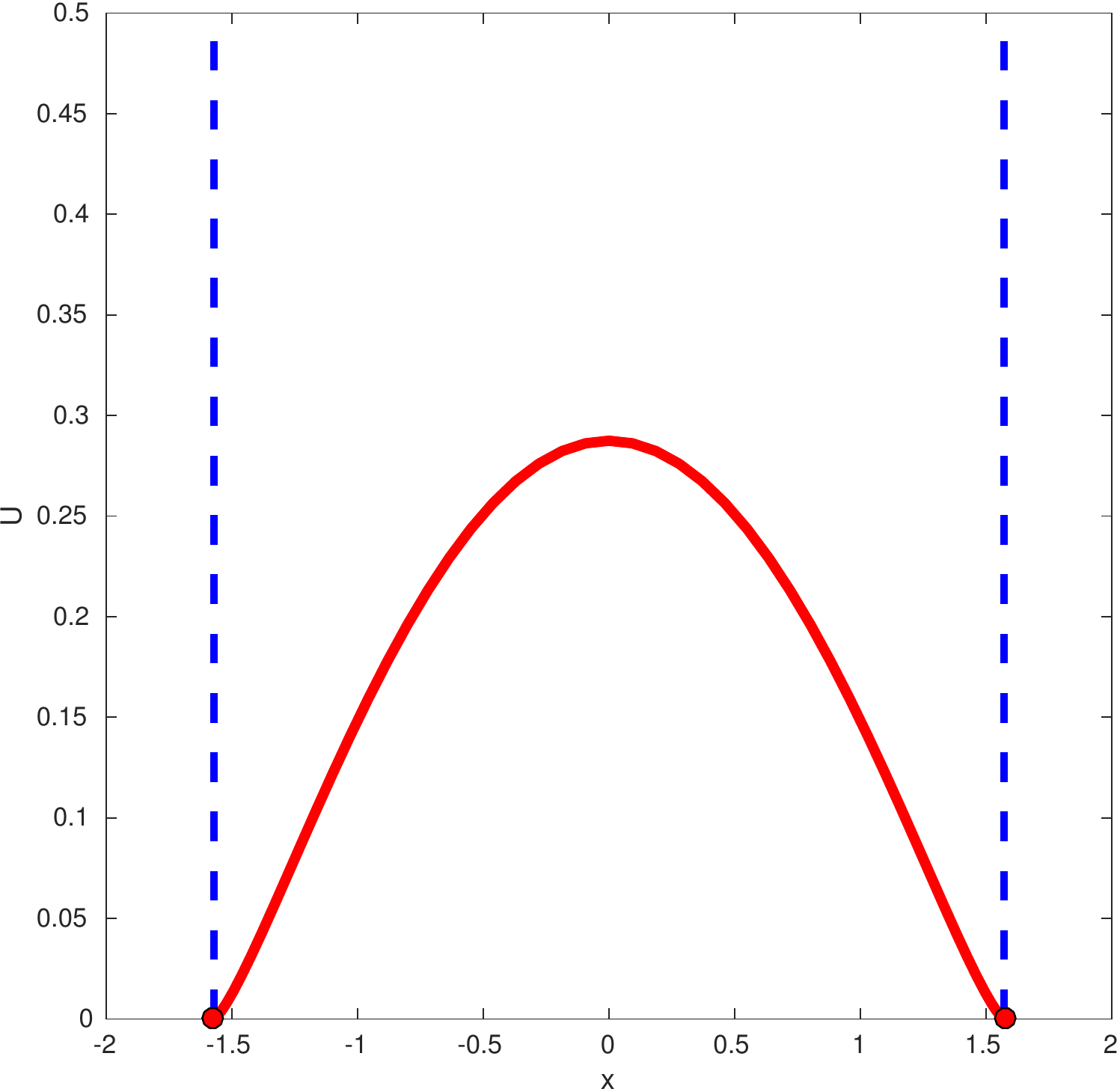}\caption{$t = 0.38$}\end{subfigure}\\%
        \begin{subfigure}[b]{0.23\linewidth}\includegraphics[width=\textwidth]{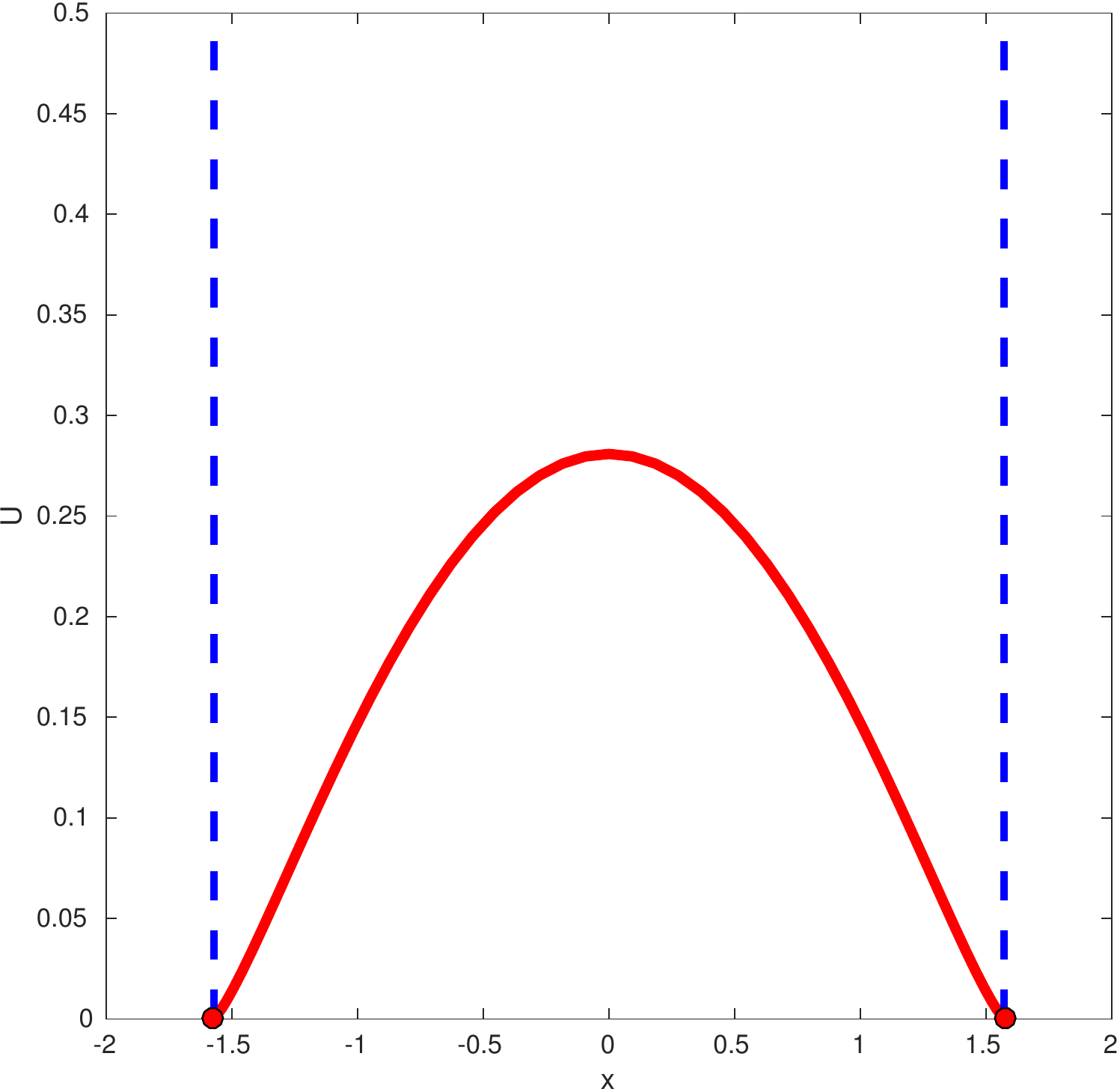}\caption{$t = 0.40$}\end{subfigure}\hspace{2mm}%
        \begin{subfigure}[b]{0.23\linewidth}\includegraphics[width=\textwidth]{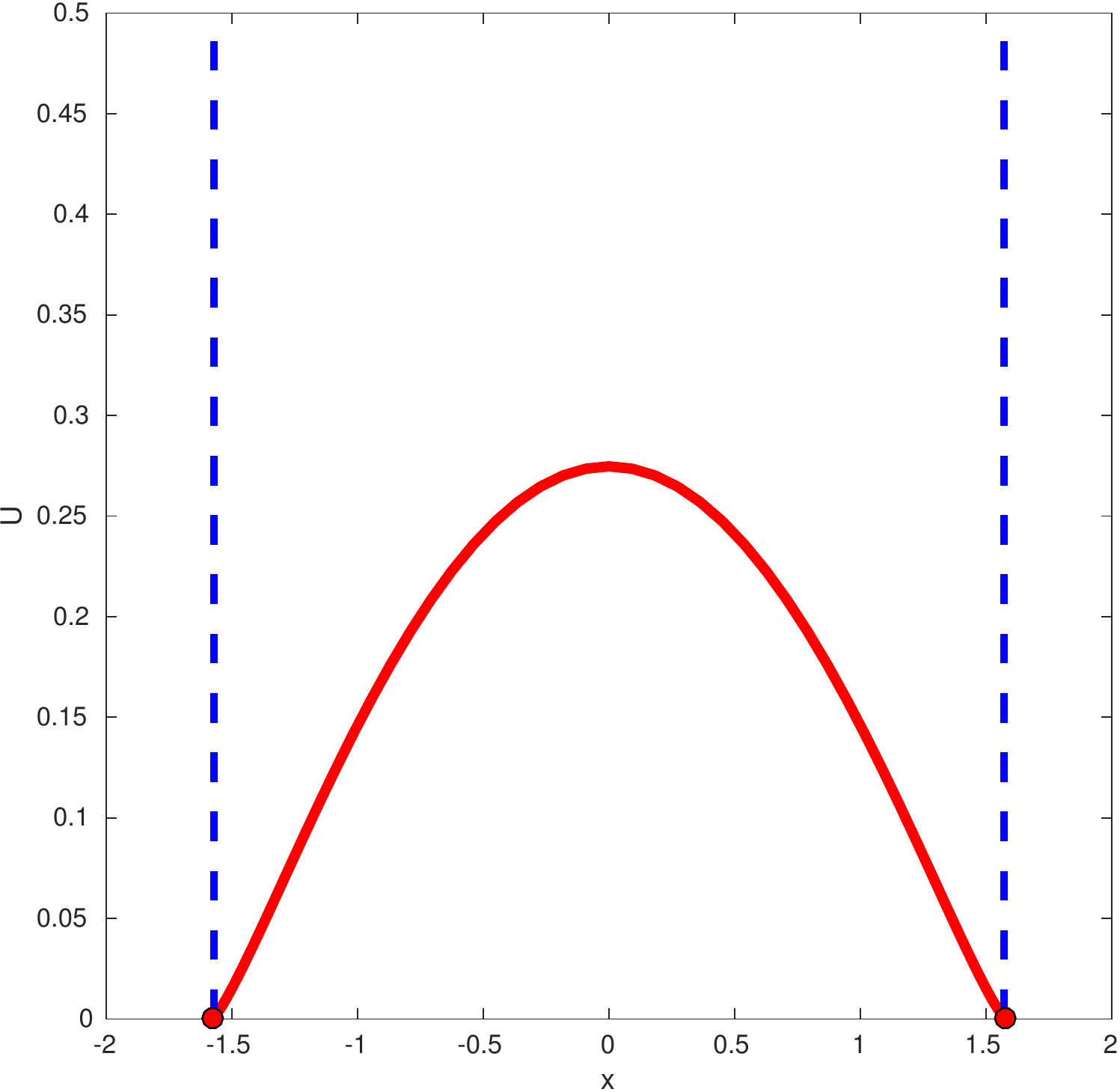}\caption{$t = 0.42$}\end{subfigure}\hspace{2mm}%
        \begin{subfigure}[b]{0.23\linewidth}\includegraphics[width=\textwidth]{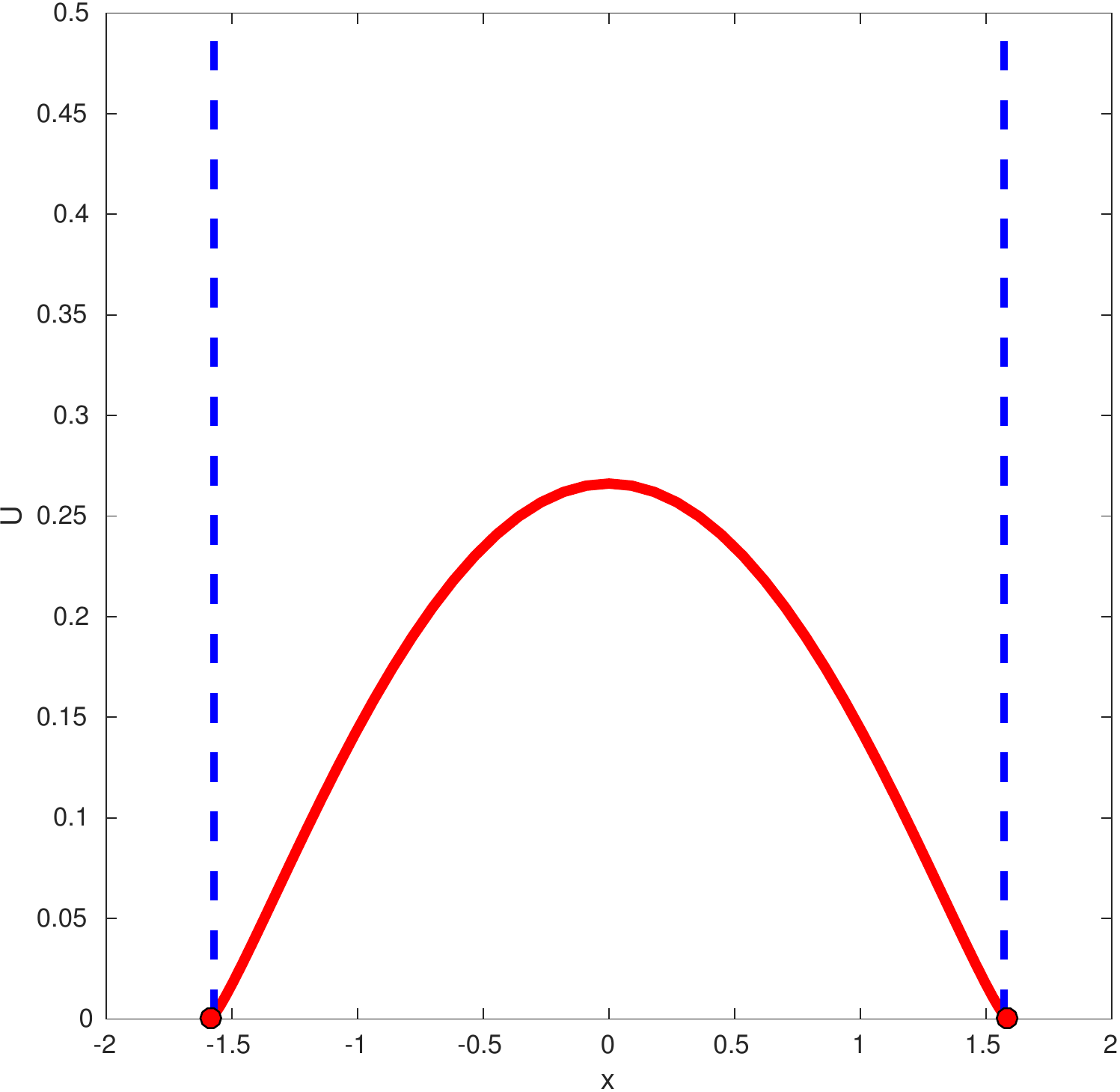}\caption{$t = 0.45$}\end{subfigure}\hspace{2mm}%
        \begin{subfigure}[b]{0.23\linewidth}\includegraphics[width=\textwidth]{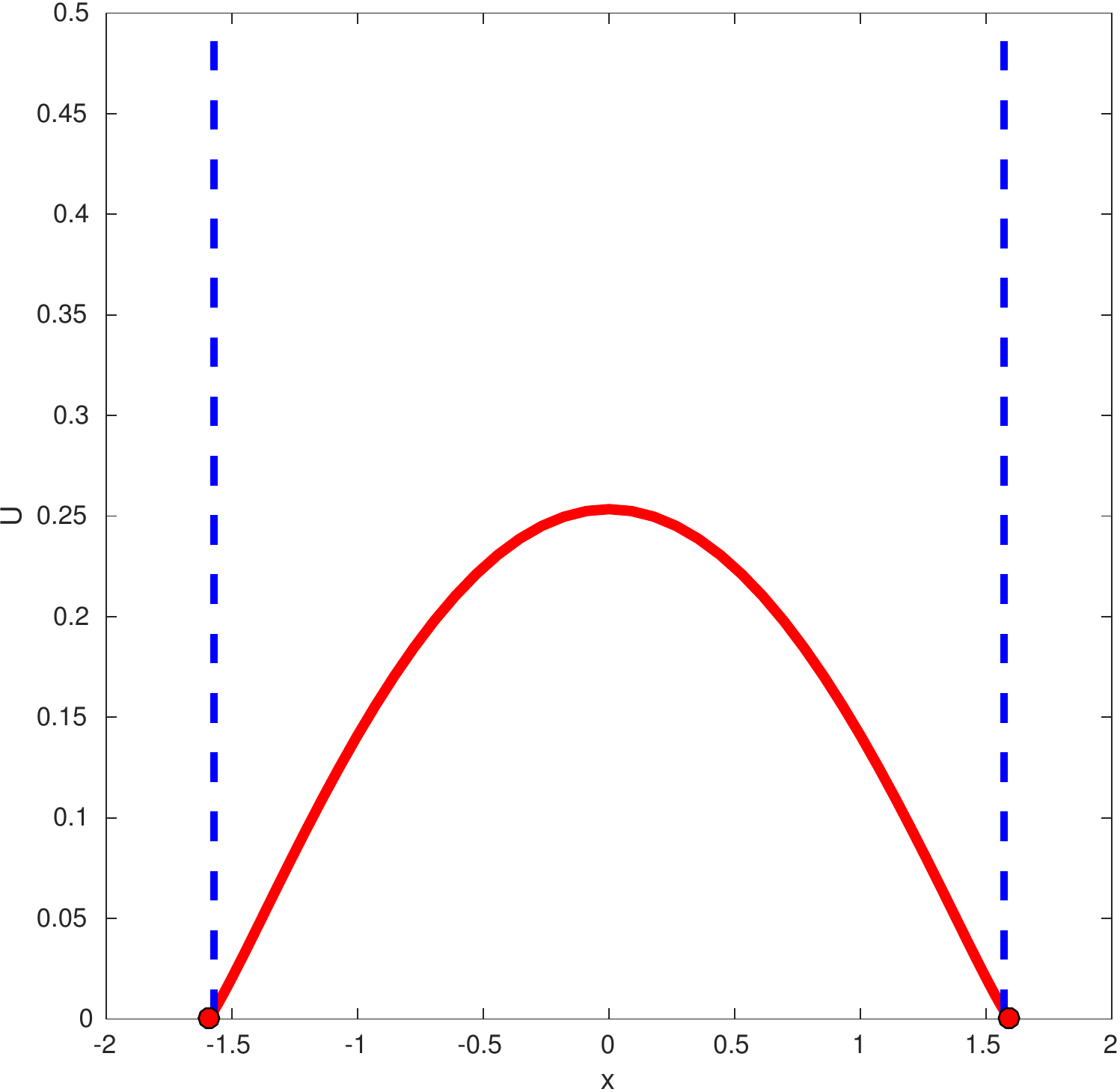}\caption{$t = 0.50$}\end{subfigure}\\%
        \begin{subfigure}[b]{0.23\linewidth}\includegraphics[width=\textwidth]{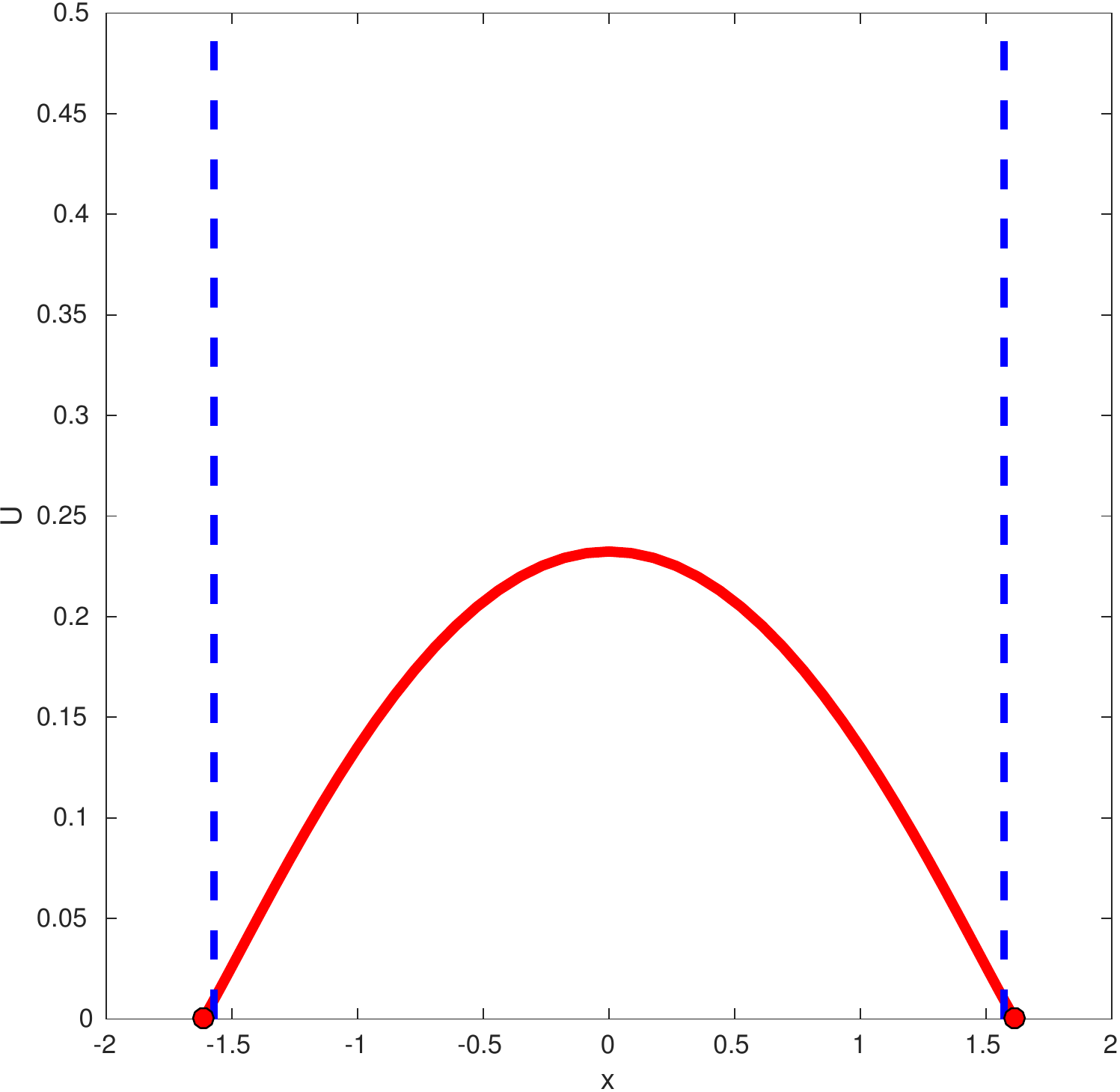}\caption{$t = 0.60$}\end{subfigure}\hspace{2mm}%
        \begin{subfigure}[b]{0.23\linewidth}\includegraphics[width=\textwidth]{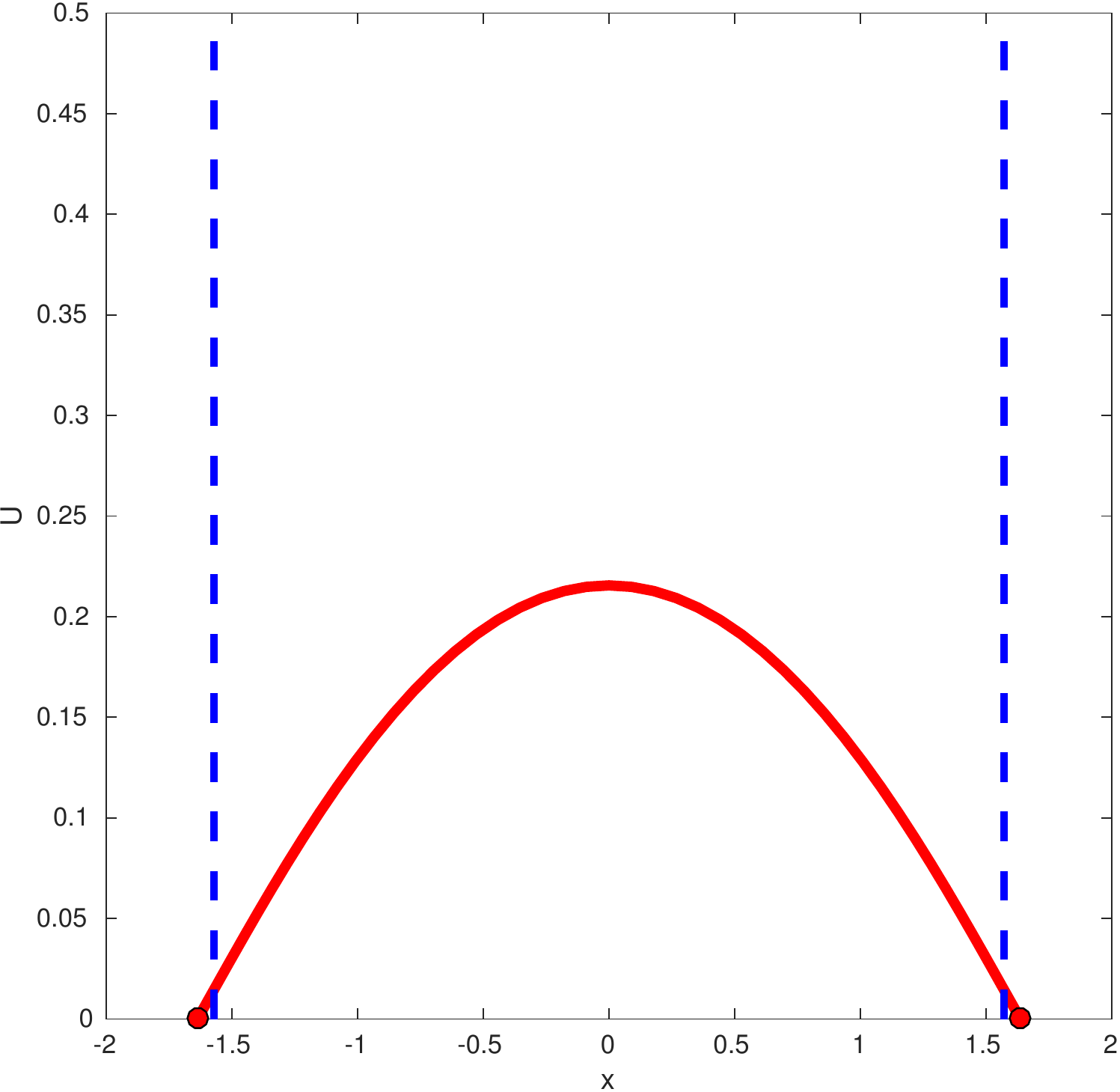}\caption{$t = 0.70$}\end{subfigure}\hspace{2mm}%
        \begin{subfigure}[b]{0.23\linewidth}\includegraphics[width=\textwidth]{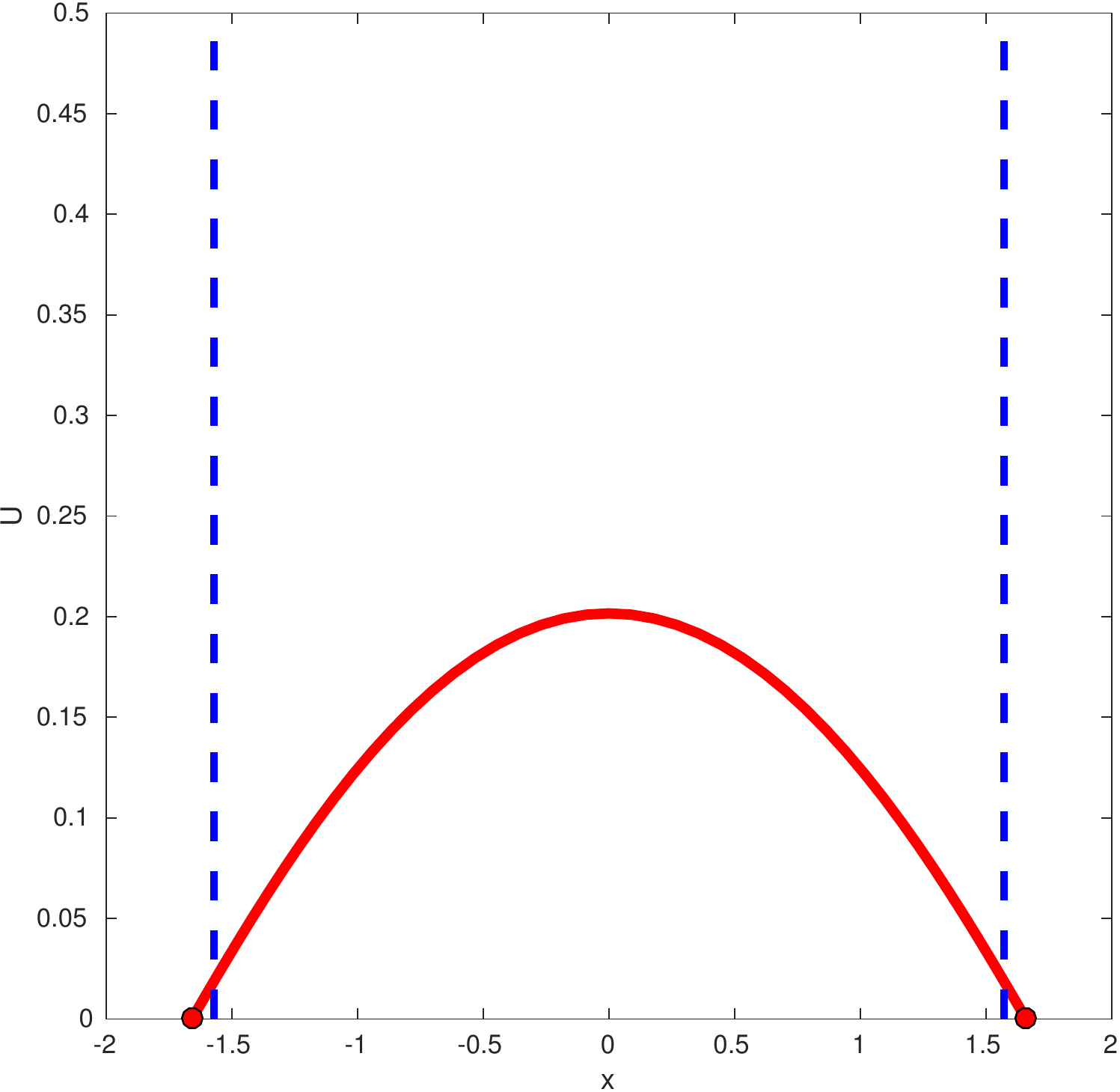}\caption{$t = 0.80$}\end{subfigure}\hspace{2mm}%
        \begin{subfigure}[b]{0.23\linewidth}\includegraphics[width=\textwidth]{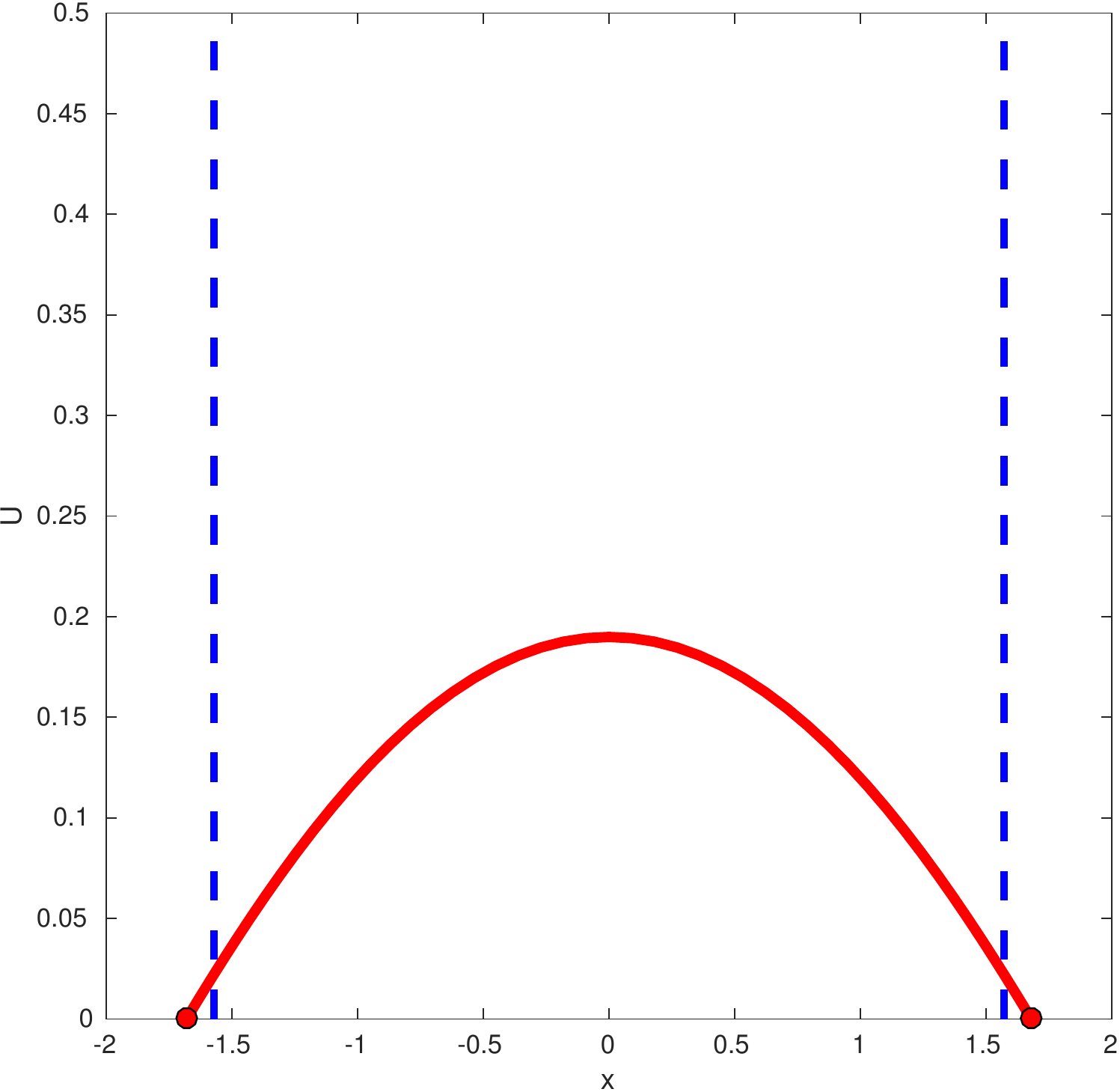}\caption{$t = 0.90$}\end{subfigure}%
        \caption{Example~\ref{ex:pme-V-waiting}. The cross section at $y = 0$ of a computed solution is shown at various time instants ($N = 4011$).}
        \label{fig:pme-V-waiting-2D-cross-section}
    \end{figure}

\begin{exam}[Waiting time phenomenon] \label{ex:pme-V-waiting}
It is known that for some initial solutions, the IBVP \eqref{eqn:PME-V-Darcy} exhibits the waiting-time
phenomenon where the free boundary does not move initially until a finite amount of time has elapsed.
Consider an example with $m = 2$ and the initial solution
    \begin{equation}
        v_0(x,y) = 
        \begin{cases}
            \frac{1}{m} \cos^m\left(\sqrt{x^2+y^2}\right) , & \quad \text{for} \quad \sqrt{x^2 + y^2} \leq \frac{\pi}{2} \\
            0 , & \quad \text{otherwise} .
        \end{cases}
    \label{eqn:pme-V-waiting}
    \end{equation}
We observe that 
 \[
     \grad \left(\frac{1}{m}\cos^m\left(\sqrt{x^2+y^2}\right)\right) = -\frac{ \cos^{m-1}(\sqrt{x^2+y^2})\, \sin(\sqrt{x^2+y^2})}{\sqrt{x^2+y^2}} 
 \begin{bmatrix} x \\ y \end{bmatrix} ,
 \]
which vanishes at $\sqrt{x^2+y^2} = \frac{\pi}{2}$. According to Darcy's law \eqref{eqn:v-darcy}, the velocity
of the free boundary is zero, and thus we should not expect the free boundary to move initially.
In Fig. \ref{fig:pme-V-waiting}, the mesh and associated solution are plotted at several instants
(with contrasting circles representing the initial boundary). To further confirm the waiting-time phenomenon,
we also plot the cross sections of the solution at various time instants
in Fig. \ref{fig:pme-V-waiting-2D-cross-section}.
A closer look suggests that the free boundary does not start moving until around $t = 0.2$.
Moreover, no oscillations are visible in the computed solutions. This is in contrast with embedding methods
which typically produce computed solutions with small oscillations near the free boundary;
e.g., see \cite[Fig.~12]{NH2016}.
\qed
\end{exam}

\begin{figure}[ht]
        \centering
        \vspace{-10mm}
        \begin{subfigure}[b]{0.35\linewidth}\includegraphics[scale=0.28]{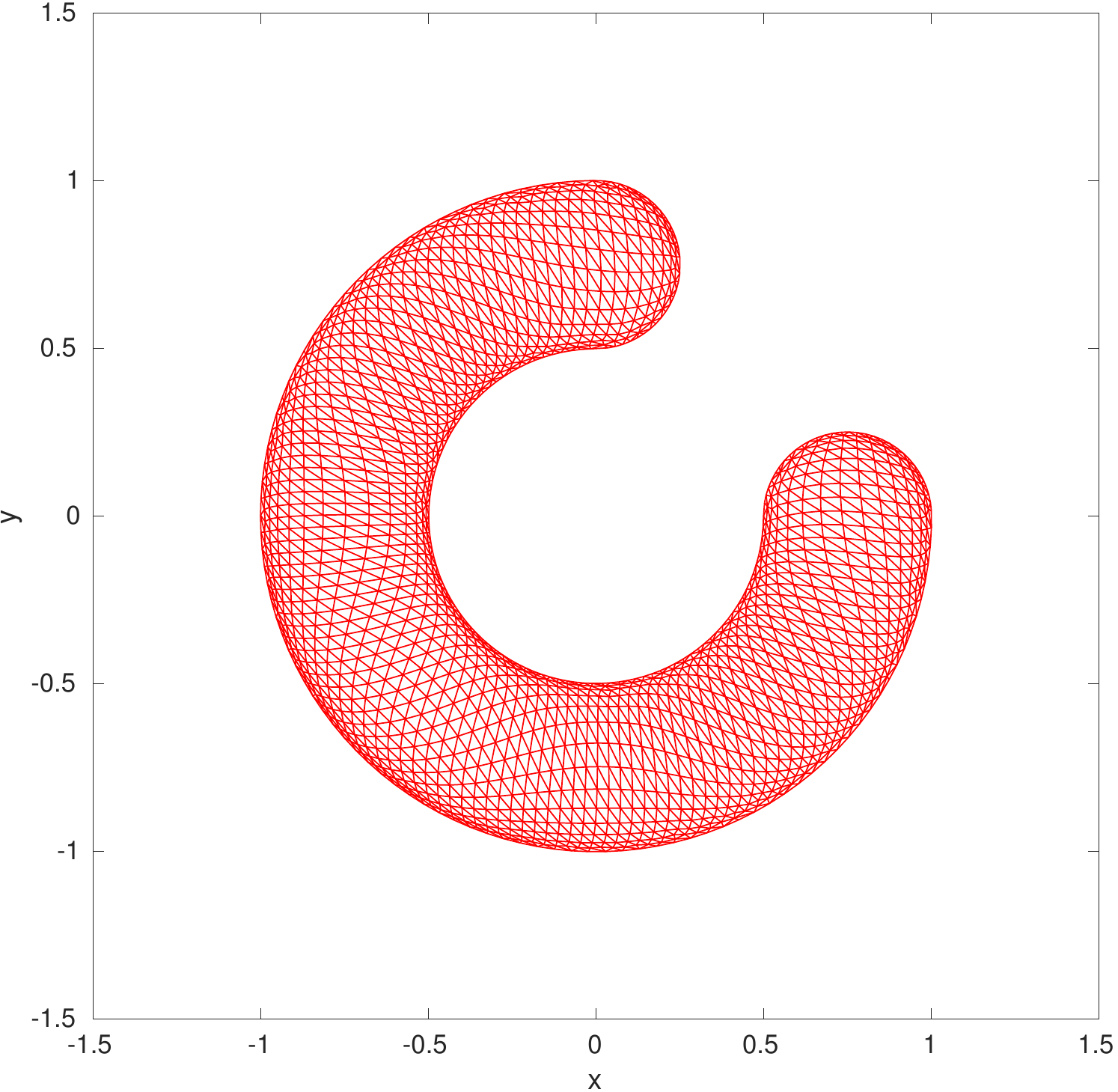}\caption{$t = 0.1$}\end{subfigure}\hspace{5mm}
        \begin{subfigure}[b]{0.35\linewidth}\includegraphics[scale=0.26]{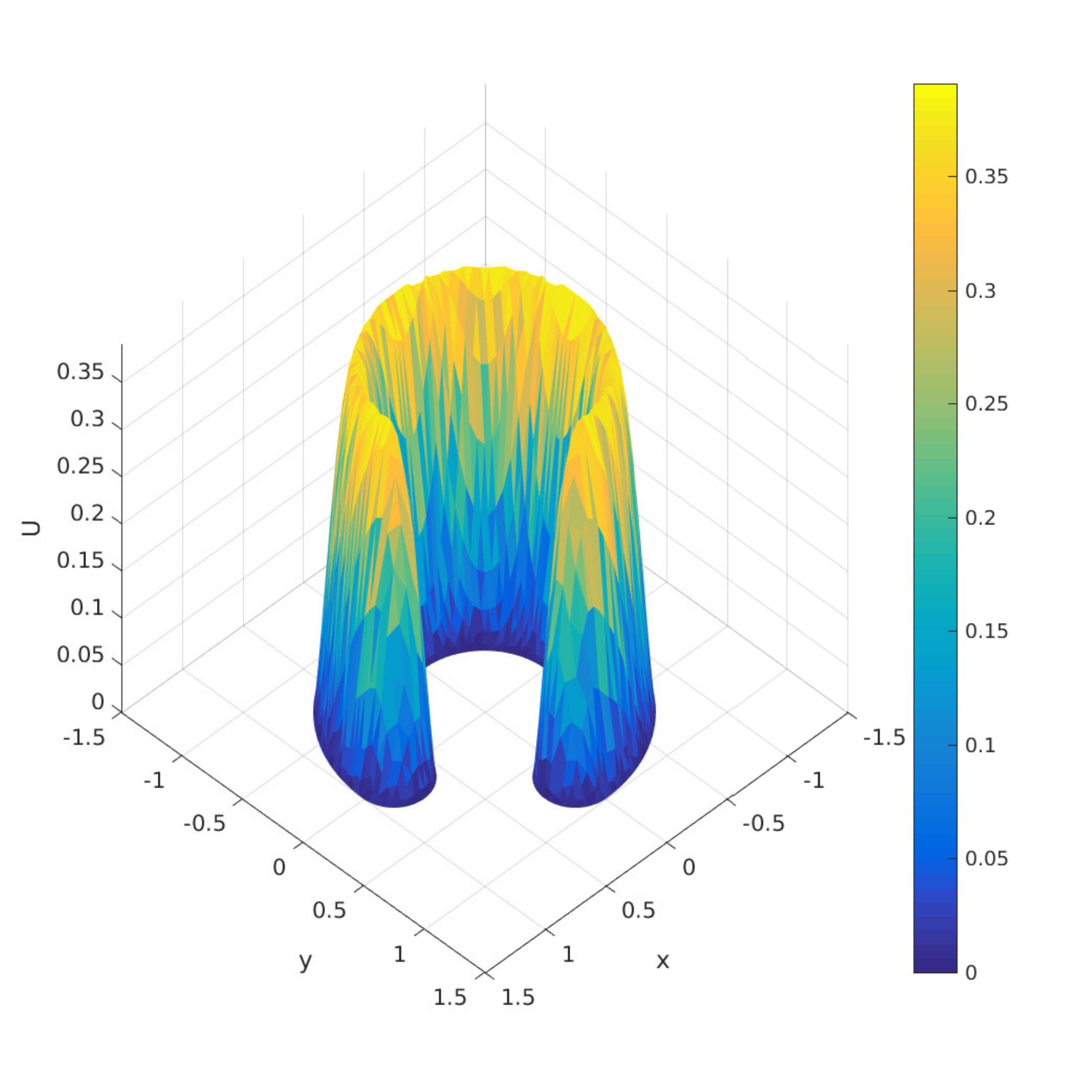}\caption{$t = 0.1$}\end{subfigure}\\
        \begin{subfigure}[b]{0.35\linewidth}\includegraphics[scale=0.28]{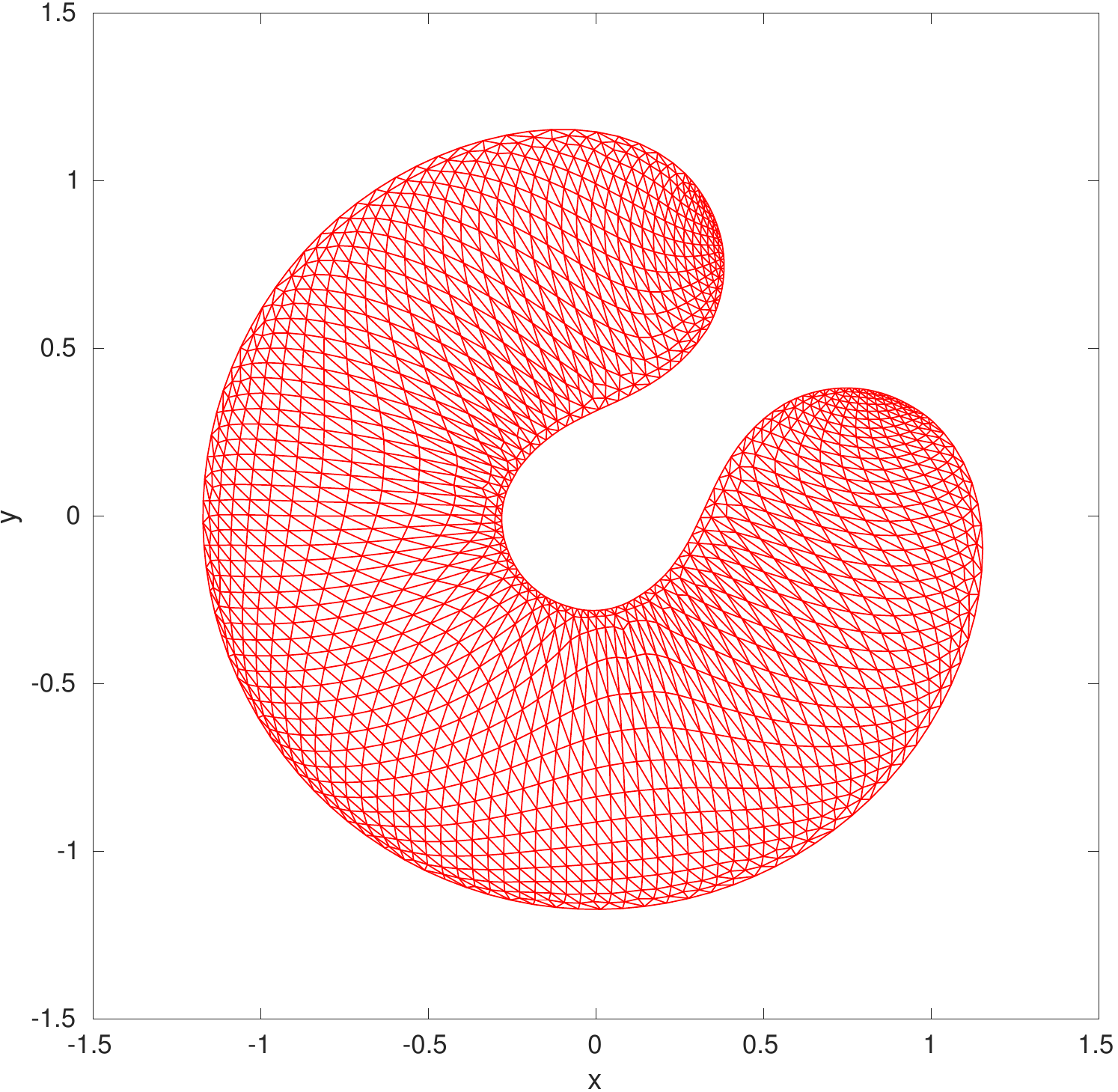}\caption{$t = 0.3$}\end{subfigure}\hspace{5mm}
        \begin{subfigure}[b]{0.35\linewidth}\includegraphics[scale=0.26]{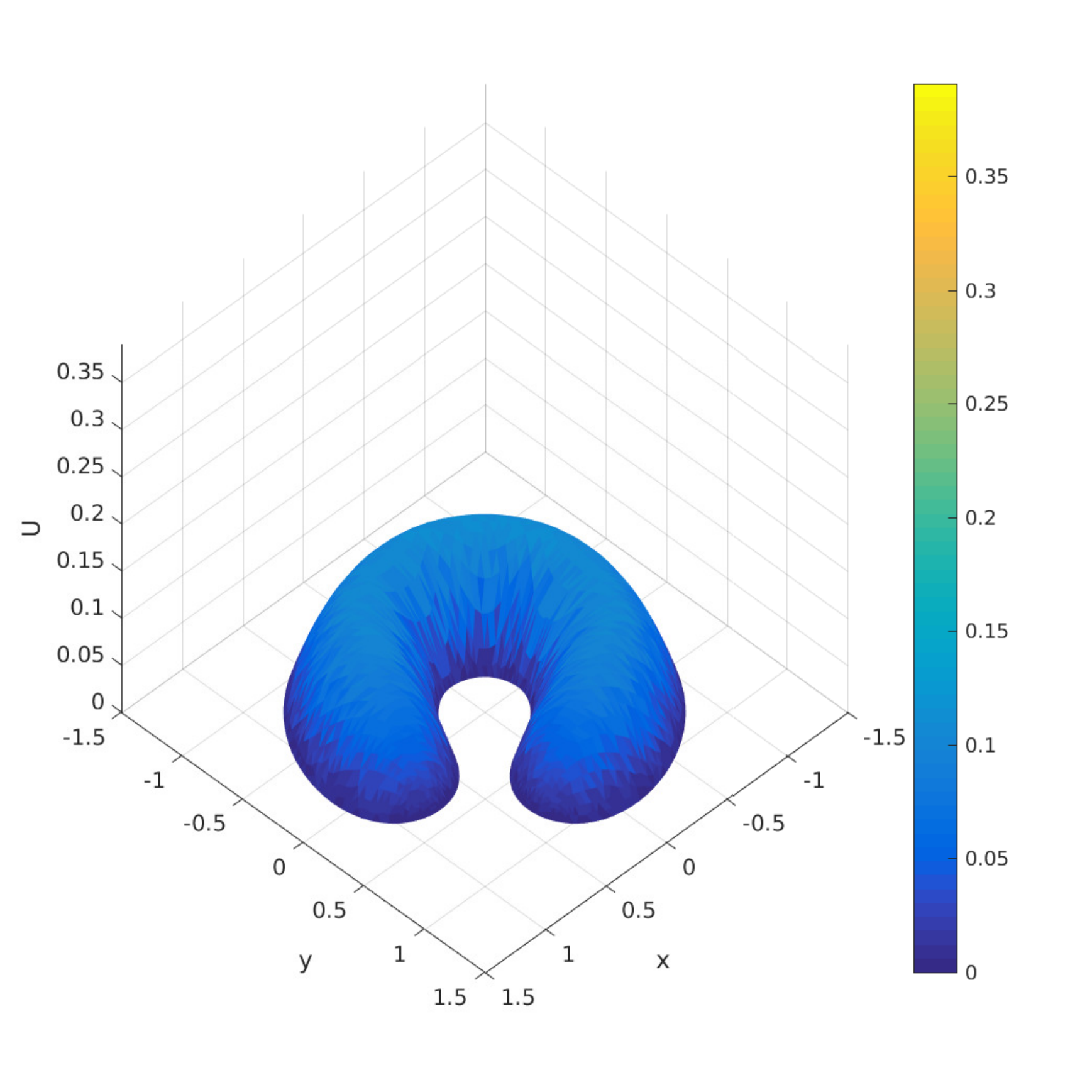}\caption{$t = 0.3$}\end{subfigure}\\
        \begin{subfigure}[b]{0.35\linewidth}\includegraphics[scale=0.28]{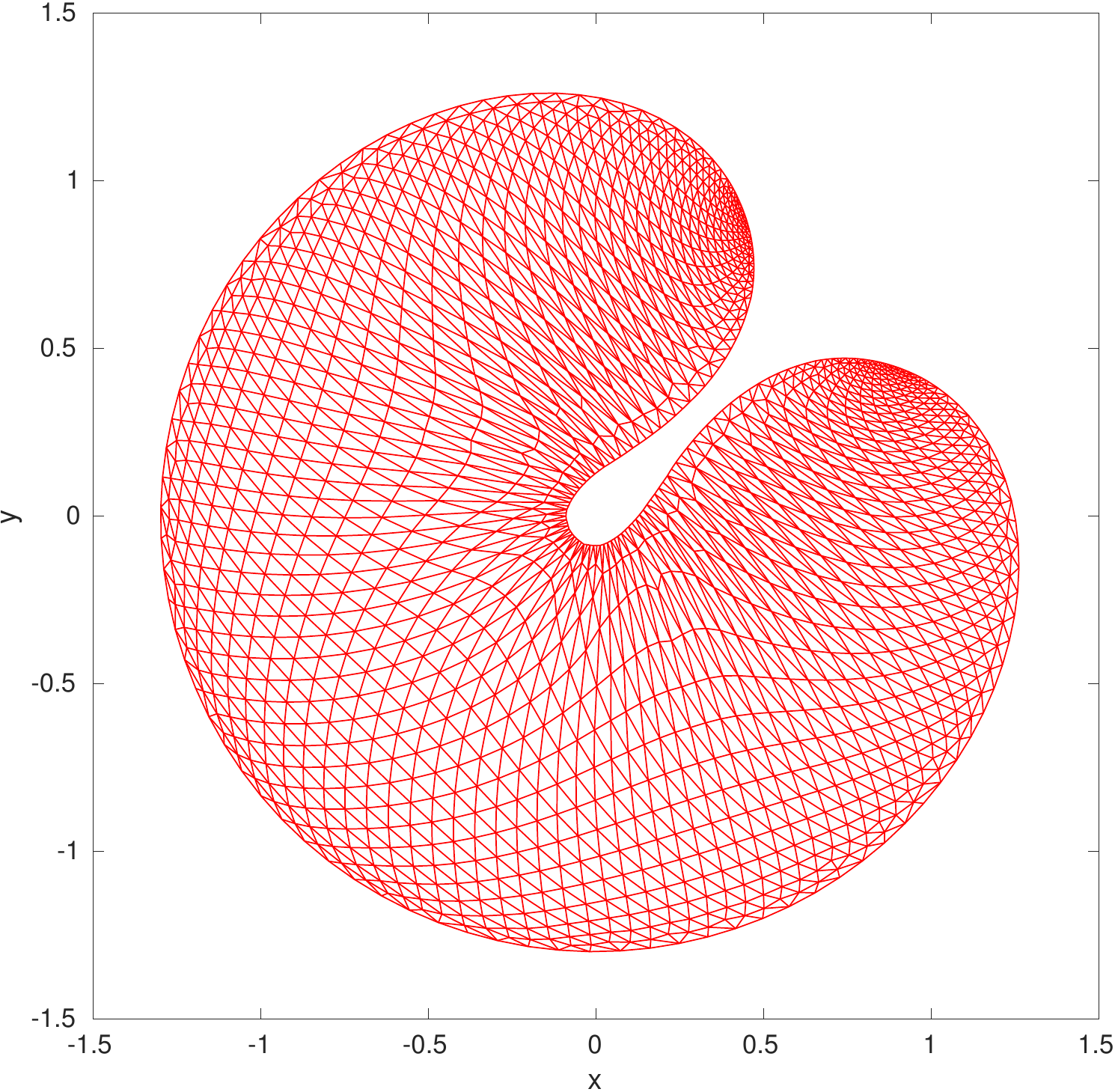}\caption{$t = 0.75$}\end{subfigure}\hspace{5mm}
        \begin{subfigure}[b]{0.35\linewidth}\includegraphics[scale=0.26]{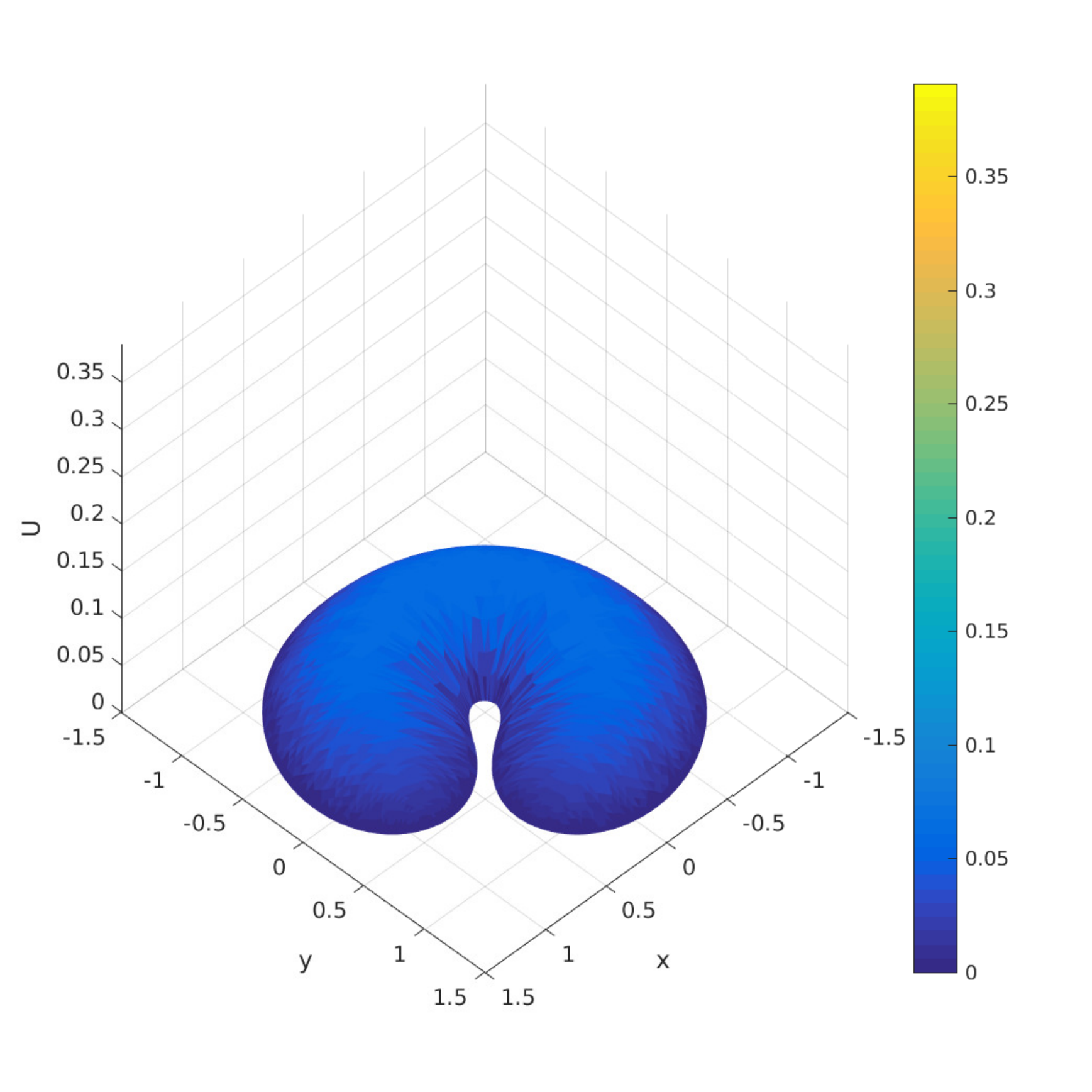}\caption{$t = 0.75$}\end{subfigure}\\
        \begin{subfigure}[b]{0.35\linewidth}\includegraphics[scale=0.28]{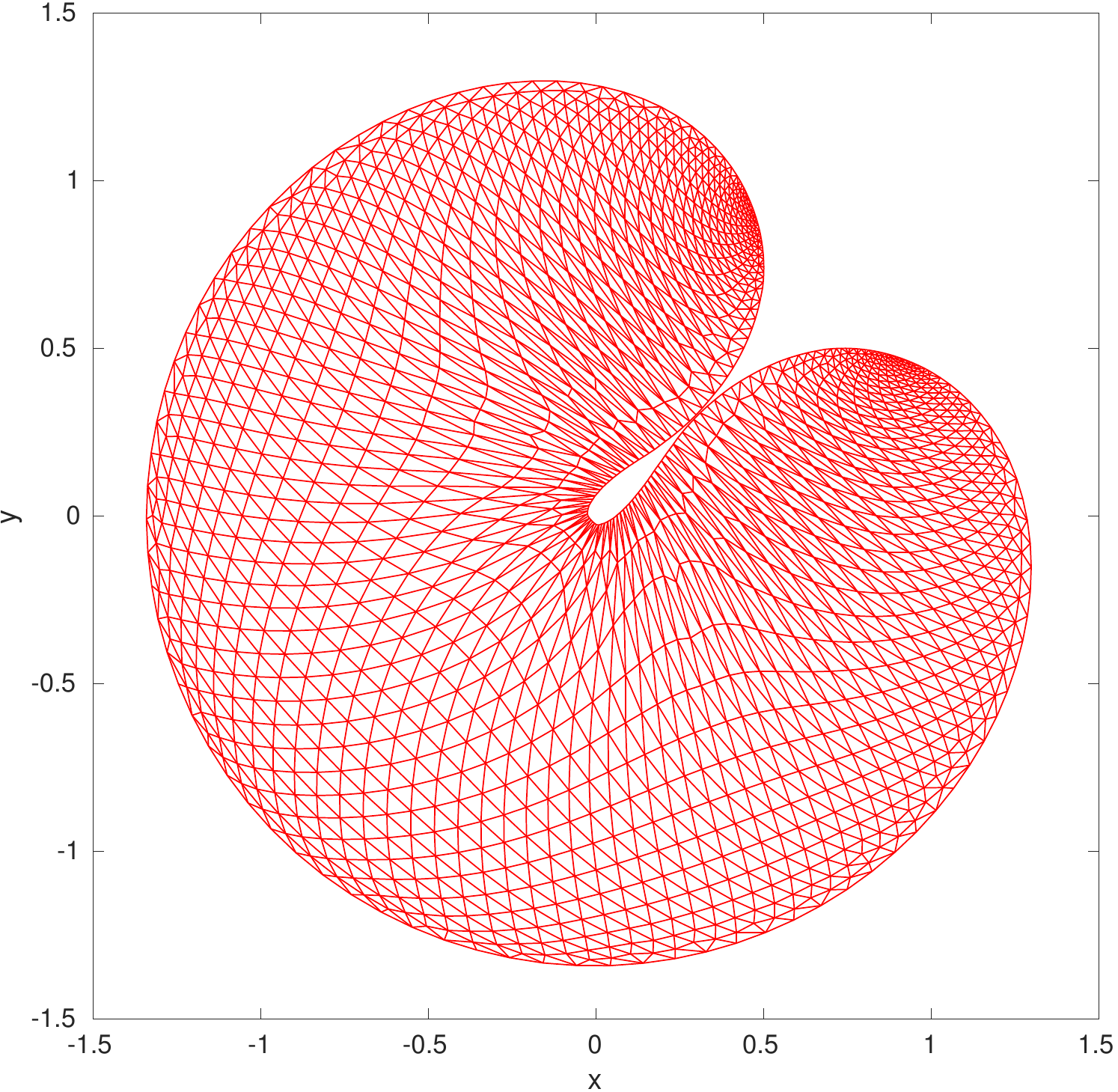}\caption{$t = 1$}\end{subfigure}\hspace{5mm}
        \begin{subfigure}[b]{0.35\linewidth}\includegraphics[scale=0.26]{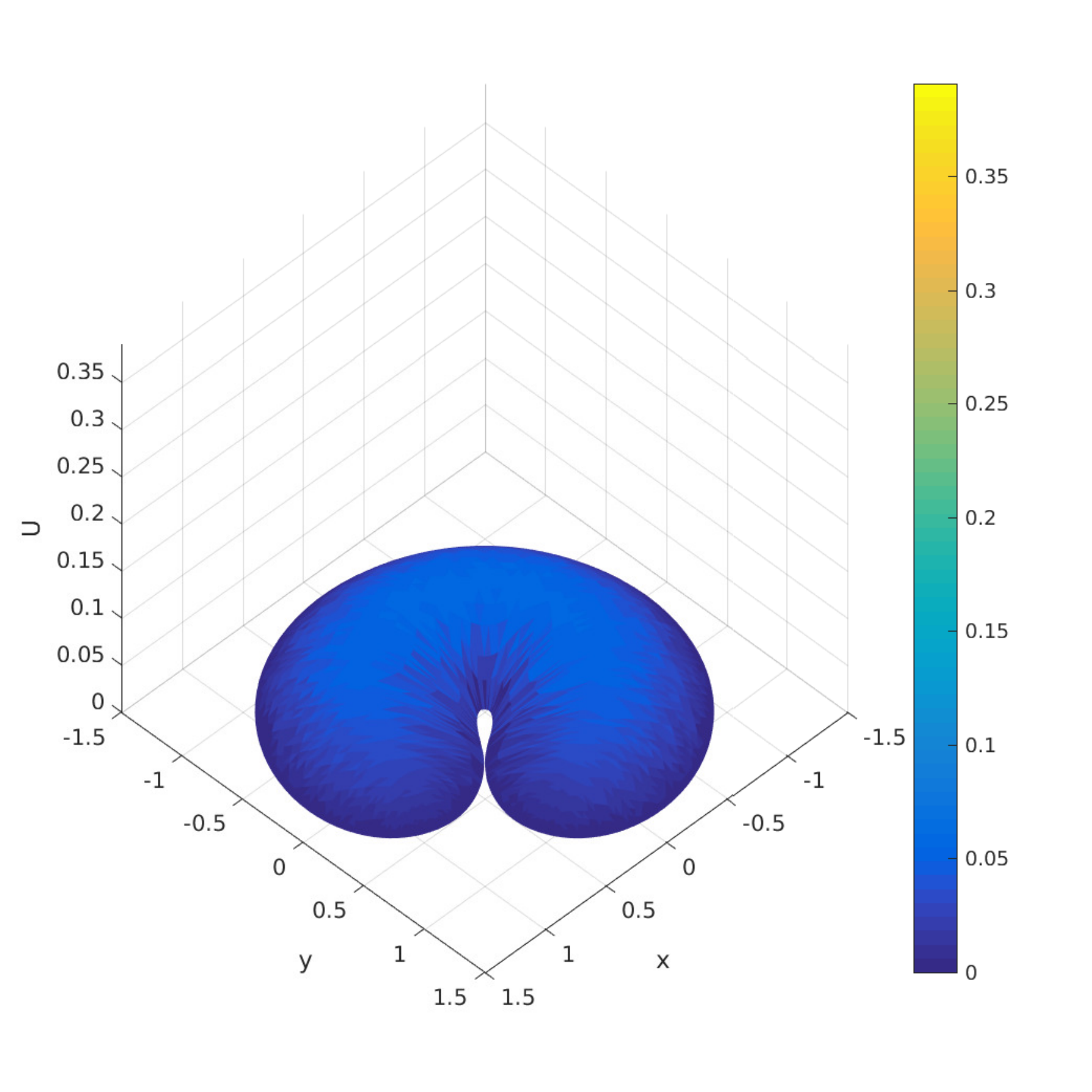}\caption{$t = 1$}\end{subfigure}%
        \caption{Example~\ref{ex:pme-V-complex}. A computed solution is shown at various time instants ($N = 4011$).}
        \label{fig:pme-V-complex}
    \end{figure}

\begin{exam}[Complex domain] \label{ex:pme-V-complex}
  In this example we apply our method to PME-V with $m = 2$ and 
    \begin{equation}
        v_0(x,y) = 
        \begin{cases}
            25 \left[0.25^2-(\sqrt{x^2+y^2}-0.75)^2\right]^{\frac{3}{2}} , \hspace{-2mm}& \quad \text{for} \quad \sqrt{x^2 + y^2} \in [0.5,1] \text{ and } (x < 0 \text{ or } y < 0) \\
            25 \left[0.25^2 - x^2 - (y-0.75)^2\right]^{\frac{3}{2}} , & \quad \text{for} \quad x^2 + (y-0.75)^2 \leq 0.25^2 \text{ and } x \geq 0  \\
            25 \left[0.25^2 - (x-0.75)^2 - y^2\right]^{\frac{3}{2}} , & \quad \text{for} \quad (x-0.75)^2 + y^2 \leq 0.25^2 \text{ and } y \geq 0  \\
            0 , & \quad \text{otherwise} .
        \end{cases}
    \label{eqn:pme-V-complex}
    \end{equation}
The partial donut-shaped support pertaining to this initial solution is first seen in Baines~et~al.~\cite{BHJ05a}.
Fig. \ref{fig:pme-V-complex} shows a representative mesh and the corresponding solution.
One can see that the mesh evolves smoothly following the evolving free boundary.
The simulation does not encounter any difficulty until the hole closes out and the boundary begins to touch.
This example also demonstrates that the MMPDE moving mesh method works well for concave domains
(without causing any mesh tangling or crossing).
\qed
\end{exam}
    
\section{Conclusions}
\label{SEC:conclusion}

In the previous sections we have studied the adaptive finite element solution of PME
in the pressure formulation (PME-V) \eqref{eqn:PME-V}. The motivation is that the pressure
variable $v = u^m/m$ has better regularity than $u$ and the equation for the movement
of the free boundary can be expressed naturally in $v$.
In fact, many theoretical works (e.g., \cite{Angenent1988,Caffarelli1980,Caffarelli1987,Daskalopoulos1998a})
have relied on PME-V to obtain mathematical properties for $u$.
But PME-V has not been used in numerical studies of PME in the past.
We have used PME-V in a nonembedding approach with which PME is discretized
only within the support of the solution and the free boundary is traced explicitly via
Darcy's law \eqref{eqn:v-darcy}.
This approach, though not as robust as the embedding methods
(e.g. \cite{Duque2013,NH2016,Socolovsky1984,Zhang2009})
in dealing with problems having complex free boundaries, offers better solution regularity,
leads to more accurate location of the free boundary, and avoids unwanted oscillations in computed solutions.
We have employed the boundary-mesh-physical splitting solution procedure where
the boundary equation (Darcy's law) is first solved using the Euler scheme,
then the mesh is updated, and finally PME is solved using piecewise line finite elements on the moving mesh.
The MMPDE moving mesh method has been used for mesh updating. Since $v$ is smooth near
the free boundary in the current situation, a $v$-value based metric tensor (\ref{eqn:M}) has been
used to control the mesh concentration near the boundary.

Numerical results obtained with the moving mesh finite element method have been presented
for three two-dimensional examples. They show that the method
can effectively concentrate mesh points near the free boundary without causing any mesh tangling
or crossing for both convex and concave domains. Moreover, no oscillations have been observed
in the computed solutions. Furthermore, the error in $v$ shows a convergence rate of
second order in space and first order in time. Finally, the error in the location of the free boundary
exhibits almost second-order convergence in the $L^\infty$ norm.
Interestingly, these convergence behaviors are essentially independent of $m$ for its tested
range from $2$ to $15$.
It is noted that the original variable $u$ can be obtained from the computed solution through
the definition $u = (m v)^{1/m}$.
However, the convergence in the so-computed $u$ is not satisfactory. As a matter of fact,
the convergence order is between 0.5th-order and first-order and closer to 0.5th-order for large $m$
in the $L^2$ norm and is between first-order and second-order and close to first-order for large $m$ 
in the $L^1$ norm.

Recall that the embedding moving mesh finite element method developed in \cite{NH2016} based on
the original formulation (PME-U) is able to handle problems with complex, emerging,
or splitting free boundaries and is second-order in space in $u$ for relatively small $m$.
For large $m$, it can run into difficulties since the mesh elements have to be stretched
extremely thin near the free boundary which can cause the time integration of PME to stop
due to too small time steps.
Comparing these with the observations made in this work,
we can conclude that PME-V can offer advantages over PME-U for situations with large $m$
or when a more accurate location of the free boundary is desired. For other situations, numerical
methods based on the original formulation can be more advantageous.



\end{document}